\documentclass{amsart}

\usepackage{amsmath,amsthm}
\usepackage {latexsym}
\usepackage{amssymb}

\newcommand{\less}{\lesssim}

\newcommand{\beal}{\begin{align}}
\newcommand{\enal}{\end{align}}
\newcommand{\bealn}{\begin{align*}}
\newcommand{\enaln}{\end{align*}}
\newcommand{\bear}{\begin{eqnarray}}
\newcommand{\eear}{\end{eqnarray}}
\newcommand{\beeq}{\begin{equation}}
\newcommand{\eneq}{\end{equation}}

\newcommand{\const}{\mbox{\rm const}}

\newcommand{\supp}{\mbox{\rm supp}}
\newcommand{\spec}{{\rm spec}}

\newcommand{\eps}{{\varepsilon}}
\newcommand{\R}{{\mathbb R}}
\newcommand{\C}{{\mathbb C}}

\newcommand{\Compl}{{\mathbb C}}

\newcommand{\calS}{{\mathcal S}}

\newcommand{\calH}{{\mathcal H}}

\newcommand{\kato}{{\mathcal K}}

\newcommand{\la}{\langle}
\newcommand{\ra}{\rangle}

\def\pr{\partial}
\def\Lap{\Delta}
\def\nn{\nonumber}
\def\bm{\left[ \begin{array}{cc}}
\def\endm{\end{array}\right]}
\def\ker{{\rm ker}}
\def\Ran{{\rm Ran}}
\def\Hil{{\mathcal H}}

\def\f{\frac}

\newtheorem{theorem}{Theorem}
\newtheorem{lemma}[theorem]{Lemma}
\newtheorem{defi}[theorem]{Definition}
\newtheorem{cor}[theorem]{Corollary}
\newtheorem{prop}[theorem]{Proposition}
\newtheorem{proposition}[theorem]{Proposition}

\theoremstyle{remark}
\newtheorem{remark}[theorem]{Remark}

\def\il{\int\limits}
\def\bin{\binom}
\def\til{\tilde}

\renewcommand{\Im}{\,{\rm Im}}

\def\Laplace{\Delta}
\def\half{\frac12}
\def\uone{u^{(1)}}
\def\aone{a^{(1)}}
\def\utwo{u^{(2)}}
\def\atwo{a^{(2)}}
\def\vone{v^{(1)}}
\def\bone{b^{(1)}}
\def\vtwo{v^{(2)}}
\def\btwo{b^{(2)}}
\def\hone{h^{(1)}}
\def\htwo{h^{(2)}}
\def\hj{h^{(j)}}

\def\uj{u^{(j)}}
\def\aj{a^{(j)}}
\def\vj{v^{(j)}}
\def\bj{b^{(j)}}
\def\pone{p^{(1)}}
\def\ptwo{p^{(2)}}
\def\qtwo{q^{(2)}}
\def\qone{q^{(1)}}
\def\tilvj{{\tilde v}^{(j)}}
\def\nj{n^{(j)}}
\def\wj{w^{(j)}}
\def\none{n^{(1)}}
\def\ntwo{n^{(2)}}
\def\wone{w^{(1)}}
\def\wtwo{w^{(2)}}
\def\Ea{E^{(a)}}
\def\Eb{E^{(b)}}
\def\ma{m^{(a)}}
\def\mb{m^{(b)}}
\def\Sa{S^{(a)}}
\def\Aa{A^{(a)}}
\def\Ab{A^{(b)}}
\def\Sa{S^{(a)}}
\def\Sb{S^{(b)}}
\def\calE{{\mathcal E}}
\def\dom{{\rm Dom}}
\def\rank{\rm rank}
\def\kj{k^{(j)}}

\begin{document}

\title {On the focusing critical semi-linear wave equation}
\author{J.\ Krieger, W.\ Schlag}
\thanks{The first author was supported by the NSF grant DMS-0401177.
The second author was partially supported by the NSF grant
DMS-0300081 and a Sloan Fellowship. He wishes to thank Sigmund
Selberg from NTNU in Trondheim, Norway, for his hospitality and the
support by the Research Council of Norway through the research
project "Partial Differential Equations and Harmonic Analysis".
Also, he wishes to thank Fritz Gesztesy for many helpful
discussions. Both authors are grateful to Carlos Kenig  and Kaihua Cai for their interest in this work and comments on a preliminary
version of this paper. They also gratefully acknowledge the support of the MSRI, Berkeley, where this work was presented
during a workshop on dispersive PDE in August of 2005.  Finally, the authors wish to thank Herbert Koch for pointing out an error in a preliminary
 formulation of Theorem~\ref{thm:main}.}
\address{Harvard University, Dept.\ of Mathematics, Science Center, 1 Oxford Street, Cambridge, MA 02138, U.S.A.}
\email{jkrieger@math.harvard.edu}
\address{Department of Mathematics, University of Chicago, 5734 South University Ave., Chicago, IL 60637, U.S.A.
}
\email{schlag@math.uchicago.edu}

\date{}
\maketitle

\begin{abstract}
The wave equation
\[ \pr_{tt} \psi -\Lap \psi - \psi^5=0  \]
in $\R^3$ is known to exhibit finite time blowup for
large data. It also admits the special static solutions
\[ \phi(x,a) =
(3a)^{\frac14} (1+a |x|^2)^{-\frac12}
\] for all $a>0$ which are linearly unstable. We view these functions as a curve in
the energy space $\dot H^1\times L^2$.
We show that in a small neighborhood of itself, which lies on
 a stable hyper-surface of radial data, this curve acts as a one-dimensional
attractor.
\end{abstract}

\section{Introduction}
\label{sec:intro}

We consider the equation \beeq\label{eq:weq} \Box \psi - \psi^5 =
\pr_{tt} \psi -\Lap \psi - \psi^5=0 \eneq in $\R^3$.
By an argument of Levine~\cite{Lev}, this equation can blow up in finite time.
In fact, negative energy leads to finite-time blowup, see Strauss~\cite{Str}.
On the other hand, this equation
admits the static solutions of Aubin
\beeq\label{eq:wave}
\phi(r,a) =
(3a)^{\frac14} (1+a r^2)^{-\frac12}
\eneq for all $a>0$, i.e., these
functions satisfy
\[ -\Lap \phi - \phi^5 =0 \]
Note that although $\phi(\cdot,a)\not\in L^2(\R^3)$, we have $\phi(\cdot,a)\in \dot{H}^1$.
In particular, $\phi(\cdot,a)$ is a finite-energy solution of~\eqref{eq:weq}.
Our goal is to understand small perturbations of $\phi$, more
precisely, to show that  stable manifolds of finite co-dimension exist.
As we shall see later, the point here is that these special solutions are linearly
unstable (i.e., the linearized operator has an unstable eigenvalue).
Seeking solutions of the form $\psi(x,t)=\phi(r,a(t))+u(x,t)$ leads
to the linearized problem \beeq\label{eq:lin} \pr_{tt} u + H(a(t)) u =
-\pr_{tt} \phi(r,a(t)) + N(u,\phi(\cdot,a(t))) \eneq where $H(a(t))= -\Lap -
5\phi^4(\cdot,a(t)) = -\Lap + V(\cdot,a(t))$ and
\[ N(\phi,u) = 10 \phi^3 u^2 + 10 \phi^2 u^3 + 5 \phi u^4 + u^5\]
By definition,
\[ H(a)\, \partial_a \phi(\cdot,a) =0, \qquad H(a) \nabla \phi(\cdot,a) =0\]
so that $H(a)$ has a resonance (due to $\pr_a \phi\not\in L^2(\R^3)$),
as well as an eigenvalue at zero energy. Moreover, since $\la
H(a)\phi,\phi\ra <0$, $H(a)$ has negative eigenvalues.

In this paper, we will only consider radial solutions of
\eqref{eq:weq}. This precludes any movement of the solution $\phi(\cdot,a)$.
More precisely, due to the translation and Lorentz invariance of~\eqref{eq:weq}, a
small, non-radial perturbation of $\phi(\cdot,a)$ will typically impart a non-zero
momentum on $\phi$ and thus lead to
a moving bulk-term, whereas in our radial case case only the dilation factor $a(t)$ can move
with time.

Let $a>0$ be fixed. 
If we consider $H(a)$ on the invariant subspace
$L^2_r(\R^3)$ of radial functions, then it is easy to describe the
spectrum as well as the resonances rigorously:  there is exactly one
negative eigenvalue and a unique resonance at zero. To see this,
write the equation $H(a)f=-k^2 f$ with a radial $f\in L^2(\R^3)$ in the
form $-g'' + V(a)g=-k^2 g$ on $L^2(0,\infty)$ where
$V(a)=V(r,a)=-5\phi^4(r,a)$ and $g(r)=rf(r)$. Then $g_1:= r\partial_a
\phi$ has a unique positive zero, which implies that there is a
simple ground-state $g>0$ with negative energy $-k^2$, and $g_1\in
L^\infty$ is the unique resonance function at zero. Note that $k>0$
 depends on the parameter $a$ through scaling: $k(a)=a^{\half}k(1)$.
This follows from $V(r,a)=aV(\sqrt{a},1)$.
 Moreover, $g$ decays
exponentially by Agmon's estimate. In what follows, we will denote the
ground state eigenfunction of $H(a)$ with $a=1$ by $g_0$ and the associated
eigenvalue by~$-k_0^2$, with $k_0>0$.

Our main result is as follows:

\begin{theorem}
\label{thm:main}  Fix $R> 1$ and let$\,$\footnote{The compact
support assumption can be replaced with power-like decay, but we
ignore this issue}
\[ X_R=\{ (f_1,f_2)\in H^3_{\rm rad}(\R^3)\times H^2_{\rm rad}(\R^3)\:|\: \supp(f_j)\subset B(0,R)\} \]
Define
\[
\Sigma_0:=\{ (f_1,f_2)\in X_R\:|\: \la k_0 f_1+f_2,g_0\ra   =0 \}
\]
and let $B_\delta(0)\subset\Sigma_0$ denote a $\delta$-ball in the topology of~$X_R$.
Then there exists $\delta=\delta(R)>0$ and a Lipschitz function
$h:B_\delta(0)\subset\Sigma_0\to \R$
with the following properties:
\begin{align}
|h(f_1,f_2)| &\less \|(f_1,f_2)\|_{X_R}^2, \qquad \forall\; (f_1,f_2)\in B_\delta(0)\nn \\
|h(f_1,f_2)-h(\til f_1, \til f_2)| &\less \delta \|(f_1,f_2)-(\til
f_1,\til f_2)\|_{X_R}\qquad \forall\; (f_1,f_2), (\til f_1,\til f_2)\in B_\delta(0) \label{eq:heta}
\end{align}
and for any $(f_1,f_2)\in \Sigma_0$ the Cauchy problem
\begin{align}
 \Box \psi - \psi^5 &= 0 \nn \\
 \psi(\cdot,0) &= \phi(\cdot,1) + f_1 + h(f_1,f_2) g_0,
 \quad \pr_t \psi(\cdot,0) = f_2
 \label{eq:mod_init}
 \end{align}
has a unique global solution of the form
\begin{equation}
\label{eq:darst}
 \psi(\cdot,t) = \phi(\cdot, a(\infty)) + v(\cdot,t)
 \end{equation}
where $|a(\infty)-1|\less\delta$.
The radiative term $v(\cdot,t)$ disperses like a free wave, i.e.,
$\|v(\cdot,t)\|_\infty \less \delta \la t\ra^{-1}$
for all $t>0$, and it also scatters like a free wave with energy data:
\[  (v,\pr_t v)(\cdot,t) = (\til v,\pr_t \til v)(\cdot,t)  + o_{\dot H^1\times L^2}(1) 
\text{\ \ as\ \  } t \to\infty\]
 with $\Box \til v =0, \; (\til v(\cdot,0),\pr_t \til v(\cdot,0))\in \dot H^1\times L^2$.
\end{theorem}

The point of \eqref{eq:mod_init} is that it describes a parametric
surface $\phi(\cdot,1)+\Sigma$ where\footnote{It is perhaps desirable to
replace $g$ by a cut-off of itself to $|x|\less R$ in~\eqref{eq:mod_init}. Because of the exponential
decay of $g$ we do not make this distinction.}
 $\Sigma\subset H^3\times H^2$ is parameterized by $(f_1,f_2)\in B_\delta(0)$.
In view of the estimates~\eqref{eq:heta}, $\Sigma$ can also be realized
as a Lipschitz graph of a function with domain~$B_\delta(0)$. Thus,
Theorem~\ref{thm:main} states the existence of a codimension one
Lipschitz graph which is the distorted image of the ball
$B_\delta(0)\subset \Sigma_0$.
Moreover, $\Sigma_0$ is the tangent plane to $\Sigma$ at zero.
Theorem~\ref{thm:main} states that
this manifold $\Sigma$ is a stable manifold for the
NLW~\eqref{eq:weq} centered at the special solution~$\phi(\cdot,1)$.
Furthermore, the theorem implies that the curve $a\mapsto\phi(\cdot,a)$ acts as a one-dimensional
attractor for data on this stable manifold. In other words,  if $(\psi(\cdot,0),\pr_t\psi(\cdot,0))\in\phi(\cdot,1)+\Sigma$, then
there is a global solution of the form
\[ \psi(\cdot,t) = \phi(\cdot,a(\infty)) + O(\delta \la t\ra^{-1})_{L^\infty} \text{\ \ as\ }t\to\infty\]
with some $|a(\infty)-1|\less \delta$.
The proof of Theorem~\ref{thm:main} proceeds by making an ansatz
\beeq
\label{eq:ansatz} \psi(\cdot,t) = \phi(\cdot, a(t)) + u(\cdot,t)
\eneq
where $a(t)$ is a path that starts with $a(0)=1$ and a dispersive part $u(\cdot, x)$,
see Section~\ref{sec:ansatz}. The solution $(a,u)$ as well as the correction $h(f_1,f_2)$
are obtained by a Banach iteration. The significance of the linear hyperplane $\Sigma_0$
is that it serves as the stable manifold for the linearized operator at the first step of the iteration. At
 each subsequent step
of the iteration we then need to make corrections to these data perpendicular to~$\Sigma_0$.
Eventually, these corrections add up to the quadratically small function~$h(f_1,f_2)$.

\noindent
Theorem~\ref{thm:main} was motivated by the stable manifold papers
\cite{Sch}, \cite{KrSch} and \cite{KrSch2}, as well as by the numerical evidence
presented in Bizo\'n, Chmaj, Tabor~\cite{Biz} and Szpak~\cite{Szp}.
It is a rigorous expression of the heuristic principle that instability
results from the simple, negative eigenvalue. It remains an open
problem to determine what happens to data in $B_\delta^+(0)$ and $B_\delta^-(0)$,
which are the two halves of $B_\delta(0)\setminus \Sigma$. According to some numerical
evidence it is conjectured in~\cite{Biz} that these
two halves should display a blow-up/scattering dichotomy.
However, our understanding of these issues is very limited\footnote{
The recent work of Merle and Zaag~\cite{MZ1}--\cite{MZ3} investigates the question of
 blow-up for semi-linear wave equations. However, their work does not concern the $H^1$ critical case
 which we study here}.
Even the question of (nonlinear)
orbital stability of $\phi(\cdot,a)$ in the radial energy space appears to be difficult.
Of course, we expect that they are orbitally unstable, but the fundamental theory of orbital stability by
Grillakis, Shatah, Strauss~\cite{GSS1}, \cite{GSS2} does
not apply to this case.

Nevertheless, it is possible to make a statement about the unique role that the manifold~$\Sigma$ plays
in $B_\delta(0)$: If data from $B_\delta(0)$ lead to a solution of the form~\eqref{eq:ansatz}, and $\dot{a}$
and $u(x,t)$ satisfy the estimates~\eqref{eq:u_decay}--\eqref{eq:adot_decay} below, then the data are
necessarily from~$\Sigma$. This follows from the proof of Proposition~\ref{prop:stab}.

Our paper is organized as follows: In Section~\ref{sec:ansatz}, we develop
the formalism of the linearized equation on which Theorem~\ref{thm:main} is based.
This will be accomplished by writing the wave equation as a Hamiltonian equation,
and by analyzing the spectrum of the linearized Hamiltonian.
Readers familiar with modulation theory will recognize some aspects of it in Section~\ref{sec:ansatz}.
However, we would like to emphasize that modulation theory does not apply
to the case of resonances, since they are not associated with
$L^2$ projections. In other words, it is not possible to project $u(\cdot,t)$ as in Theorem~\ref{thm:main}
onto $\pr_a \phi(\cdot,a)$ as a means of deriving the ODE for $a(t)$. However, the basic principle
of modulation theory still applies here: {\em the ODE for $a(t)$ ensures that the non-dispersive part of the
linearized operator resulting from the resonance is removed}.

The remedy here will be a careful analysis of the evolution operator $\f{\sin(t\sqrt{H})}{\sqrt{H}}$ for
a Schr\"odinger operator $H=-\Delta+V$ which has a resonance, but no eigenvalue at zero energy.
Recall that this means that
\[ (-\Delta+V-z^2)^{-1} = \frac{B_{-1}}{z} + B_0(z) \text{\ \ as\ \ }\Im z>0, \; z\to0 \]
where $B_{-1}$ is a rank-one operator and $B_0(z)$ is uniformly bounded as $z\to0$ in the operator norm
$L^{2,-1-\eps}(\R^3)\to L^{2,1+\eps}(\R^3)$.
Equivalently, in three dimensions, a zero energy resonance (but no zero eigenvalue) means that there is a solution $f$ of
\[ (-\Delta+V)f=0 \text{\ \ with\ \ } f\in L^{2,-\sigma}(\R^3)\setminus L^2(\R^3)\;\; \forall\, \sigma>\half
\]
In our application, $f=\pr_a \phi(\cdot,a)$ plays this role.
This analysis will be carried out in Section~\ref{sec:pointwise} and one of our main results there is
the representation formula
\[ \f{\sin(t\sqrt{H})}{\sqrt{H}}P_c = c_0\,(\psi\otimes\psi) +
\calS(t), \quad \|\calS(t)f\|_\infty \less t^{-1}
\|f\|_{W^{1,1}(\R^3)}  \]
where $P_c$ is the projection onto the continuous spectral subspace,
$c_0\ne0$, and $\psi$ is the resonance function of $H$ at zero normalized to $\int V\psi\,dx=1$.
The point here is that the singularity of the spectral measure of $H$ at zero produces
the contribution $c_0\,(\psi\otimes\psi)$, which is formally a projection, but {\em not} in $L^2$.
Particular attention here needs to be paid to the fact that the decay of the potential is exactly $\la x\ra^{-4}$.
For more details, as well as other results (e.g.,~it is essential to understand the kernel of $\cos(t\sqrt{H})$
and estimates on it)
we refer the reader to Section~\ref{sec:pointwise}.
As explained above, our $H$
does have eigenvalues at zero when considered as an operator on $L^2(\R^3)$. However,
its restriction to the subspace of radial functions $L^2_{\rm rad}(\R^3)$
does not. Since that subspace is an invariant subspace of $H$,
the results from Section~\ref{sec:pointwise} apply
as long as we restrict the evolution to radial functions.
Equipped with the formalism of Section~\ref{sec:ansatz} and the estimates from Section~\ref{sec:pointwise},
we find $a(t)$, $u(x,t)$, as well as $h(u_0)$ by means of
a contraction argument in Sections~\ref{sec:stability} and~\ref{sec:2}.

\section{The linearized problem}
\label{sec:ansatz}

Our goal is  to solve the Cauchy problem
\begin{align}
\pr_{tt} u + H(a(\infty)) u &= -\pr_{tt} \phi(\cdot,a(t)) +
(V(\cdot,a(\infty))-V(\cdot,a(t)) u + N(u,\phi(\cdot,a(t))) \label{eq:lin2}\\
u(\cdot,0) &= w_1, \quad \pr_t u(\cdot,0) = w_2 \label{eq:lin_init}
\end{align}
globally in time for small radial data $(w_1,w_2)$ that lie on a suitable
manifold. Of course, $(w_1,w_2)$ will eventually equal the expression on the right-hand side of
\eqref{eq:mod_init}. We are assuming here that $a(t)\to a(\infty)$ as
$t\to\infty$ with $0<a(\infty)<\infty$. We will use $H_\infty$ and
$H(a(\infty))$ interchangeably. As noted previously, $H_\infty$ has a
(unique) negative (radial) ground-state $g_\infty>0$ with
$\|g_\infty\|_2=1$. It decays exponentially by Agmon's estimate.
There is no other negative spectrum. Let $H_\infty g_\infty = -k_\infty^2
g_\infty$. Here $k_\infty=\sqrt{a(\infty)}\,k_0$  because
of scaling, see above. We recast the linearized equation~\eqref{eq:lin2} as a
Hamiltonian system, with $J=\bm 0&1\\-1&0\endm$:
\begin{align}
\pr_t \binom{u_1}{u_2} &= J\bm H_\infty &0\\ 0& 1\endm
\binom{u_1}{u_2} + \binom{0}{-\pr_{tt}\phi(\cdot,a(t)) +(V(\cdot,a(\infty))-V(\cdot,a(t))) u_1+
N(u_1,\phi(\cdot,a(t)))}
\label{eq:sys}\\
\binom{u_1}{u_2}(0) &= \binom{w_1}{w_2}\nn
\end{align}
which we write more succinctly as
\[ \dot{U}=J\Hil_\infty U + W, \quad U(0)=U_0\]
It is easy to check that\footnote{Recall that we are only
considering radial functions.}
\[ \spec(J\Hil_\infty)=i\R\cup\{\pm k_\infty\}\]
with $\pm k_\infty$ being simple eigenvalues (there are no other
eigenvalues). We write $k$ instead of $k_\infty$ for simplicity.
The eigenfunctions are\footnote{Recall $\|g_\infty\|_2=1$}
\[G_{\pm}:= (2k)^{-\f12}\bin{g_\infty}{\pm kg_\infty}\] They satisfy
\[ J\Hil_\infty G_{\pm} = \pm k G_{\pm} \]
Similarly, one easily checks that $G_{\pm}^*:= JG_{\mp}$ are the
dual eigenfunctions, i.e., \[ (J\Hil_\infty)^* G_{\pm}^* = \pm k
G_{\pm}^*\] Since we have normalized $G_{\pm}$ such that $\la
G_{\pm}, JG_{\mp}\ra=\mp 1$, this implies that the Riesz projections
$P_\pm$  onto the discrete spectrum are
\[ P_{\pm } = \mp\la \cdot, G_{\pm}^*\ra G_{\pm} = \mp\la \cdot, JG_{\mp} \ra G_{\pm}
\]
Now define $n_{\pm}=n_{\pm}(t)$ by \beeq \label{eq:npm_def} n_{\pm}(t)
G_{\pm}(\cdot) := P_{\pm } U(\cdot,t) \eneq where $U(\cdot,a)$ solves \eqref{eq:sys}. Then
one checks that $n_{\pm} = \mp \la U,JG_{\mp}\ra$ solve
\beeq\label{eq:n_ode} \dot{n}_{\pm}(t) \mp k n_{\pm}(t) =   \pm \la
W,JG_{\mp}\ra=: F_{\pm}(t) \eneq We now recall a trivial fact from
ODEs:

\begin{lemma}
\label{lem:ODE_stable} Consider the two-dimensional ODE
\[ \dot{x}(t)-A_0 x(t) = y(t),\qquad x(0)=\binom{x_1(0)}{x_2(0)}\]
where $y=\binom{y_1}{y_2}\in L^\infty([0,\infty),\Compl^2)$ and $
A_0 = \bm k & 0 \\ 0 & -k \endm $ where $k>0$. Then
$x(t)=\binom{x_1(t)}{x_2(t)}$ remains bounded for all times iff
\beeq \label{eq:stable} 0 = x_1(0)+\int_0^\infty e^{-k t} y_1(t)\,
dt. \eneq Moreover, in that case
\beeq
 x_1(t) = - \int_t^\infty e^{-(s-t)k} y_1(s)\, ds, \quad \label{eq:Duh_stab} 
 x_2(t) = e^{-tk}x_2(0)+ \int_0^t e^{-(t-s)k} y_2(s)\, ds. 
\eneq
for all $t\ge0$. In particular, if $y_1(s),y_2(s)$ decay like $\la
s\ra^{-\beta}$ with some $\beta>0$, then $x_1(t),x_2(t)$ decay at
least as fast.
\end{lemma}
\begin{proof}
Clearly, $x_1(t) = e^{t\sigma}x_1(0) + \int_0^t e^{(t-s)\sigma}
y_1(s)\, ds $ and $ x_2(t) = e^{-t\sigma}x_2(0) + \int_0^t
e^{-(t-s)\sigma} y_2(s)\, ds$. If
\[\lim_{t\to\infty}
e^{-t\sigma}x_1(t)=0\] then \[ 0 = x_1(0) + \int_0^\infty
e^{-s\sigma} y_1(s)\, ds\] which is \eqref{eq:stable}. Conversely,
if this holds, then \[ x_1(t) = -e^{t\sigma} \int_t^\infty
e^{-s\sigma } y_1(s)\, ds\] and the lemma is proved.
 \end{proof}

This lemma leads to the {\em stability condition} \beeq
\label{eq:stab_cond} 0= n_+(0) + \il_0^\infty e^{-sk} F_+(s)\, ds
\eneq for the linearized problem. We pass to the decomposition
\beeq\label{eq:U_decomp} U = n_+ G_+ + n_- G_- + \tilde{U} \eneq
where $\tilde{U}$ is the projection of $U$ onto the essential
spectrum of $J\Hil$, i.e., $\tilde U=(I-P_+-P_-)U=:P_e U$ (here "e"
refers to the essential spectrum). Written out in components,
\begin{align*} (P_++P_-)\bin{u}{v} &= (2k)^{-1}\Big\la \bin{u}{v},
\bin{kg_\infty}{g_\infty}\Big\ra\bin{g_\infty}{kg_\infty} +
(2k)^{-1}\Big\la \bin{u}{v}, \bin{kg_\infty}{-g_\infty}\Big\ra\bin{g_\infty}{-kg_\infty}\\
&=\bin{\la u,g_\infty\ra g_\infty}{\la v,g_\infty\ra
g_\infty}=\bin{P_{g_\infty}u}{P_{g_\infty}v}\\
P_e\bin{u}{v} &= \bin{u-\la u,g_\infty\ra g_\infty}{v-\la v,g_\infty\ra g_\infty}=\bin{P_{g_\infty}^\perp
u}{P_{g_\infty}^\perp v}\end{align*}
 Projecting~\eqref{eq:sys}
yields\footnote{Note that $\Hil_\infty$ and $P_{\pm}$ are associated with
$a=a(\infty)$ and thus do not depend on time. In particular,
$P_{g_\infty}^\perp \phi_a(\cdot,a(\infty))=\phi_a(\cdot,a(\infty))$. } \beeq
\label{eq:dis_pde} \pr_t \tilde U = J\Hil_\infty \tilde U + P_e W,\quad
\tilde U(0)=P_e U(0) \eneq In view of the preceding, the
 propagator takes the form
 \beeq\label{eq:free_prop}
e^{tJ\Hil_\infty}P_e = \bm \cos(t\sqrt{H_\infty})\,P_{g_\infty}^\perp
& \f{\sin(t\sqrt{H_\infty})}{\sqrt{H_\infty }}\,P_{g_\infty}^\perp \\
-\sqrt{H_\infty}\sin(t\sqrt{H_\infty})\,P_{g_\infty}^\perp & \cos(t\sqrt{H_\infty})\,P_{g_\infty}^\perp \endm \eneq
with $I-P_{g_\infty}=P_{g_\infty}^\perp$ being the projection onto the continuous
spectrum of~$H_\infty$.  The solution of \eqref{eq:dis_pde} is
\begin{equation}
\label{eq:sys_duh}
\tilde{U}(t) = e^{tJ\Hil_\infty} P_e\tilde U(0) + \int_0^t e^{(t-s)J\Hil_\infty}
P_e W(s)\, ds
\end{equation}
 If we set $\tilde U =\bin{\til u}{\til v}$
in~\eqref{eq:dis_pde}, then $\til v=\pr_t \til u$ and
\eqref{eq:dis_pde} is equivalent with
\begin{align}
\pr_{tt} \til u + H_\infty \til u &= P_{g_\infty}^\perp [-\pr_{tt}\phi(\cdot,a(t))
+(V(\cdot,a(\infty))-V(\cdot,a(t))) u+
N(u,\phi(\cdot,a(t)))] \label{eq:proj_sys}\\
\til u(0) &= P_{g_\infty}^\perp w_1,\; \pr_t \til u(0) = P_{g_\infty}^\perp w_2\nn
\end{align}
Here $u$ is the solution of the full equation~\eqref{eq:lin2}, i.e.,
\[ u(\cdot,t) = (2k)^{-\f12}(n_+(t)+n_-(t))g_\infty + \til u(\cdot,t)\]
and $H_\infty$ as well as $g_\infty$ are to be taken relative to $a=a(\infty)$.

The solution of~\eqref{eq:proj_sys} is
\begin{align} \til u(\cdot,t) &=
\cos(t\sqrt{H_\infty})P_{g_\infty}^\perp w_1 +
\f{\sin(t\sqrt{H_\infty})}{\sqrt{H_\infty}}P_{g_\infty}^\perp w_2 \nn \\
& \quad + \il_0^t \f{\sin((t-s)\sqrt{H_\infty})}{\sqrt{H_\infty}} P_{g_\infty}^\perp
[-\pr_{s}(\dot{a}(s)\phi_a(\cdot,a(s))) +(V(\cdot,a(\infty))-V(\cdot,a(s)))
u(\cdot,s)+ N(u,\phi(\cdot,a(s)))]\, ds\label{eq:sys2}
\end{align}
Now
\begin{align}
& \il_0^t \f{\sin((t-s)\sqrt{H_\infty})}{\sqrt{H_\infty}} P_{g_\infty}^\perp
\pr_{s}(\dot{a}(s)\phi_a(\cdot,a(s))) \, ds \nn \\
& = - \dot{a}(0)\,\f{\sin(t\sqrt{H_\infty})}{\sqrt{H_\infty}} P_{g_\infty}^\perp
\phi_a(\cdot,a(0)) + \il_0^t {\cos((t-s)\sqrt{H_\infty})} P_{g_\infty}^\perp
\,\phi_a(\cdot,a(s)) \,\dot{a}(s)\, ds \label{eq:sys3}
\end{align}
Let us first continue with a model case. Since (with
$\psi=\phi_a(\cdot,a(\infty))$)
\[ (\pr_{tt} + H_\infty)\psi =0, \qquad (\pr_{tt} + H_\infty)
\;t\psi =0\] we obtain
\[ \cos(t\sqrt{H_\infty})\, \psi =
\psi, \quad \f{\sin(t\sqrt{H_\infty})}{\sqrt{H_\infty}}\,\psi = t\psi
\]
from which we conclude that
\begin{align*}
& - \dot{a}(0)\,\f{\sin(t\sqrt{H_\infty})}{\sqrt{H_\infty}}P_{g_\infty}^\perp \psi +
\il_0^t {\cos((t-s)\sqrt{H_\infty})} P_{g_\infty}^\perp
\,\psi \,\dot{a}(s)\, ds \\
&= -  t\dot{a}(0)\, \psi + \il_0^t \psi \,\dot{a}(s)\, ds
 = - t\dot{a}(0)\, \psi +
(a(t)-a(0))\psi
\end{align*}
To pass from the model case to the real one, we use that for all
Schwartz functions~$f$
\begin{align}
\f{\sin(t\sqrt{H_\infty})}{\sqrt{H_\infty}}P_{g_\infty}^\perp &= c_0\,(\psi\otimes\psi) +
\calS(t), \quad \|\calS(t)f\|_\infty \less t^{-1}
\|f\|_{W^{1,1}(\R^3)} \label{eq:sin_2} \\
\|\cos(t\sqrt{H_\infty}) P_{g_\infty}^\perp f\|_\infty & \less
t^{-1}\|f\|_{W^{2,1}(\R^3)}\label{eq:cos_2}
\end{align}
We will need the following more precise statement concerning the
$L^1$ norm on the right-hand side of~\eqref{eq:cos_2}: There exists
a kernel $K_t(x,y)$ so that \beeq\label{eq:Kt} |K_t(x,y)|\less
(\chi_{[|x|+|y|>t]}+\la t\ra^{-1}) (\la x\ra \la y\ra)^{-1} \eneq
and
\[ \|(\cos(t\sqrt{H_\infty}) P_{g_\infty}^\perp - K_t) f\|_\infty \less
t^{-1}(\|\nabla f\|_{L^1(\R^3)}+\|D^2 f\|_{L^1(\R^3)})
\]
Notice that $\|K_t f\|_\infty \less t^{-1} \|f\|_1 $, but we will
also need to apply $K_t$ to functions not in $L^1$. More precisely,
we shall use that
\[ \sup_{x\in\R^3}\Big|\il_{\R^3} K_t(x,y)\, \la y\ra^{-3}\, dy \Big|\less \la
t\ra^{-1} \] The bounds \eqref{eq:sin_2} and \eqref{eq:cos_2} are
proved in Section~\ref{sec:pointwise}. Since
\[ \phi_a(r,a) = -\f14\,3^{\f14}a^{-\f54} r^{-1} + O(r^{-3}) \qquad \text{as\
} r\to\infty \] we conclude that (at least if $a_1,a_2 \in (1/2,2)$)
\[
\phi_a(r,a_1) = (a_2/a_1)^{\f54} \phi_a(r,a_2) + O(|a_1-a_2|\la
r\ra^{-3}) \qquad \text{as\ }r\to\infty
\]
and the $O$-term satisfies symbol-type estimates under
differentiation.  Hence, returning to~\eqref{eq:sys3}, we obtain
\begin{align*}
\f{\sin(t\sqrt{H_\infty})}{\sqrt{H_\infty}}P_{g_\infty}^\perp \phi_a(\cdot,a(0)) & =
\f{\sin(t\sqrt{H_\infty})}{\sqrt{H_\infty}}P_{g_\infty}^\perp [ (a(\infty)/a(0))^{\f54}\psi
+ O(\la r\ra^{-3}) ]
\\
& =t (a(\infty)/a(0))^{\f54}\psi + (c_0\,
(\psi\otimes \psi) + \calS(t)) O(\la r\ra^{-3}) \\
& =t (a(\infty)/a(0))^{\f54}\psi + \Omega_1(t)
\end{align*}
In view of \eqref{eq:sin_2} and the bounds on $K_t$,
\begin{equation}
\label{eq:omega1}\sup_{t\ge0} \|\Omega_1(t)\|_\infty < \infty
\end{equation}
Similarly,
\[
\cos((t-s)\sqrt{H_\infty}) P_{g_\infty}^\perp \,\phi_a(\cdot,a(s))
= (a(\infty)/a(s))^{\f54}\psi + \Omega_2(t,s)
\]
where
\[ \Omega_2(t,s) = \cos((t-s)\sqrt{H_\infty}) P_{g_\infty}^\perp
[\phi_a(\cdot,a(s))-(a(\infty)/a(s))^{\f54}\psi]
\]
is again bounded (and small) since
\[ |\phi_a(x,a(s))-(a(\infty)/a(s))^{\f54}\psi(x) | \less |a(s)-a(\infty)|\la x\ra^{-3}
\]
Therefore,
\begin{align*} \til u(\cdot,t) &=
\cos(t\sqrt{H_\infty})P_{g_\infty}^\perp w_1 +
\f{\sin(t\sqrt{H_\infty})}{\sqrt{H_\infty}}P_{g_\infty}^\perp w_2  \nn \\
& \quad   + \dot{a}(0) (t
(a(\infty)/a(0))^{\f54}\psi
 + \Omega_1(t))
 - \il_0^t [ (a(\infty)/a(s))^{\f54}\psi + \Omega_2(t,s)]\dot{a}(s)\, ds \nn \\
 & \quad -
\il_0^t \f{\sin((t-s)\sqrt{H_\infty})}{\sqrt{H_\infty}} P_{g_\infty}^\perp
[(V(\cdot,a(\infty))-V(\cdot,a(s))) u(\cdot,s)+ N(u,\phi(\cdot,a(s)))]\, ds \nn
\end{align*}
which we can further rewrite as
\begin{align}
\til u(\cdot,t) & = \cos(t\sqrt{H_\infty})P_{g_\infty}^\perp w_1 +
\calS(t) P_{g_\infty}^\perp w_2
 - \il_0^t \cos((t-s)\sqrt{H_\infty}) P_{g_\infty}^\perp
[\phi_a(\cdot,a(s))-(a(\infty)/a(s))^{\f54}\psi]\,\dot{a}(s)\, ds \nn \\
& \quad-\il_0^t \calS(t-s) P_{g_\infty}^\perp
[(V(\cdot,a(\infty))-V(\cdot,a(s))) u(\cdot,s)+ N(u,\phi(\cdot,a(s))) \big]\, ds \nn \\
 & \quad +  \dot{a}(0) (t
(a(\infty)/a(0))^{\f54}\psi
 + \Omega_1(t)) + \psi \Big\{ c_0 \la \psi, w_2 \ra + 4 (a(\infty))^{\f54} (a(t)^{-\f14}-a(0)^{-\f14}) \label{eq:nondisp1}\\
& \qquad\qquad\qquad\qquad- c_0\,
\il_0^t \la \psi,(V(\cdot,a(\infty))-V(\cdot,a(s))) u(\cdot,s)+ N(u,\phi(\cdot,a(s)))\ra\, ds \Big\} \label{eq:nondisp2}
\end{align}
The equation for $a(t)$ is now determined from the requirement that
\eqref{eq:nondisp1} and~\eqref{eq:nondisp2} need to decay in time. In particular, this forces
$\dot{a}(0)=0$. The equation for $a(t)$ can be determined by setting the expression in braces equal to zero.
However, this would impose  compatibility conditions at $t=0$ which we cannot fulfil without further assumptions
on the data (e.g., $\la \psi,w_2\ra=0$). Instead, we only require that it vanishes for times $>1$. More precisely,
the equation for $a(t)$ is
\begin{align}
m_1\, \omega(t) + m_2\, t\omega(t) &= c_0\la \psi,w_2\ra + 4\;a(\infty)^{\f54}
(a^{-\f14}(t)-a^{-\f14}(0)) \nn \\
& \qquad - c_0 \il_0^t \la \psi,
(V(\cdot,a(\infty))-V(\cdot,a(s)))u(\cdot,s)+ N(u(\cdot,s),\phi(\cdot,a(s)))\ra\, ds \label{eq:a_eq}
\end{align}
where $\omega(t)$ is a fixed smooth function on $[0,\infty)$ with $\omega(t)=1$ for $0\le t\le \half$,
and $\omega(t)=0$ if $t\ge1$. The constants $m_1$ and $m_2$ will be specified shortly.
Despite the fact that $\psi$ only decays like $r^{-1}$, the
scalar product on the right-hand side of~\eqref{eq:a_eq} is well-defined due  to the
decay  of $V,\psi, u$. Indeed, we will show below that
\[ |u(x,t)|\less \delta \la x\ra^{-1} \]
see \eqref{eq:u_est}.

\noindent The constants $m_1$ and $m_2$ are determined as follows. First, setting $t=0$ leads to
the condition \beeq\label{eq:orth_cond} m_1 = c_0\la w_2,\psi \ra  \eneq
Second, the requirement $\dot{a}(0)=0$ implies
\begin{equation}
\label{eq:u(0)} m_2  = -c_0\,\la \psi, (V(\cdot,a(\infty))-V(0))w_1 + N(w_1,
\phi(\cdot,a(0)))\ra
\end{equation}
Hence the equation for $\tilde u$ now reads
\begin{align} \til u(\cdot,t)
&= \cos(t\sqrt{H_\infty})P_{g_\infty}^\perp w_1 + \calS(t)  P_{g_\infty}^\perp\,w_2
 - \il_0^t \dot{a}(s)\, \cos((t-s)\sqrt{H_\infty}) P_{g_\infty}^\perp
[\phi_a(\cdot,a(s))-(a(\infty)/a(s))^{\f54}\psi]\, ds \nn \\
&\qquad \qquad  - \il_0^t \calS(t-s)\, P_{g_\infty}^\perp [(V(\cdot,a(\infty))-V(\cdot,a(s)))
u(\cdot,s)+ N(u,\phi(\cdot,a(s)))]\, ds \label{eq:u_disp} \\
& \qquad \qquad+ c_0\big[ \la w_2,\psi \ra - t \,\la \psi, (V(\cdot,a(\infty))-V(0))w_1 + N(w_1,
\phi(\cdot,a(0)))\ra \big]\omega(t)\psi \label{eq:fix}
\end{align}
and we need to prove that the right-hand side here is dispersive in a
suitable sense. This will require proving the basic collection of
estimates for the wave equation, i.e., energy, dispersive, and
possibly also Strichartz (for finer results than the ones presented
here). Energy is the same as usual: Let $(\pr_{tt}+H_\infty)u=0$. Then
\begin{align*}
\f{d}{dt} \f12\;[\|\pr_t u(\cdot,t)\|_2^2 + \|\sqrt{H_\infty}\,
u(\cdot,t)\|_2^2] &= \la
u_{tt},u_t(t)\ra + \la \sqrt{H_\infty}\, u(\cdot,t), \sqrt{H_\infty}\, u_t(t) \ra \\
&= \la u_{tt} + H_\infty\, u , u_t(t) \ra =0
\end{align*}
so that for all $t\ge0$
\[ \|\pr_t u(\cdot,t)\|_2^2 + \|\sqrt{H_\infty }\, u(\cdot,t)\|_2^2 = \const
\]
We further remark that for any $f$, with $H_\infty=-\Delta+V_\infty$,
\begin{align}
 \|\sqrt{H_\infty} f\|_2^2
&= \la H_\infty\, f,f\ra = \|\nabla f\|_2^2 + \la V_\infty\, f,f\ra \nn \\
&\less \|\nabla f\|_2^2 + \|V_\infty\,\|_{\f32} \|f\|_6^2 \less \|\nabla
f\|_2^2 \label{eq:grad1}
\end{align}
by Sobolev imbedding. Note, however, that the reverse inequality
here cannot hold because of the resonance function. However, there
is a replacement:
\begin{align}
\|\nabla f\|_2^2 &= \la -\Laplace f,f\ra = \la H_\infty\, f,f \ra - \la V_\infty\,
f,f\ra \nn \\
&\less \|\sqrt{H_\infty}\, f\|_2^2 + \| |V_\infty|^{\half} f\|_2^2 \label{eq:grad2}
\end{align}
Similarly, \beeq\label{eq:D2} \|D^2 f\|_2 \less \|H_\infty\, f\|_2 + \|V_\infty\, f\|_2,
\qquad \|H_\infty f\|_2 \less \|D^2 f\|_2 + \|V_\infty\, f\|_2 \eneq The dispersive
estimates are proved in Section~\ref{sec:pointwise}, but the
Strichartz estimates are still lacking. It seems that one needs to
develop a suitable Littlewood-Paley theorem in the perturbed
setting. For the radial case and $H_\infty$ this can be done, albeit only
in the range $3/2<p<3$, see~\cite{Sch3}. This range is optimal due
to the resonance at zero energy. On the other hand, we should be
looking for a Littlewood-Paley theory not relative to all of $H_\infty$,
but only for functions which are orthogonal to the resonance, or
more precisely, which correspond to the regular part of the spectral
measure. In that case, the full range $1<p<\infty$ should again be
available.

\section{The contraction scheme: stability}
\label{sec:stability}

We will keep the radius $R> 1$ in Theorem~\ref{thm:main} fixed.
Constants will be allowed to depend on it.

\begin{defi} Let $Y_{R,\delta}$ denote the metric space of
$a\in C^1([0,\infty),\R^+)$ and $u\in C([0,\infty),H^2_{\rm
rad}(\R^3))$  satisfying the following  properties: For all $t\ge0$,
\begin{align}
\| u(\cdot,t)\|_\infty &\le \delta \la t\ra^{-1} \label{eq:u_decay}\\
\|\nabla u(\cdot,t) \|_{2+\infty} &\le \delta \la t\ra^{-\eps} \label{eq:nablau_decay}\\
\|\nabla u(\cdot,t)\|_2 + \|D^2 u(\cdot,t)\|_2 &\le \delta \label{eq:L2}\\
|u(x,t)| &\le C_1\, \delta \la x\ra^{-1} \label{eq:u_est}\\
|\dot{a}(t)| &\le  \delta \la t\ra ^{-2}
\label{eq:adot_decay}
\end{align}
Here $\delta>0, \eps>0$ are small$\,$\footnote{$\eps>0$ is a small
positive constant that is fixed once and for all, whereas $\delta>0$
is small but arbitrary.}, and $C_1>1$ is some constant that
does not depend on $\delta$. In addition, $a(0)=1$ and $\dot{a}(0)=0$.
\end{defi}

Let $u_0:=(f_1,f_2) \in B_{\delta_0}(0)\subset\Sigma_0$, where
$B_{\delta_0}(0)$ is a $\delta_0$-ball in $\Sigma_0$ centered at
zero of radius $\delta_0\ll \delta$, see Theorem~\ref{thm:main}. Our
goal is to find a fixed point for the map
\[\Phi=\Phi_{u_0}:
Y_{R,\delta}\to Y_{R,\delta}, \; (u,a)\mapsto (v,b)
\] which we now describe (and which depends on the choice of $u_0$).
We intend to show that for a suitable -- and unique -- choice of
$h(u_0;u,a)$  the solution $v,b$ of (with
$\psi_0= \pr_a\phi(\cdot,1)$ and $c_0$ as in~\eqref{eq:sin_2})
\begin{align}
&\pr_{tt} v + H(a(\infty)) v = -\pr_{tt} \phi(\cdot,b(t)) +
(V(\cdot,a(\infty))-V(\cdot,a(t))) u + N(u,\phi(\cdot,a(t))) \label{eq:v_eq}\\
& 4\;a(\infty)^{\f54} (b^{-\f14}(t)-1) =\nn\\
& =  c_0 \il_0^t \la \phi_a(\cdot,a(\infty)),
(V(\cdot,a(\infty))-V(\cdot,a(s)))u(\cdot,s)+
N(u(\cdot,s),\phi(\cdot,a(s)))\ra\, ds  - c_0\la f_2,\phi_a(\cdot,a(\infty)) \ra\label{eq:b_eq}\\
&\quad  + c_0\big[ \la f_2,\phi_a(\cdot,a(\infty)) \ra - t \,\la \phi_a(\cdot,a(\infty)), (V(\cdot,a(\infty))-V(\cdot,a(0)))w_1 + N(w_1,
\phi(\cdot,a(0)))\ra \big]\omega(t) \nn \\
& v(0) = f_1 + h(u_0;u,a) g_0 =: w_1, \quad \pr_t v(0) =  f_2
  \label{eq:init_cond1}\\
& b(0) =1, \quad \dot{b}(0)=0 \label{eq:init_cond2}
\end{align}
satisfies the same bounds \eqref{eq:u_decay}--\eqref{eq:adot_decay}
as $u,a$. It is with this choice of $h$ that the map
$\Phi_{u_0}$ is being defined. This system should be compared
to~\eqref{eq:a_eq} and~\eqref{eq:u_disp}. The choice of $h$
will  be based on the stability condition~\eqref{eq:stab_cond}.

\begin{prop}
\label{prop:stab} There exists $0<\delta_0\ll \delta$
small\mbox{$\,$}\footnote{This means that $\delta_0=c\delta$ where
$c$ is a small absolute constant.} so that for any $u_0 =(f_1,f_2)\in
B_{\delta_0}(0)\subset\Sigma_0$ the following holds: For any
$(u,a)\in Y_{R,\delta}$ there is a unique choice of $h(u_0;u,a)$
 so that $\Phi_{u_0}(u,a)\in Y_{R,\delta}$. Moreover,
\[ |h(u_0;u,a)|\less \delta^2 \]
as well as \beeq \label{eq:hetadiff} |h(u_0;u,a)-h(\tilde u_0;u,a)|
\less \delta\|u_0-\tilde
u_0\|_{H^3\times H^2} \eneq for all $u_0,\tilde u_0\in B_{\delta_0}(0)$ and
$(u,a)\in Y_{R,\delta}$. In particular,
\beeq\label{eq:Philip}\|\Phi_{u_0}(u,a)-\Phi_{\til
u_0}(u,a)\|_{Y_{R,\delta}} \less \|u_0-\til u_0\|_{H^3\times H^2}\eneq for all
$u_0,\til u_0\in B_{\delta_0}(0)$ and $(u,a) \in Y_{R,\delta}$.
\end{prop}
\begin{proof}
We begin by checking that \eqref{eq:b_eq} reproduces the decay of
$\dot{a}$ under these assumptions. In view of \eqref{eq:u_decay} and
\eqref{eq:u_est}, as well as the definition of $\phi$ and $V$,
\beeq\label{eq:ab1} |(V(\cdot,a(\infty))-V(\cdot,a(t)))u(x,t)|\less
\delta^2 \la t\ra^{-2} \la x\ra^{-4}  \eneq as well as
\beeq\label{eq:ab2} |N(u(\cdot,t),\phi(\cdot,a(t)))| \less C_1^3 \,
\delta^2 \la t\ra^{-2}\, \la x\ra^{-3} \eneq Let us first discuss
the solvability of~\eqref{eq:b_eq}. The equation (in this proof $\psi=\phi_a(\cdot,a(\infty))$)
\begin{align*}
&a(\infty)^{\f54} (b^{-\f14}(\infty)-1) \\
& =  \f{c_0}{4} \il_0^\infty \la \phi_a(\cdot,a(\infty)),
(V(\cdot,a(\infty))-V(\cdot,a(s)))u(\cdot,s)+
N(u(\cdot,s),\phi(\cdot,a(s)))\ra\, ds  - c_0\la f_2,\psi\ra
\end{align*}
has a unique solution $b(\infty)=1+O(\delta)$ since the right-hand
side here is $O(\delta)$. Hence, \eqref{eq:b_eq} has a
well-defined solution $b(t)$ for all $t\ge0$ with the property that
$|b(t)-1|\less \delta$ for all $t\ge0$. By construction, $b(0)=1$
and $\dot{b}(0)=0$.
Second,
\begin{align}\label{eq:adot}
\dot{b}(t) &= - c_0\, (b(t)/a(\infty))^{\frac54}\Big\{ \la
\psi, (V(\cdot,a(\infty))-V(\cdot,a(t)))u(\cdot,t)+
N(u(\cdot,t),\phi(\cdot,a(t)))\ra \\
& \quad - \big[ \la f_2,\psi \ra - t \,\la \psi, (V(\cdot,a(\infty))-V(\cdot,a(0)))w_1 + N(w_1,
\phi(\cdot,a(0)))\ra \big]\dot\omega(t) \nn\\
& \quad +\la \psi, (V(\cdot,a(\infty))-V(\cdot,a(0)))w_1 + N(w_1,
\phi(\cdot,a(0)))\ra \omega(t) \Big\} \nn
\end{align}
The bounds \eqref{eq:ab1}, \eqref{eq:ab2} in conjunction with~\eqref{eq:adot}
reproduce~\eqref{eq:adot_decay} for small~$\delta$. Observe that on the support of $\omega$
the equation for $\dot{b}$ contains the unknown $h$. However, only a very crude bound $|h|\less \delta$,
say,
is required to obtain the estimate we need. We will comment on this issue later, when we
solve for~$h$.
We remark that $|\dot{b}(t)|\less\delta^2$ outside
of the support of $\dot{\omega}$. It is because of $|\la f_2,\psi\ra|\less \delta$
that we only obtain $\delta$ on the support of $\dot{\omega}$. In passing, we also
remark that we cannot replace $t^{-1}$ by $t^{-\alpha}$ with
$\alpha<1$ in~\eqref{eq:u_decay} since that would mean that at best
$|\dot{a}(t)|\less \delta \la t\ra^{-2\alpha}$. However, that rate
of decay cannot be reproduced from~\eqref{eq:a_eq} because of the
term $(V(\cdot,a(\infty))-V(\cdot,a(s))) u(\cdot,s)$. Indeed, that term would
contribute $s^{1-3\alpha}$ which is worse than $s^{-2\alpha}$
because of $\alpha<1$. So we work with $t^{-1}$ for the decay of
$\|u(\cdot,t)\|_\infty$ (which is the best possible rate of decay).

\noindent Next, we turn to $v$. As in Section~\ref{sec:ansatz} we
write
\begin{equation} \label{eq:v_split} v(t) = (2k_\infty)^{-\half}
(n_+(t)+ n_-(t))g(\cdot,a(\infty)) + \tilde{v}(t) \end{equation}
with $\tilde v(t)= P_{g(\cdot,a(\infty))}^\perp\,v(t)$
(recall that $k_\infty=k(a(\infty))=a(\infty)k_0$). The finite-dimensional
part satisfies, cf.~\eqref{eq:n_ode},
\begin{equation}
\label{eq:npm2} \dot{n}_{\pm}(t) \mp k_\infty n_{\pm}(t) = \mp\big\la -\pr_{tt}
\phi(\cdot,b(t)) + (V(\cdot,a(\infty))-V(\cdot,a(t))) u +
N(u,\phi(\cdot,a(t))), g(\cdot,a(\infty)) \big\ra =: \mp F(t)
\end{equation}
whereas the remaining part $\tilde v$ satisfies (with $H=-\Delta+V(\cdot,a(\infty))$
 throughout this proof)
\begin{align}
\til v(t) &= \cos(t\sqrt{H})P_{g(\cdot,a(\infty))}^\perp\, w_1 + \calS(t)  P_{g(\cdot,a(\infty))}^\perp\, f_2 \nn \\
& \quad  - \il_0^t \dot{b}(s)\, \cos((t-s)\sqrt{H})
 P_{g(\cdot,a(\infty))}^\perp\;
[\phi_b(\cdot,b(s))-(a(\infty)/b(s))^{\f54}\psi]\, ds \nn \\
&\quad   - \il_0^t \calS(t-s)\, P_{g(\cdot,a(\infty))}^\perp
[(V(\cdot,a(\infty))-V(\cdot,a(s))) u(\cdot,s)+
N(u,\phi(\cdot,a(s)))]\, ds \label{eq:v_disp} \\
& \quad + c_0\big[ \la f_2,\psi \ra - t \,\la \psi, (V(\cdot,a(\infty))-V(\cdot,a(0)))w_1 + N(w_1,
\phi(\cdot,a(0)))\ra \big]\omega(t)\psi \nn
\end{align}
cf.~\eqref{eq:u_disp}. We remark that the final term on the right-hand side,
which is a multiple of $\psi$ of size $\ll\delta$, does not affect any of the estimates
\eqref{eq:u_decay}--\eqref{eq:u_est}. Hence, we will ignore it from now on.

\noindent In order to avoid exponential
growth of $n_+(t)$ it is both necessary and sufficient that
\beeq\label{eq:con2} 0= n_+(0) + \il_0^\infty e^{-tk_\infty} F(t)\, dt
\eneq see~\eqref{eq:stab_cond}. We need to transform~\eqref{eq:con2}
into a condition on $(w_1,f_2)$. To do so, note that \eqref{eq:v_split}
and~\eqref{eq:init_cond1} imply that
\begin{align*}
n_+(0) + n_-(0) & = (2k_\infty)^\half \la w_1, g(\cdot,a(\infty))\ra \\
\dot{n}_+(0) + \dot{n}_-(0)&= (2k_\infty)^\half \la f_2, g(\cdot,a(\infty))\ra
\end{align*}
whereas we deduce from \eqref{eq:npm2} that
\[ \dot{n}_+(0) + \dot{n}_-(0) = k_\infty(n_+(0)-n_-(0)) \]
The conclusion is that
\begin{align} n_+(0) &=
(2k_\infty)^{-\half} \big[\la
g(\cdot,a(\infty)), k_\infty f_1 + f_2  \ra + k_\infty h(u_0;u,a) \la g_0, g(\cdot,a(\infty))\ra\big] \nn \\
&= k_\infty\il_0^\infty \dot{b}(t)\, e^{-tk_\infty} \la \phi_b(\cdot,b(t)),
g(\cdot,a(\infty)) \ra\, dt \nn \\
&\qquad \qquad - \il_0^\infty e^{-tk_\infty} \la
(V(\cdot,a(\infty))-V(\cdot,a(t))) u + N(u,\phi(\cdot,a(t))),
g(\cdot,a(\infty)) \big\ra \, dt \label{eq:n+(0)}
\end{align}
Observe that $(f_1,f_2)\in \Sigma_0$ implies that
\begin{align*} |\la
g(\cdot,a(\infty)), k_\infty f_1 + f_2  \ra|
\le |\la
g(\cdot,a(\infty))-g_0, k_0 f_1 + f_2\ra|
+ |\la g(\cdot,a(\infty)), (k_\infty-k_0) f_1  \ra|
 \less \delta^2
\end{align*}
and
\[ \la g_0, g(\cdot,a(\infty))\ra = 1 + O(\delta), \qquad |\la \phi_b(\cdot,b(t)),
g(\cdot,a(\infty)) \ra| = |\la \phi_b(\cdot,b(t))-\phi_a(\cdot,a(\infty)),
g(\cdot,a(\infty)) \ra|\less \delta \]
As mentioned above, $\dot{b}$ depends on $h(u_0;u,a)$ through $w_1$. Thus, $h(u_0;u,a)$ also
appears on the right-hand side of~\eqref{eq:n+(0)}, and not just on the left.
However, this dependence
occurs either with a small coefficient (in fact,~$\delta$), or to higher order.
Hence, we can still solve for $h(u_0;u,a)$ by means of the implicit function theorem.

In view of our assumptions on $u,a$ and the bounds we proved on $b$,
it follows for small $\delta$ that there exists a unique choice of
$h(u_0;u,a)$  so that~\eqref{eq:n+(0)} is satisfied. Moreover,
\[ |h(u_0;u,a)|\less \delta^2  \]
as well as \begin{equation} \label{eq:npm_decay}
|{n}_+(t)|+\delta|{n}_-(t)|\less \delta^2  \la t\ra^{-2}
\end{equation}
see \eqref{eq:Duh_stab}. To bound $n_-(t)$ we use that
$ |n_+(0)| \less\delta^2$ which implies that $|n_-(0)|\less \delta$.
It is also an easy matter to check
that~\eqref{eq:hetadiff} holds, which we leave to the reader.

\medskip
We now turn to estimating $\tilde v$.  First, let $0<t\less 1$.
Then, using $\| H^{-\half} \sin(t\sqrt{H}) \;P_c\|_{2\to2} \le t$,
we read off from the equation for $\tilde v$, viz.
\begin{align*}
\til v(t) &= \cos(t\sqrt{H}) P_{g(\cdot,a(\infty))}^\perp\,w_1 +
\Big[\f{\sin(t\sqrt{H})}{\sqrt{H}}-c_0(\psi\otimes\psi)\Big]
P_{g(\cdot,a(\infty))}^\perp\,f_2  \\
&\qquad  - \il_0^t \dot{b}(s)\, \cos((t-s)\sqrt{H}) P_{g(\cdot,a(\infty))}^\perp
[\phi_b(\cdot,b(s))-(a(\infty)/b(s))^{\f54}\psi]\, ds \nn \\
&\qquad  - \il_0^t
\Big[\f{\sin((t-s)\sqrt{H})}{\sqrt{H}}-c_0(\psi\otimes\psi)\Big]\,
P_{g(\cdot,a(\infty))}^\perp [(V(\cdot,a(\infty))-V(\cdot,a(s))) u(\cdot,s)+
N(u,\phi(\cdot,a(s)))]\, ds \\
& \qquad + c_0\big[ \la f_2,\psi \ra - t \,\la \psi, (V(\cdot,a(\infty))-V(\cdot,a(0))w_1 + N(w_1,
\phi(\cdot,a(0)))\ra \big]\omega(t)\psi \nn
\end{align*}
that
\[ \tilde v(t) = w(t) + \psi\, \eta(t), \quad \|w(t)\|_2 + |\eta(t)|\ll \delta \]
for $0<t\less 1$. In particular, \beeq\label{eq:u_split}
\sup_{0<t\less 1}\|\tilde v(t)\|_{2+\infty}\ll \delta \eneq Later we
will show that~\eqref{eq:u_split}, combined with the equation for
$\tilde v$, yields
\[ \|\nabla \tilde v(t)\|_2 + \|D^2 \tilde v(t)\|_2 \ll \delta \]
for all $0<t<1$. Since $\nabla\psi\in L^2$ and $D^2\psi\in L^2$ we
conclude that in fact $w(t)\in H^2$ for small times. Therefore,
\[ \|\tilde v(t)\|_\infty \less \|w(t)\|_\infty + |\eta(t)| \less \|w(t)\|_{H^2} + |\eta(t)| \ll \delta \quad \forall
\; 0<t\less 1 \] by Sobolev imbedding. Hence, it suffices to
consider $t\gg 1$. Using the dispersive decay of
$\cos((t-s)\sqrt{H})$ in the integral involving $\dot{b}$, we obtain
the bound
\[ \delta^2\il_0^t \la t-s\ra^{-1} \la s\ra^{-2} \, ds \less \delta^2 t^{-1}\]
see Remark~\ref{rem:weakL1}. Next, the dispersive
bound on $\calS(t-s)$ yields
\[ \delta^2\il_0^{t-t^{-10}} (t-s)^{-1}  \la s\ra^{-1-\eps/2}\, ds \less \delta^2 t^{-1} \]
We are using here that
\[ \|(V(\cdot,a(\infty))-V(\cdot,a(s))) u(\cdot,s) \|_{W^{1,1}(\R^3)} \less \delta^2 \la s\ra^{-1-\eps} \] as well as
\[ \| N(u,\phi(\cdot,a(s))) \|_{W^{1,1}(\R^3)} \less \delta^2 \la s\ra^{-1-\eps/2} \]
This latter bound in turn reduces to four terms of which we consider
only the extreme cases $u^2 \phi^3$ and $u^5$ (the fact that there
is $\eps/2$ and not $\eps$ is due to a small loss through
interpolation). Since $\phi^3\in L^p$ for all $p>1$ and because
of~\eqref{eq:nablau_decay}, we have
\[ \| u^2 (s)\phi^3 \|_1 +  \|\nabla u^2(s) \phi^3 \|_1 \less \delta^2 \la s\ra^{-1-\eps/2}\] as claimed.
The $u^5$ term requires more care due to the possible growth of
$\|u(\cdot,s)\|_2$ in time. However, in view of~\eqref{eq:u_est} we
have $u^3\in L^p$ with $p>1$ so that the same arguments apply as in
the case of $\phi^3$. For the integral over $[t-t^{-10},t]$ we use
Sobolev imbedding. More precisely, we write
\begin{align}
& \il_{t-t^{-10}}^t \calS(t-s)\, P_{g(\cdot,a(\infty))}^\perp
[(V(\cdot,a(\infty))-V(\cdot,a(s))) u(\cdot,s)+
N(u,\phi(\cdot,a(s)))]\, ds \nn \\
&= \il_{t-t^{-10}}^t \frac{\sin((t-s)\sqrt{H})}{\sqrt{H}}\,
P_{g(\cdot,a(\infty))}^\perp [(V(\cdot,a(\infty))-V(\cdot,a(s)))
u(\cdot,s)+ N(u,\phi(\cdot,a(s)))]\, ds \label{eq:10one}\\
&\quad - c_0\, \il_{t-t^{-10}}^t (\psi\otimes \psi)\,
[(V(\cdot,a(\infty))-V(\cdot,a(s))) u(\cdot,s)+ N(u,\phi(\cdot,a(s)))]\, ds
\label{eq:10two}
\end{align}
The final integral \eqref{eq:10two} is estimated directly in
$L^\infty$, leading to a bound of $\delta^2 t^{-12}$. The one
in~\eqref{eq:10one} is estimated in $H^2$ by means of the same $L^2$
based arguments as above, leading to a contribution of $\delta^2 t^{-22}$.
 The conclusion  is that
\eqref{eq:u_decay} is regained for $\tilde{v}$. To extend to $v$,
simply use the bounds~\eqref{eq:npm_decay}. Next, we deal with the
$L^2$ estimates~\eqref{eq:L2}. In view of \eqref{eq:grad2},
\[ \|\nabla \tilde v\|_2 \less \|\sqrt{H} \tilde v\|_2 +
\|V\|_1 \|\tilde v\|_\infty \less \|\sqrt{H} \tilde v\|_2 + \delta
\la t\ra^{-1} \] at least for $t>1$. If $0<t<1$, then we use
$\|\tilde v(t)\|_{2+\infty}\ll \delta$, see~\eqref{eq:u_split}
instead of $\|\tilde v(t)\|_\infty$ (since we used the energy bound
above to control $\|\tilde v(t)\|_\infty$ for small times and thus
have to avoid going in circles). Bounding $\|\sqrt{H} \tilde
v(t)\|_2$ amounts to checking that
\begin{align*}
& \il_0^t |\dot{b}(s)|\, \|
[\phi_b(\cdot,b(s))-(a(\infty)/b(s))^{\f54}\psi]\|_{H^1}\, ds
 + \il_0^t \| [(V(\cdot,a(\infty))-V(\cdot,a(s)))
u(\cdot,s)+ N(u,\phi(\cdot,a(s)))]\|_2\, ds \\
& \less \delta^2 \il_0^\infty \la s\ra^{-2} \, ds \less \delta^2
\end{align*}
Here we used \eqref{eq:grad1} as well as $
\sqrt{H}\calS(t) = \calS(t) \sqrt{H} = \sin(t\sqrt{H})$ which
is bounded on $L^2$. This relation is of course a
consequence of $\sqrt{H}\psi=0$. However, since $\psi$ does not lie in
$\dom(\sqrt{H})=W^{1,2}(\R^3)$, this latter claim
requires some care and needs to be interpreted weakly. More precisely, we claim that
\[ \sqrt{H}(\chi(\cdot/R)\psi) \rightharpoonup 0   \]
in the sense of weak convergence on $L^2(\R^3)$ as $R\to\infty$ (here $\chi$ is a smooth cut-off around zero).
First, note that
\[ \|\sqrt{H}(\chi(\cdot/R)\psi)\|_2 \less \|\nabla (\chi(\cdot/R)\psi)\|_2 + \| \,|V|^{\f12}\,(\chi(\cdot/R)\psi)\|_2
<\infty
\]
uniformly in $R>1$. Consequently, it suffices to check the weak convergence
against a family of functions which is dense in~$L^2$. One such family is
$\Ran(\sqrt{H})$. This is dense
in $L^2$ since $\ker(\sqrt{H})=\{0\}$ and since $\sqrt{H}$ restricted to $\{g\}^\perp$ is self-adjoint.
Now, by the explicit decay of $\psi$ and $H\psi=0$ it follows that for any $f\in\dom(\sqrt{H})$
\[ \la \sqrt{H}(\chi(\cdot/R)\psi),\sqrt{H} f\ra = \la H (\chi(\cdot/R)\psi), f\ra \to 0 \]
as $R\to\infty$. Hence the claim.

\noindent For $\|D^2 \tilde v\|_2$ we face an additional derivative
on the right-hand side. More precisely, using \eqref{eq:D2} as well
as~\eqref{eq:grad1}, the main estimate is
\begin{align*}
& \il_0^t |\dot{b}(s)|\, \|P_{g(\cdot,a(\infty))}^\perp
[\phi_b(\cdot,b(s))-(a(\infty)/b(s))^{\f54}\psi]\|_{H^2}\, ds \\
& + \il_0^t \|\nabla P_{g(\cdot,a(\infty))}^\perp
[(V(\cdot,a(\infty))-V(\cdot,a(s)))
u(\cdot,s)+ N(u,\phi(\cdot,a(s)))]\|_{L^2}\, ds \\
& \less \delta^2 \il_0^\infty \la s\ra^{-2} \, ds +
\delta^2\il_0^\infty \la s\ra^{-1} s^{-\eps}\, ds \less \delta^2
\end{align*}
Hence, \eqref{eq:L2} follows. Finally, for \eqref{eq:nablau_decay},
we of course use the dispersive estimate as before, but with one
extra derivative. More precisely, we invoke the estimates
\begin{align}
 \|\nabla \calS(t) P_c f\|_\infty &\less t^{-1}
\|f\|_{W^{2,1}(\R^3)} \label{eq:singrad}\\
\| [\nabla \cos(t\sqrt{H}) P_c -K_t]f\|_\infty &\less t^{-1}
\sum_{1\le |\alpha|\le 3}\|D^\alpha
f\|_{L^{1}(\R^3)}\label{eq:cosgrad}
\end{align}
where the kernel $K_t$ satisfies~\eqref{eq:Kt}. These inequalities
are obtained by passing the gradient through the various expansions
in Section~\ref{sec:pointwise}. Doing so leads to commutators
between the gradient and the potential, which requires smoothness of
the potential --- which we have in our case (it is also convenient
that $V$ has a definite sign). Hence these commutators are harmless.
More details will be presented in Section~\ref{sec:pointwise}, and we now
use~\eqref{eq:singrad}, \eqref{eq:cosgrad}. We start with the initial data $f_2$
in~\eqref{eq:v_disp} (we leave the analogous estimation of~$w_1$ to the reader). If $t>1$, then
\[
\|\nabla \calS(t)  P_{g(\cdot,a(\infty))}^\perp\,f_2 \|_{2+\infty}
\less \|\nabla \calS(t) P_{g(\cdot,a(\infty))}^\perp\,f_2\|_{\infty}
\less t^{-1} \|f_2\|_{W^{2,1}(\R^3)} \less t^{-1} \|f_2\|_{H^2}
\]
by the compact support assumption on the data. Next, if $0<t<1$,
then
\begin{align*}
&\|\nabla \calS(t)  P_{g(\cdot,a(\infty))}^\perp\,f_2\|_{2+\infty}
\less \|\nabla
\calS(t)  P_{g(\cdot,a(\infty))}^\perp\,f_2\|_{2} \\
& \less  \|\sqrt{H} \calS(t)
P_{g(\cdot,a(\infty))}^\perp\,f_2 \|_{2} + \| |V|^{\half}
\calS(t)  P_{g(\cdot,a(\infty))}^\perp\,f_2 \|_{2} \\
&\less \| \sin(t\sqrt{H}) P_{g(\cdot,a(\infty))}^\perp\, f_2 \|_2 +
\| |V|^{\half} (H^{-\half}\sin(t\sqrt{H})-c_0(\psi\otimes\psi))
P_{g(\cdot,a(\infty))}^\perp\,f_2 \|_{2}
\\
&\less \|f_2 \|_2
\end{align*}
The $\cos((t-s)\sqrt{H})$ term is treated basically in the same way
as in the dispersive estimate for $u$, so it will suffice to bound
the $\calS(t-s)$ integral. First, suppose that $0<t\less 1$. In that
case we use only $L^2$:
\begin{align*}
& \Big\| \nabla \il_0^t \calS(t-s) P_{g(\cdot,a(\infty))}^\perp
[(V(\cdot,a(\infty))-V(\cdot,a(s))) u(\cdot,s)+ N(u,\phi(\cdot,a(s)))]\, ds \Big\|_2 \\
&\less \Big\| \sqrt{H} \il_0^t \calS(t-s)
P_{g(\cdot,a(\infty))}^\perp
[(V(\cdot,a(\infty))-V(\cdot,a(s))) u(\cdot,s)+ N(u,\phi(\cdot,a(s)))]\, ds \Big\|_2 \\
& + \Big\| |V|^{\half} \il_0^t \calS(t-s)
P_{g(\cdot,a(\infty))}^\perp [(V(\cdot,a(\infty))-V(\cdot,a(s))) u(\cdot,s)+
N(u,\phi(\cdot,a(s)))]\, ds \Big\|_2
\end{align*}
This can be further
bounded by \begin{align*}
 &\less \Big\| \il_0^t \sin((t-s)\sqrt{H})
P_{g(\cdot,a(\infty))}^\perp
[(V(\cdot,a(\infty))-V(\cdot,a(s))) u(\cdot,s)+ N(u,\phi(\cdot,a(s)))]\, ds \Big\|_2 \\
& + \Big\| |V|^{\half} \il_0^t
(\sin((t-s)\sqrt{H})-c_0\psi\otimes\psi)
P_{g(\cdot,a(\infty))}^\perp [(V(\cdot,a(\infty))-V(\cdot,a(s)))
u(\cdot,s)+ N(u,\phi(\cdot,a(s)))]\, ds \Big\|_2 \\
&\less \delta^2
\end{align*}
Next, let $t\gg 1$. To bound the integral over $[0,t-t^{-10}]$ we
use the estimates
\[
\|D^2(V(\cdot,a(\infty))-V(\cdot,a(s))) u(\cdot,s) \|_{L^{1}(\R^3)} + \| D^2
N(u,\phi(\cdot,a(s))) \|_{L^{1}(\R^3)} \less \delta^2 \la
s\ra^{-3\eps/2}
\]
We discuss the contribution by $u^5$:
\[ \| |Du|^2 u^3(s) \|_1 + \| u^4 D^2 u(\cdot,s)\|_1  \less \|Du(\cdot,s)\|_{2+\infty}^2
\|u^3(s)\|_{1\cap \infty} + \| u^3(s) \|_2 \|D^2 u(\cdot,s)\|_2
\|u(\cdot,s)\|_\infty \less \delta^5 \la s\ra^{-3\eps/2}
\]
Strictly speaking, \eqref{eq:u_est} only gives $u^3\in L^p$ for
$p>1$. This, however, suffices since we can lower the $\infty$ in
$Du\in L^2+ L^\infty$ by interpolation (this explains the loss in
going from $2\eps$ to $3\eps/2$). Hence, we are dealing with the
integral
\[ \delta^2 \il_0^{t-t^{-10}} (t-s)^{-1} \la s\ra^{-3\eps/2}\, ds \less \delta^2
t^{-\eps}
\]
as desired. Finally, the contribution of $[t-t^{-10},t]$ is dealt
with in basically the same way as the case of small times. We skip
the details.

\noindent It remains to prove the important decay
estimate~\eqref{eq:u_est}. In view of~\eqref{eq:u_decay} it will
suffice to consider the case $|x|>A\la t\ra$ with $A$ large
depending on~$R$. We write the equation for $v$ in the form
\begin{align}
\Box v & = -V(\cdot,a(\infty))v - \pr_t (\dot{b}\phi_b(\cdot,b(t))) +
(V(\cdot,a(\infty)) - V(\cdot,a(t))) u + N(u,\phi(\cdot,a(t))) \label{eq:vrhs}\\
v(0) & = w_1, \qquad \pr_t v(0) = f_2 \nn
\end{align}
Then, with $S_0(t)=\f{\sin(t\sqrt{-\Delta})}{\sqrt{-\Delta}}$, we
have
\begin{align*}
& u(x,t) =
S_0(t) u_0(x) + \il_0^t S_0(t-s)
\big[-V(\cdot,a(\infty))v(\cdot,s) -
\pr_s (\dot{b}\phi_b(\cdot,b(s))) \\
&\qquad \qquad\qquad \qquad\qquad \qquad + (V(\cdot,a(\infty))-V(\cdot,a(s))
u(\cdot,s) + N(u,\phi(\cdot,a(s)))\big]
\, ds
\end{align*}
By our assumption on $x$, $S_0(t) u_0(x)=0$.  Due to the nature of
$S_0$ as an averaging operator, and our assumption on $x$,  we
obtain the bound
\begin{align*}
&|S_0(t-s) \big[-V(\cdot,a(s)) v(\cdot,s) +(V(\cdot,a(\infty))-V(\cdot,a(s))
u(\cdot,s) + N(u,\phi(\cdot,a(s)))\big]|(x) \\
&\less (t-s) (\delta  \la x\ra^{-4}\la s\ra^{-1} + \delta^2\la s\ra^{-2}\la x\ra^{-4} + C_1^5\, \la x\ra^{-5}
\delta^2 ) 
\end{align*}
Here we used the estimate $\|v(\cdot,s)\|_\infty \le \delta$, which we derived above
independently of the point-wise decay on $v$ being proven now.
Thus,
\begin{align*}
& \Big|\il_0^t S_0(t-s) \big[-V(\cdot,a(s)) v(\cdot,s) +(V(\cdot,a(\infty))-V(\cdot,a(s))
u(\cdot,s) +
N(u,\phi(\cdot,a(s)))\big] \,
ds\Big|(x) \\
& \less C_1\, \delta t^2 \la x\ra^{-4} + C_1^5\,\delta^2 t^2 \la
x\ra^{-5}\\
&\less C_1\, \delta t^2 (A\la t\ra)^{-3} \la x\ra^{-1} +
C_1^5\,\delta^2 t^2 (A\la t\ra)^{-4} \la x\ra^{-1}\\
&\ll C_1\,\delta \la x\ra^{-1}
\end{align*}
provided $A$ (and thus $C_1$) is large  and $\delta$ is small.
Furthermore, by the nature of $\cos((t-s)\sqrt{-\Delta})$,
\begin{align*}
& \Big| \il_0^t \dot{b}(s)\, \cos((t-s)\sqrt{-\Delta})\,
\phi_b(\cdot,b(s))\, ds \Big|(x) \less \il_0^t \delta \la
s\ra^{-2} \la x\ra^{-1}\, ds  \less \delta \la x\ra^{-1} \ll
C_1\,\delta \la x\ra^{-1}
\end{align*}
where we have again exploited that $|x|$ is large relative to $t$.
It is also important to note that the bound on $\dot{b}$ does not
deteriorate when $C_1$ becomes large, provided we make $\delta$ small.
Thus, the
conclusion is that
\[ |u(x,t)|\le C_1\,\delta\la x\ra^{-1} \]
as desired. The proposition is proved.
\end{proof}

\section{The contraction scheme: the fixed-point}
\label{sec:2}

We now show that the map $\Phi_{u_0}$, which we constructed in the
previous section, has a fixed-point. We will show that $\Phi_{u_0}$
contracts relative to the following distance:

\begin{defi}
\label{defi:norm} For any two points $\pone=(\uone,\aone),
\ptwo=(\utwo,\atwo)\in Y_{R,\delta}$ we define their distance to be
\begin{align}
\label{eq:weak_norm} d(\pone,\ptwo) &:= \sup_{t\ge0}
\big[\|D\uone(\cdot,t)-D\utwo(\cdot,t)\|_2 + \la
t\ra^{\eps}\|\uone(\cdot,t)-\utwo(\cdot,t)\|_{2+\infty} \big] \\
&\qquad\qquad\qquad\qquad +
\sup_{t\ge0} \la t\ra^{1+\eps} |\dot{a}^{(1)}(t)-\dot{a}^{(2)}(t)|\,
\nn
\end{align}
where $\eps>0$ is small and fixed (and no larger than the one
in~\eqref{eq:nablau_decay}).
\end{defi}

We record a simple technical fact which will need in the main
argument. Throughout this section, we will make use of the bounds
\eqref{eq:u_decay} -- \eqref{eq:u_est} without further mention.

\begin{lemma}
\label{lem:Vdiff} Fix $\pone=(\uone,\aone), \ptwo=(\utwo,\atwo)\in
Y_{R,\delta}$ and define
\[ W(\cdot,s) := [V(\cdot,\aone(\infty))-V(\cdot,\aone(s))]-
[V(\cdot,\atwo(\infty))-V(\cdot,\atwo(s))]
\]
where $V(x,a)=-5\phi(x,a)^4$.  Then
\begin{align}
& \sup_{s\ge0} |\aone(s)-\atwo(s)|  \less d(\pone,\ptwo) \label{eq:adiff}\\
& |W(x,s)| \less \la x\ra^{-4} \la s\ra^{-\eps} d(\pone,\ptwo)
\label{eq:West}
\end{align}
for all $x\in\R^3$, $s\ge0$.
\end{lemma}
\begin{proof}
First, for all $s\ge0$,
\[ |\aone(s)-\atwo(s)| \le \il_0^s \la
\sigma\ra^{1+\eps}
|\dot{a}^{(1)}(\sigma)-\dot{a}^{(2)}(\sigma)|\,\la
\sigma\ra^{-1-\eps}\, d\sigma \less  d(\pone,\ptwo)
\]
Second,
\begin{align*}
V(\cdot,\aone(s))- V(\cdot,\atwo(s)) &= \il_0^1 \f{d}{d\tau}
V(\cdot,\tau \aone(s) + (1-\tau) \atwo(s))\, d\tau \\
&= \il_0^1 \pr_a V(\cdot,\tau \aone(s) + (1-\tau) \atwo(s))\,
d\tau\; (\aone(s)-\atwo(s))
\end{align*}
Therefore,
\begin{align*}
W(\cdot,s) &= \il_0^1 [\pr_a V(\cdot,\tau \aone(\infty) + (1-\tau)
\atwo(\infty)) -\pr_a V(\cdot,\tau \aone(s) + (1-\tau) \atwo(s))] \,
d\tau\, (\aone(s)-\atwo(s)) \\
&\qquad + \il_0^1 \pr_a V(\cdot,\tau \aone(\infty) + (1-\tau)
\atwo(\infty))\, d\tau \; \il_s^\infty
(\dot{a}^{(1)}(\sigma)-\dot{a}^{(2)}(\sigma))\, d\sigma
\end{align*}
which implies that \begin{align*} |W(x,s)| &\less  \delta \la
x\ra^{-4} \la s\ra^{-1} |\aone(s)-\atwo(s)|   + \la x\ra^{-4}
\il_s^\infty
|\dot{a}^{(1)}(\sigma)-\dot{a}^{(2)}(\sigma)|\, d\sigma \\
&\less \la x\ra^{-4} \la s\ra^{-\eps} d(\pone,\ptwo)
\end{align*}
The lemma follows.
\end{proof}

The main result of this section is

\begin{prop}
\label{prop:fp} Let $u_0=(f_1,f_2)\in B_{\delta_0}(0)\subset \Sigma_0$ with
some $\delta_0\ll\delta$. Then there exists a unique fixed-point
$(u,a)\in Y_{R,\delta}$ of $\Phi_{u_0}$. Moreover, if we define
$ h(u_0):= h(u_0;u,a) $
with this choice of $(u,a)$, then
\[ |h(u_0)| \less \|u_0\|^2_{H^3\times H^2},\qquad
|h(u_0)-h(\til u_0)|\less \delta
\|u_0-\til u_0\|_{H^3\times H^2}
\] for all $u_0,\til u_0\in B_{\delta_0}(0)$.
\end{prop}
\begin{proof}
Let $\pone=(\uone,\aone)\in Y_{R,\delta}$ and
$\ptwo=(\utwo,\atwo)\in Y_{R,\delta}$. We set
\[ \qone=(\vone,\bone)=\Phi_{u_0}(\uone,\aone),\qquad
\qtwo=(\vtwo,\btwo)=\Phi_{u_0}(\utwo,\atwo) \] as well as
\[ \hone = h(u_0;\uone,\aone),\;
\qquad \htwo = h(u_0;\utwo,\atwo)\]
In other words, we have the following equations:
First, the equation for $\vj$
\begin{align}
& \pr_{tt} \vj + H(\aj(\infty)) \vj \nn \\
&= -\pr_{tt} \phi(\cdot,\bj(t)) +
(V(\cdot,\aj(\infty))-V(\cdot,\aj(t))) \uj  + N(\uj,\phi(\cdot,\aj(t))) \nn
\end{align}
Next, the equation for $\bj$
\begin{align}
0 & = 4\;(\aj(\infty))^{\f54} ((\bj)^{-\f14}(t)-1) - c_0\la f_2,\phi_a(\cdot,\aj(\infty))\ra\nn\\
&  - c_0 \il_0^t \la \phi_a(\cdot,\aj(\infty)),
(V(\cdot,\aj(\infty))-V(\cdot,\aj(s)))\uj(\cdot,s)+
N(\uj(\cdot,s),\phi(\cdot,\aj(s)))\ra\, ds  \nn\\
&  + c_0\big[ \la f_2,\phi_a(\cdot,\aj(\infty)) \ra - t \,\la \phi_a(\cdot,\aj(\infty)),
(V(\cdot,\aj(\infty))-V(\cdot,\aj(0)))\wj_1 + N(\wj_1,
\phi(\cdot,\aj(0)))\ra \big]\omega(t)  \nn
\end{align}
and finally the initial conditions:
\begin{align}
& \vj(0) = f_1 + \hj g_0 =: \wj_1 , \quad \pr_t \vj(0) =  f_2 \nn \\
& \bj(0) =1, \quad \dot{b}^{(j)}(0)=0 \nn
\end{align}
for $j=1,2$. In order to compare $\dot{b}^{(1)}(t)$ with
$\dot{b}^{(2)}(t)$, we need to following estimate
\begin{align*}
& \Big|\la \phi_a(\cdot,\aone(\infty)),
(V(\cdot,\aone(\infty))-V(\cdot,\aone(t)))\uone(\cdot,t)+
N(\uone(\cdot,t),\phi(\cdot,\aone(t)))\ra \\
&- \la \phi_a(\cdot,\atwo(\infty)),
(V(\cdot,\atwo(\infty))-V(\cdot,\atwo(t)))\utwo(\cdot,t)+
N(\utwo(\cdot,t),\phi(\cdot,\atwo(t)))\ra \Big| \\
&\less \Big|\la
\phi_a(\cdot,\aone(\infty))-\phi_a(\cdot,\atwo(\infty)),
(V(\cdot,\aone(\infty))-V(\cdot,\aone(t)))\uone(\cdot,t)+
N(\uone(\cdot,t),\phi(\cdot,\aone(t)))\ra \Big| \\
& \quad + \big| \la \phi_a(\cdot,\atwo(\infty)), W(\cdot,t)
\uone(\cdot,t) \ra \big| + \big| \la \phi_a(\cdot,\atwo(\infty)),
(V(\cdot,\atwo(\infty))-V(\cdot,\atwo(t)))(\uone(\cdot,t)-\utwo(\cdot,t))
\ra
\big| \\
&+ \big|  \la
\phi_a(\cdot,\atwo(\infty)),N(\uone(\cdot,t),\phi(\cdot,\aone(t))) -
N(\utwo(\cdot,t),\phi(\cdot,\atwo(t))) \ra\big| \\
&=: (A+B+C+D)(t)
\end{align*}
where $W$ is as in Lemma~\ref{lem:Vdiff}. By~\eqref{eq:adiff},
\[ |\phi_a(\cdot,\aone(\infty))-\phi_a(\cdot,\atwo(\infty))| \less
\la x\ra^{-1} d(\pone,\ptwo) \] Hence,
\[ A\less  \delta^2 \la
t\ra^{-2}d(\pone,\ptwo) 
\]
Furthermore, in view of \eqref{eq:West},
\[ B\less \delta \la t\ra^{-1-\eps} d(\pone,\ptwo)
\]
as well as
\[ C\less \delta \la t\ra^{-1-\eps} d(\pone,\ptwo) \]
Next, we have
\begin{align*}
 & \big|  \la
\phi_a(\cdot,\atwo(\infty)),(\uone(\cdot,t))^2\phi(\cdot,\aone(t))^3
-
(\utwo(\cdot,t))^2\phi(\cdot,\atwo(t))^3\ra\big| \\
&\less \big|  \la \phi_a(\cdot,\atwo(\infty)),((\uone(\cdot,t))^2 -
(\utwo(\cdot,t))^2)\phi(\cdot,\atwo(t))^3\ra\big| \\
& \qquad  + \big|  \la
\phi_a(\cdot,\atwo(\infty)),(\utwo(\cdot,t))^2
(\phi(\cdot,\aone(t))^3
- \phi(\cdot,\atwo(t))^3)\ra\big| \\
& \less \delta \la t\ra^{-1-\eps} d(\pone,\ptwo)
\end{align*}
and for the quintic term
\begin{align*}
 & \big|  \la
\phi_a(\cdot,\atwo(\infty)),(\uone(\cdot,t))^5 -
(\utwo(\cdot,t))^5\ra\big| \\
&\less  \delta^5 \la t\ra^{-1-\eps} d(\pone,\ptwo)
\end{align*}
In conclusion,
\[ (A+B+C+D)(t)\less \delta \la
t\ra^{-1-\eps} d(\pone,\ptwo)\] From this and the bound
\[
|\la f_2,\phi_a(\cdot,\aone(\infty))\ra - \la f_2,\phi_a(\cdot,\atwo(\infty))\ra | \less \delta\, d(\pone,\ptwo)
\]
we infer  that
\[ |\bone(\infty)-\btwo(\infty)|\less \delta d(\pone,\ptwo)\]
as well as
\[ \sup_{t\ge0} |\bone(t)-\btwo(t)|\less \delta d(\pone,\ptwo)\]
The latter requires some care, as on the support of $\omega$ it also involved the difference
$\hone-\htwo$. However, if we borrow the estimate~\eqref{eq:hdiffest} for now, then we obtain
our desired result.
Using these bounds we obtain furthermore that
\begin{equation}
\label{eq:bdiffest}  \sup_{t\ge0}\la t\ra^{1+\eps} |\dot{b}^{(1)}(t)
-\dot{b}^{(2)}(t)|\, dt \less \delta d(\pone,\ptwo)
\end{equation}
Next, we turn to the initial conditions. In view of~\eqref{eq:n+(0)},
\begin{align}
 &(2\kj)^{-\half} \big[\la
g(\cdot,\aj(\infty)), \kj\wj_1+f_2 \ra+\kj \hj \la g_0, g(\cdot,\aj(\infty))\ra\big] \nn \\
& = \kj\il_0^\infty \dot{b}^{(j)}(t)\,
e^{-t\kj} \la \phi_b(\cdot,\bj(t)),
g(\cdot,\aj(\infty)) \ra\, dt \nn \\
&\qquad  - \il_0^\infty e^{-t\kj} \la
(V(\cdot,\aj(\infty))-V(\cdot,\aj(t))) \uj(\cdot,t) +
N(\uj(\cdot,t),\phi(\cdot,\aj(t))), g(\cdot,\aj(\infty)) \big\ra \,
dt \label{eq:S2}
\end{align}
with $j=1,2$.
It is important to note that $\hj$ also appears on the right-hand side of~\eqref{eq:S2}
because of the dependence of $\dot{b}^{(j)}$ on~$\hj$. However, as remarked in the proof
of Proposition~\ref{prop:stab}, this dependence occurs with a small factor of $\delta$,
or to a higher order. Hence, we can solve for $\hj$ via the implicit function theorem, and then
estimate the difference. Alternatively, we can solve for the difference and then estimate it.
Either way, we obtain that
\begin{equation}
|\hone-\htwo|\less \delta d(\pone,\ptwo)
\label{eq:hdiffest}
\end{equation}
It follows from Lemma~\ref{lem:ODE_stable} and \eqref{eq:S2} that
\begin{align*} n^{(j)}_+(t) &
= \kj\il_t^\infty \dot{b}^{(j)}(s)\, e^{-(s-t)\kj} \la
\phi_b(\cdot,\bj(s)), g(\cdot,\aj(\infty)) \ra\, ds + \kj
\dot{b}^{(j)}(t)\, \la \phi_b(\cdot,\bj(t)),
g(\cdot,\aj(\infty)) \ra \nn \\
&- \il_t^\infty e^{-(s-t)\kj} \la
(V(\cdot,\aj(\infty))-V(\cdot,\aj(s))) \uj(\cdot,s) +
N(\uj(\cdot,s),\phi(\cdot,\aj(s))), g(\cdot,\aj(\infty)) \big\ra \,
ds
\end{align*}
Note that $\kj=\sqrt{\aj(\infty)}\, k_0$ and thus
\[ |k^{(1)}-k^{(2)}|\less d(\pone,\ptwo),\qquad |e^{-tk^{(1)}}-e^{-tk^{(2)}}|\less d(\pone,\ptwo)\,e^{-tk^{(1)}/2}
\qquad \forall \;t\ge0
\]
By this and the estimates which we have just derived,
\[ \la t\ra^{1+\eps} |n^{(1)}_+(t) - n^{(2)}_+(t)| \less \delta
d(\pone,\ptwo) \] In a similar fashion, we derive
\[ \la t\ra^{1+\eps} |n^{(1)}_-(t) - n^{(2)}_-(t)| \less \delta
d(\pone,\ptwo) \] from the representations
\begin{align*}
& n^{(j)}_-(t)  = e^{-t\kj}\;n^{(j)}_-(0) \\
& +  \kj\il_0^t \dot{b}^{(j)}(s)\, e^{-(t-s)\kj} \la
\phi_b(\cdot,\bj(s)), g(\cdot,\aj(\infty)) \ra\, ds + \kj
\dot{b}^{(j)}(t)\, \la \phi_b(\cdot,\bj(t)),
g(\cdot,\aj(\infty)) \ra \nn \\
&- \il_0^t e^{-(t-s)\kj} \la (V(\cdot,\aj(\infty))-V(\cdot,\aj(s)))
\uj(\cdot,s) + N(\uj(\cdot,s),\phi(\cdot,\aj(s))),
g(\cdot,\aj(\infty)) \big\ra \, ds
\end{align*}
and the fact that
\begin{align*}
|n^{(1)}_-(0)- n^{(2)}_-(0)| &\less  (2k^{(1)})^{-\half} |\la
g(\cdot,\aone(\infty)), \wone_1 \ra-\la g(\cdot,\atwo(\infty)),
\wtwo_1 \ra| \\
&\quad+ |k^{(1)}-k^{(2)}| |\la g(\cdot,\atwo(\infty)),
\wtwo_1 \ra|+ |\none_+(0)-\ntwo_+(0)|
\less \delta d(\pone,\ptwo)
\end{align*}
In view of these bounds, estimate~\eqref{eq:npm_decay}, and the
fact that
\[ \vj(\cdot,t) = (2\kj)^{-\half} (\nj_+(t)+\nj_-(t))g(\cdot,\aj(\infty)) + \tilvj(\cdot,t)
\]
it will suffice to estimate the difference of the $\tilvj$ which are
given by
\begin{align} \tilvj(t)
&= \cos(t\sqrt{H_j}) P_{g(\cdot,\aj(\infty))}^\perp\,\wj_1+\calS_j(t)  P_{g(\cdot,\aj(\infty))}^\perp\, f_2 \nn\\
&\quad - \il_0^t \dot{b}^{(j)}(s)\, \cos((t-s)\sqrt{H_j})
 P_{g(\cdot,\aj(\infty))}^\perp\;
[\phi_b(\cdot,\bj(s))-(\aj(\infty)/\bj(s))^{\f54}\phi_a(\cdot,\aj(\infty))]\, ds \label{eq:bj_disp}\\
&\quad  - \il_0^t \calS_j(t-s)\, P_{g(\cdot,\aj(\infty))}^\perp
[(V(\cdot,\aj(\infty))-V(\cdot,\aj(s))) \uj(\cdot,s)+
N(\uj,\phi(\cdot,\aj(s)))]\, ds \label{eq:vj_disp} \\
& \quad + c_0\big[ \la f_2,\psi \ra - t \,\la \psi, (V(\cdot,\aj(\infty))-V(\cdot,\aj(0)))\wj_1 + N(\wj_1,
\phi(\cdot,\aj(0)))\ra \big]\omega(t)\psi \label{eq:neues}
\end{align}
Here, $H_j:=-\Delta - 5\phi^4(\cdot,\aj(\infty))$ and
$\calS_j(t) := \f{\sin(t\sqrt{H_j})}{\sqrt{H_j}}P_{g(\cdot,\aj(\infty))}^\perp-c_0(\psi_j\otimes\psi_j).$

\noindent By the bounds \eqref{eq:hdiffest},
\[ \|\wone_1 - \wtwo_1\|_{H^3} \less \delta d(\pone,\ptwo) \]
and also
\[ \|P_{g(\cdot,\aone(\infty))}^\perp\, \wone_1 - P_{g(\cdot,\atwo(\infty))}^\perp\, \wtwo_1
 \|_{H^3} \less
 \delta d(\pone,\ptwo) \]
By the estimates from the proof of Proposition~\ref{prop:stab} we
conclude that \begin{align*}
\|\cos(t\sqrt{H_1}) [P_{g(\cdot,\aone(\infty))}^\perp\, \wone_1 -
P_{g(\cdot,\atwo(\infty))}^\perp\, \wtwo_1]\|_\infty &\less \delta
\la t\ra^{-1} d(\pone,\ptwo) \\
\|\nabla
\cos(t\sqrt{H_1})[P_{g(\cdot,\aone(\infty))}^\perp\, \wone_1 -
P_{g(\cdot,\atwo(\infty))}^\perp\, \wtwo_1]\|_2 &\less \delta
d(\pone,\ptwo)\\
 \|
\calS_1(t)[P_{g(\cdot,\aone(\infty))}^\perp\, f_2 -
P_{g(\cdot,\atwo(\infty))}^\perp\, f_2]\|_\infty &\less \delta
\la t\ra^{-1} d(\pone,\ptwo) \\
 \|\nabla
\calS_1(t)[P_{g(\cdot,\aone(\infty))}^\perp\, f_2 -
P_{g(\cdot,\atwo(\infty))}^\perp\, f_2]\|_2 &\less \delta
d(\pone,\ptwo)
\end{align*}
We also need to consider terms which involve the difference of the evolutions:
\[
 \cos(t\sqrt{H_1})P_{g(\cdot,\aone(\infty))}^\perp -\cos(t\sqrt{H_2})P_{g(\cdot,\atwo(\infty))}^\perp
 \text{\ \ and\ \ }
 \calS_1(t) - \calS_2(t)
\]
However, these operators lead to the desired bounds because of Corollary~\ref{cor:V_stab}.
The difference of the integrals in~\eqref{eq:bj_disp} is of the form
\begin{align}
&\il_0^t (\dot{b}^{(1)}(s)-\dot{b}^{(2)}(s))\, \cos((t-s)\sqrt{H_1})
 P_{g(\cdot,\aone(\infty))}^\perp\;
[\phi_b(\cdot,\bone(s))-(\aone(\infty)/\bone(s))^{\f54}\phi_a(\cdot,\aone(\infty))]\,
ds \label{eq:int1_diff}\\
&+\il_0^t \dot{b}^{(2)}(s)\, \cos((t-s)\sqrt{H_1}) P_{g(\cdot,\aone(\infty))}^\perp\;\Big\{
[\phi_b(\cdot,\bone(s))-(\aone(\infty)/\bone(s))^{\f54}\phi_a(\cdot,\aone(\infty))]\nn
\\
& \quad -
[\phi_b(\cdot,\btwo(s))-(\atwo(\infty)/\btwo(s))^{\f54}\phi_a(\cdot,\atwo(\infty))]
\Big\} \, ds \label{eq:int2_diff} \\
&\quad +\il_0^t \dot{b}^{(2)}(s)\,
 [\cos((t-s)\sqrt{H_1})P_{g(\cdot,\aone(\infty))}^\perp-\cos((t-s)\sqrt{H_2})P_{g(\cdot,\atwo(\infty))}^\perp]\cdot\nn \\
& \qquad\qquad \qquad\qquad\qquad\qquad\qquad\qquad\qquad\qquad \cdot[\phi_b(\cdot,\btwo(s))-(\atwo(\infty)/\btwo(s))^{\f54}\phi_a(\cdot,\atwo(\infty))]
 \, ds \label{eq:vergessen}
\end{align}
First,
\begin{align*}
 \| \eqref{eq:int1_diff}\|_\infty &\less \delta\,
d(\pone,\ptwo) \il_0^t \la s\ra^{-1-\eps} \la t-s\ra^{-1} \la
s\ra^{-1} \, ds \less \delta\, \la t\ra^{-1}\, d(\pone,\ptwo) \\
\| \eqref{eq:int2_diff}\|_\infty &\less \il_0^t \delta^2\, \la
s\ra^{-2} \la t-s\ra^{-1} \delta\, d(\pone,\ptwo) \, ds \less \delta
\la t\ra^{-1} \, d(\pone,\ptwo)
\end{align*}
and second,
\begin{align*}
 \| \nabla \eqref{eq:int1_diff}\|_2 &\less \delta\,
d(\pone,\ptwo) \il_0^t \la s\ra^{-1-\eps}  \la
s\ra^{-1} \, ds \less \delta \, d(\pone,\ptwo) \\
\| \nabla\eqref{eq:int2_diff}\|_2 &\less \il_0^t \delta^2 \la
s\ra^{-2}  \delta d(\pone,\ptwo) \, ds \less \delta \,
d(\pone,\ptwo)
\end{align*}
As far as the term \eqref{eq:vergessen} is concerned, we remark that
it, too, satisfies the desired bounds due to the stability result in Section~\ref{sec:pointwise},
see Corollary~\ref{cor:V_stab}.
Finally, we turn to the difference
\begin{align*}
E(t) &:=\il_0^t \calS_1(t-s)\, P_{g(\cdot,\aone(\infty))}^\perp
[(V(\cdot,\aone(\infty))-V(\cdot,\aone(s))) \uone(\cdot,s)+
N(\uone,\phi(\cdot,\aone(s)))]\, ds \\
& - \il_0^t \calS_2(t-s)\, P_{g(\cdot,\atwo(\infty))}^\perp
[(V(\cdot,\atwo(\infty))-V(\cdot,\atwo(s))) \utwo(\cdot,s)+
N(\utwo,\phi(\cdot,\atwo(s)))]\, ds
\end{align*}
By the same type of arguments which we have used repeatedly up to
this point the reader will check that
\[ \|\nabla E(t)\|_2 + \la t\ra^{\eps} \|E(t)\|_{2+\infty} \less
\delta\, d(\pone,\ptwo)
\]
for all $t\ge0$. This concludes the proof of the estimate
\[ d(\Phi_{u_0}(\pone),\Phi_{u_0}(\ptwo)) \less \delta\,
d(\pone,\ptwo) \] and therefore of the existence of a fixed-point
\[ (u,a)(u_0)\in Y_{R,\delta} \]
Since $\Phi_{u_0}$ is Lipschitz in $u_0$ by
Proposition~\ref{prop:stab}, we conclude that the fixed-point is
also Lipschitz in $u_0$, see Lemma~\ref{lem:fp} below. Let $u_0,\til
u_0\in B_\delta(0)$ and denote their fixed-points by $(u,a), (\til
u,\til a)$. Then, by \eqref{eq:hetadiff} and~\eqref{eq:hdiffest},
\begin{align*}
|h(u_0;u,a)-h(\til u_0;\til u,\til a)| & \less |h(u_0;u,a)-h(\til
u_0; u, a)| + |h(\til u_0;u,a)-h(\til u_0;\til u,\til a)|| \\
& \less \delta \|u_0-\til u_0\|_{H^2} + \delta d(\pone,\ptwo) \\
&\less \delta \|u_0-\til u_0\|_{H^2}
\end{align*}
and we are done.
\end{proof}

The following lemma is completely standard, we present it for the
sake of completeness.

\begin{lemma}
\label{lem:fp} Let $S$ be a complete metric space  and $T$ an
arbitrary metric space. Suppose that $A:S\times T \to S$ so that
with some $0<\gamma<1$
\begin{align*}
 \sup_{t\in T}\;d_X(A(x,t),A(y,t)) &\le \gamma\; d_X(x,y) \text{\ \ for all \ \ }x,y\in S,\\
  \sup_{x\in S}\; d_X(A(x,t_1),A(x,t_2)) &\le C_0 \;d_Y(t_1,t_2) \text{\ \ for all \ \ }t_1,t_2\in T.
\end{align*}
Then for every $t\in T$ there exists a unique fixed-point $x(t)\in
S$ such that $A(x(t),t)=x(t)$. Moreover, these points satisfy the
bounds
\[ d_X(x(t_1),x(t_2)) \le \frac{C_0}{1-\gamma}\;d_Y(t_1,t_2)\]
for all $t_1,t_2\in T$.
\end{lemma}
\begin{proof}
Clearly, $x(t)= \lim_{n\to\infty} A(x_n(t),t)$ where for some fixed
(i.e., independent of $t$) $x_0$
\[ x_0(t):= x_0, \quad x_{n+1}(t)= A(x_n(t),t).\]
Then inductively,
\begin{align*}
d_X(x_{n+1}(t_1),x_{n+1}(t_2)) & \le d_X( A(x_n(t_1),t_1),
A(x_n(t_2),t_1)) + d_X(A(x_n(t_2),t_1), A(x_n(t_2),t_2)) \\
&\le \gamma d_X(x_n(t_1),x_n(t_2)) + C_0 d_Y(t_1,t_2) \\
&\le C_0 \sum_{k=0}^n\gamma^{k}\;d_Y(t_1,t_2)
\end{align*}
for all $n\ge0$. Passing to the limit $n\to\infty$ proves the lemma.
\end{proof}

It is now easy to prove Theorem~\ref{thm:main}.

\begin{proof}[Proof of Theorem~\ref{thm:main}:]
Let $u_0=(f_1,f_2)\in B_\delta(0)\subset \Sigma_0$.
By Proposition~\ref{prop:fp} there exists a fixed-point $(u,a)\in Y_{R,\delta}$
of the map $\Phi_{u_0}$. By construction, this means that there exists $h(u_0)$
as in Proposition~\ref{prop:fp} so that the modified initial data~\eqref{eq:mod_init}
lead to a global solution of~\eqref{eq:v_eq}. I.e.,
\[ \pr_{tt} u + H(a(\infty)) u = -\pr_{tt} \phi(\cdot,a(t)) + (V(\cdot,a(\infty))-V(\cdot,a(t)))u + N(u,\phi(\cdot,a(t))) \]
which  is the same as
\[ \pr_{tt} u + H(a(t)) u = -\pr_{tt} \phi(\cdot,a(t))  + N(u,\phi(\cdot,a(t))) \]
This in turn implies that
\[\psi(\cdot,t)=\phi(\cdot,a(t))+u(\cdot,t)= \phi(\cdot,a(\infty))+v(\cdot,t) \]
with $v(\cdot,t) := \phi(\cdot,a(t)) - \phi(\cdot,a(\infty)) + u(\cdot,t) $
solves
\[ \Box \psi - \psi^5 =0\]
with initial conditions \eqref{eq:mod_init}. Finally, $a(0)=1$, $\dot{a}(0)=0$ by construction
and $h(u_0)$ and $u$ satisfy the bounds from Propositions~\ref{prop:stab} and~\ref{prop:fp}.
Therefore, we also have
\[ \|v(\cdot,t)\|_\infty \less \delta \la t\ra^{-1} \]
To derive the scattering statement, we write the equation for~$u(x,t)$
as a Hamiltonian systems with Hamiltonian~$J\Hil_\infty$, see~\eqref{eq:sys}.
 Thus, set $U=\binom{u}{\pr_t u}$ and write
\[ U(\cdot,t) = n_+(t) G_+(\cdot) + n_-(t) G_-(\cdot) + \til U(\cdot,t) \]
where $\til U= P_e U$, see \eqref{eq:U_decomp}.
Our goal is to find initial data $(\til f_1,\til f_2)\in \dot H^1\times L^2$  so that
\begin{equation}
\label{eq:wave_op}
  U(t) = U_0(t) + \binom{0}{-\dot{a}(t)\phi_a(\cdot,a(t))} + o_{\calE}(1)  \text{\ \ as\  \ }t\to\infty
\end{equation}
where $U_0$ is the solution vector of the free wave equation with data $U_0(0)=\binom{\til f_1}{\til f_2}$.
Here $\calE=\dot{H}^1\times L^2$ refers to the energy space
with norm
\[ \Big\| \binom{u_1}{u_2} \Big\|_{\calE}^2 :=  \|\nabla\; u_1\|_2^2 + \|u_2\|_2^2  \]
Notice that \eqref{eq:wave_op} yields the scattering claim of Theorem~\ref{thm:main} simply because
\[ v(\cdot,t) = u(\cdot,t) + \phi(\cdot,a(t))-\phi(\cdot,a(\infty)) \]
satisfies, see \eqref{eq:wave_op},
\begin{align*}
 \binom{v}{\pr_t v}(t) &= \binom{u}{\pr_t u}(t) + \binom{\phi(\cdot,a(t))-\phi(\cdot,a(\infty))}{\dot a(t)\phi_a(\cdot,a(t))} \\
 &= U_0(t) + \binom{\phi(\cdot,a(t))-\phi(\cdot,a(\infty))}{0} + o_{\calE}(1) \\
 & = U_0(t) + o_{\calE}(1) \text{\ \ as\  \ }t\to\infty
\end{align*}
Since we have shown in Proposition~\ref{prop:stab}
that the coefficients $n_{\pm}(t)$ decay like $\la t\ra^{-2}$, it will suffice to prove~\eqref{eq:wave_op}
 with $\til U=P_e U$ instead of $U$. Thus, we need to find initial data $(\til f_1,\til f_2)\in \calE$ so that
\begin{equation}
\label{eq:wave_op'}
 \| \til U(t)  - (0,-\dot{a}(t)\phi_a(\cdot,a(t)))^\dagger - U_0(t)\|_{\calE} \to 0 \text{\ \ as\  \ }t\to\infty
\end{equation}
Here $\dagger$ means transposition.
This will be done in two steps. First, we will find initial data $(f_1',f_2')\in \calE$ so that
\begin{equation}
\label{eq:wave_op''}
 \| \til U(t)  - (0,-\dot{a}(t)\phi_a(\cdot,a(t)))^\dagger - e^{tJ\Hil_\infty}P_e\,(f_1',f_2')^{\dagger} \|_{\calE} \to 0 \text{\ \ as\  \ }t\to\infty
\end{equation}
see \eqref{eq:free_prop}.
Because of the dispersive estimates on $\til U(t)$ and
$e^{tJ\Hil_\infty}P_e$, and in view of~\eqref{eq:grad1}, \eqref{eq:wave_op''} will follow from
\begin{equation}
\label{eq:wave_op'''}
 \| \til U(t)  - (0,-\dot{a}(t)\phi_a(\cdot,a(t)))^\dagger - e^{tJ\Hil_\infty}P_e\,( f_1', f_2')^{\dagger} \|_{\calE_\infty} \to 0 \text{\ \ as\  \ }t\to\infty
\end{equation}
where
\[ \Big\| \binom{u_1}{u_2} \Big\|_{\calE_\infty}^2 := \|\sqrt{H_\infty}\, u_1\|_2^2 +  \|u_2\|_2^2 \]
for all $(u_1,u_2)\in \,P_{ g(\cdot,a(\infty))}^\perp\,[\dot H^1\times L^2]$.
Strictly speaking, we have only derived dispersive estimates on $\til u$, and not on~$\pr_t \til u$.
However, it follows from the explicit form of $\til u$, see~\eqref{eq:nondisp1}, that for large times
\begin{align*}
\pr_t \til u(\cdot,t) &= -\sin(t\sqrt{H_\infty})P_{g_\infty}^\perp \sqrt{H_\infty}\; w_1
+ \cos(t\sqrt{H_\infty}) P_{g_\infty}^\perp w_2
 - \dot{a}(t) P_{g_\infty}^\perp [\phi_a(\cdot,a(t)) - (a(\infty)/a(t))^{\f54} \psi] \\
& + \il_0^t \dot{a}(s) \sin((t-s)\sqrt{H_\infty}) P_{g_\infty}^\perp\,\sqrt{H_\infty} [\phi_a(\cdot,a(s)) - (a(\infty)/a(s))^{\f54} \psi ]\, ds
\\
& - \psi(\cdot)\,\la\psi , (V(\cdot,a(\infty))-V(\cdot,a(t)))u(\cdot,t) + N(u,\phi(\cdot,a(t))) \ra \\
& - \il_0^t \cos((t-s)\sqrt{H_\infty}) P_{g_\infty}^\perp
[(V(\cdot,a(\infty))-V(\cdot,a(t)))u(\cdot,t) + N(u,\phi(\cdot,a(t)))]  \, ds
\end{align*}
where $\psi(\cdot)=\phi_a(\cdot,a(\infty))$ and $g_\infty = g(\cdot,a(\infty))$ as usual. By the estimates which we have derived we easily conclude that
\[ \| \pr_t \til u(\cdot,t)\|_\infty \less \delta \la t\ra^{-1} \]
Thus $\til U$ is dispersive as claimed and \eqref{eq:wave_op''} reduces to~\eqref{eq:wave_op'''}.
Now we remark that the group $e^{tJ\Hil_\infty}$ is unitary on $P_{ g(\cdot,a(\infty))}^\perp\,[\dot H^1\times L^2]$
relative to the norm $\calE_\infty$. Hence, \eqref{eq:wave_op'''} is the same as showing that
\begin{equation}
\label{eq:wave_op''''}
 \| e^{-tJ\Hil_\infty}[\til U(t) - (0,-\dot{a}(t)\phi_a(\cdot,a(t)))^\dagger] - P_e\,( f_1',f_2')^{\dagger} \|_{\calE_\infty} \to 0 \text{\ \ as\  \ }t\to\infty
\end{equation}
Note that
\[ \|\dot{a}(t)[\phi_a(\cdot,a(t))-(a(\infty)/a(t))^{\f54}\phi_a(\cdot,a(\infty))] \|_2 \less \delta^4 t^{-3} \]
as $t\to\infty$ since the term in brackets decays like $\la x\ra^{-3}$.
Thus, in view of the unitarity of $e^{-tJ\Hil_\infty}$, it follows that~\eqref{eq:wave_op''''}
is equivalent with
\begin{equation}
\label{eq:wave_opv}
 \| e^{-tJ\Hil_\infty}[\til U(t) - (a(\infty)/a(t))^{\f54}(0,-\dot{a}(t)\phi_a(\cdot,a(\infty)))^\dagger] - P_e\,( f_1',f_2')^{\dagger} \|_{\calE_\infty} \to 0 \text{\ \ as\  \ }t\to\infty
\end{equation}
First, we note that by \eqref{eq:free_prop}
\begin{align*}
& e^{-tJ\Hil_\infty}\binom{0}{-\dot{a}(t)(a(\infty)/a(t))^{\f54}\psi}
= \dot{a}(t)(a(\infty)/a(t))^{\f54}\binom{t\psi}{-\psi} \\
\end{align*}
Hence, \eqref{eq:wave_opv} is equivalent with
\begin{equation}
\label{eq:wave_opvi}
 \| e^{-tJ\Hil_\infty}\til U(t) - \dot{a}(t)(a(\infty)/a(t))^{\f54}(t\psi,-\psi)^\dagger
  - P_e\,( f_1',f_2')^{\dagger} \|_{\calE_\infty} \to 0 \text{\ \ as\  \ }t\to\infty
\end{equation}
Second, by \eqref{eq:sys_duh},
\[ e^{-tJ\Hil_\infty}\til U(t) = P_e \til U(0) + \il_0^t e^{-sJ\Hil_\infty} P_e W(s)\, ds \]
with
\begin{align*}
 W(s) &= \binom{0}{-\pr_{s}(\dot{a}(s)\phi_a(\cdot,a(s))) +(V(\cdot,a(\infty))-V(\cdot,a(s))) u+
N(u,\phi(\cdot,a(s)))} \\
&= \binom{0}{-\pr_{s}(\dot{a}(s)\phi_a(\cdot,a(s)))} + \tilde W(s)
\end{align*}
see \eqref{eq:sys} (we define $\til W$ by the second line). Integrating by parts and using \eqref{eq:free_prop} yields
\begin{align}
&\il_{0}^{t} e^{-sJ\Hil_\infty} P_e \binom{0}{-\pr_{s}(\dot{a}(s)\phi_a(\cdot,a(s)))}\, ds \nn \\
&=
\dot{a}(t)(a(\infty)/a(t))^{\f54} \binom{t\psi}{-\psi} + \dot{a}(t)\binom{\f{\sin(t\sqrt{H_\infty})}{\sqrt{H_\infty}}
P_{g_\infty}^\perp [\phi_a(\cdot,a(t))-(a(\infty)/a(t))^{\f54}\psi]}{-\cos(t\sqrt{H_\infty})
P_{g_\infty}^\perp [\phi_a(\cdot,a(t))-(a(\infty)/a(t))^{\f54}\psi]} \label{eq:corr1}\\
& \quad -\il_0^t \dot{a}(s) \binom{\cos(s\sqrt{H_\infty})
P_{g_\infty}^\perp \phi_a(\cdot,a(s)) }{\f{\sin(t\sqrt{H_\infty})}{\sqrt{H_\infty}}
P_{g_\infty}^\perp H_\infty \phi_a(\cdot,a(s))} \, ds
\nn
\end{align}
Observe that the first term in \eqref{eq:corr1} is identical with the middle term in~\eqref{eq:wave_opvi}.
The second term in \eqref{eq:corr1} is $o(1)$ in the energy norm $\calE_\infty$ as $t\to\infty$
and can hence be ignored for the purposes of~\eqref{eq:wave_opvi}.
We define
\begin{equation}
\nn \binom{f_1'}{f_2'} := P_e \til U(0) -\il_0^\infty \dot{a}(s) \binom{\cos(s\sqrt{H_\infty})
P_{g_\infty}^\perp \phi_a(\cdot,a(s)) }{\f{\sin(t\sqrt{H_\infty})}{\sqrt{H_\infty}}
P_{g_\infty}^\perp H_\infty \phi_a(\cdot,a(s))} \, ds
+ \il_0^\infty e^{-sJ\Hil_\infty} P_e \til W(s)\, ds
\end{equation}
This definition is justified, since the integrals are absolutely convergent in the norm of~$\calE_\infty$.
In addition, $ \binom{f_1'}{f_2'}\in \calE_\infty\cap L^\infty\subset \calE$.
Having established~\eqref{eq:wave_opvi} and therefore~\eqref{eq:wave_op''}, we now carry out the second step.
It consists of finding initial data $(\til f_1,\til f_2)\in \calE$ so that
\begin{equation}
\label{eq:wave_op_v}
 \| e^{tJ\Hil_\infty}P_e\,(f_1',f_2')^{\dagger} - e^{tJ\Hil_{\rm free}}\,(\til f_1,\til f_2)^{\dagger} \|_{\calE} \to 0 \text{\ \ as\  \ }t\to\infty
\end{equation}
However, this latter property is a consequence of the asymptotic completeness of $-\Delta+V$,
which is standard. The theorem is proved.
\end{proof}

\section{Linear theory: Point-wise decay}
\label{sec:pointwise}

In this section $H=-\Lap+V$ where $V$ is real-valued and decays
faster than a third power: $|V(x)|\less \la x\ra^{-\kappa}$ with
$\kappa>3$. Although we are of course only interested in the special
potential $V$ from the previous sections, we will keep this
discussion more general.  We emphasize that we work on all of
$L^2(\R^3)$ here  and assume that $H$ has no eigenvalue at zero. In
the case of $H$ as above this is false, but as explained there, it
is true when $H$ is restricted to the radial functions. Hence, in
order to apply the results from this section we need to restrict $H$
to the invariant subspace of radial functions.

\subsection{The sine evolution}

We study the evolution \beeq \label{eq:evol}
\frac{\sin(t\sqrt{H})}{\sqrt{H}} P_c = \il_0^\infty
\frac{\sin(t\sqrt{\lambda})}{\sqrt{\lambda}} E(d\lambda) \eneq
 where
$E$ is the spectral resolution of $H$ and $P_c=\chi_{[0,\infty)}(H)$
the projection onto the continuous spectrum. It arises as solution
of the Cauchy problem
\[ (\partial_{tt}+H)u=0, \quad u(0)=0, \pr_t u(0)=f \]
In this section our goal is to prove
Proposition~\ref{prop:infty_decay}. The question of dispersive decay
for the wave equation with a potential has received much attention
in recent years, see the papers by Beals, Strauss, Cuccagna,
Georgiev, Viscilia, Yajima,  d'Ancona, Pierfelice in the references.
However, none of these refences apply here since they either assume
that $V\ge0$, $V$ small, or that zero is neither an eigenvalue nor a
resonance.

\begin{prop}
\label{prop:infty_decay} Assume that $V$ is a real-valued potential
such that  $|V(x)|\less \la x\ra^{-\kappa}$ for some $\kappa>3$.
If $H$ has neither a resonance nor an eigenvalue at zero, then
\[ \Big\| \f{\sin(t\sqrt{H})}{\sqrt{H}} P_c f \Big\|_\infty \less t^{-1} \|f\|_{W^{1,1}(\R^3)}
\]
for all $t>0$.
Now assume that zero is a resonance  but not an
eigenvalue of $H=-\Delta+V$. Let $\psi$ be the unique resonance
function normalized so that $\int V\psi(x)\,dx=1$. Then there exists
a constant $c_0\ne0$ such that
\begin{equation}
\label{eq:1inf} \Big\| \f{\sin(t\sqrt{H})}{\sqrt{H}} P_c f - c_0
(\psi\otimes \psi) f\Big\|_\infty \less t^{-1} \|f\|_{W^{1,1}(\R^3)}
\end{equation}
for all $t>0$.
\end{prop}

We will only prove the second part which is harder. The case when zero is
neither a resonance nor an eigenvalue is implicit in our proof below and we will
henceforth assume that we are in the second case.
In~\eqref{eq:1inf}
\[ (\psi\otimes \psi) f (x) = \psi(x) \int f(y)\psi(y)\, dy \]
which is well-defined for all $f\in L^1$ since $\psi\in
L^\infty(\R^3)$. Indeed, one has $\psi+(-\Delta)^{-1}V\psi=0$ which
is the same as
\beeq\label{eq:res_form} \psi(x)=-\il_{\R^3} \frac{V(y)\psi(y)}{4\pi|x-y|}\, dy
\eneq
Hence, $\psi$ is bounded provided we can show that $V\psi\in
L^{\f32-\eps}\cap L^{\f32+\eps}(\R^3)$ for some $\eps>0$. This,
however, follows from the decay of $V$ and
\[ \|V\psi\|_{\f32\pm} \less \|\la x\ra^{\sigma}V\|_{6\pm} \|\la
x\ra^{-\sigma} \psi\|_2 \] where $\sigma=\f12+$.
We also remark that~\eqref{eq:res_form} implies that well-known fact that
\beeq \label{eq:intne0}
\il_{\R^3} V\psi(x)\, dx \ne0
\eneq
Indeed, if this vanished, then we could write
\[ \psi(x)=-\il_{\R^3} \Big(\frac{V(y)\psi(y)}{4\pi|x-y|}-\frac{V(y)\psi(y)}{4\pi|x|}\Big)\, dy
\]
which would imply that $|\psi(x)|\less \la x\ra^{-2}$ in contradiction to $\psi\not\in L^2(\R^3)$.


We now start with the detailed argument for the proposition.
The evolution~\eqref{eq:evol} can be written as
\beeq\label{eq:sintrans} \f{1}{i\pi}\il_0^\infty
\frac{\sin(t\lambda)}{\lambda}
[R_V^+(\lambda^2)-R_V^{-}(\lambda^2)]\,\lambda d\lambda =
\f{1}{i\pi} \il_{-\infty}^\infty \sin(t\lambda) R(\lambda) d\lambda
\eneq where we have set $R(\lambda):=R_V^+(\lambda^2)$ if
$\lambda>0$ and $R(\lambda)=\overline{R(-\lambda)}$ if
$\lambda<0$. For the free resolvent, we write this as $R_0$. Then,
by the usual resolvent expansions, \beeq \label{eq:res} R=
\sum_{k=0}^{2n-1} (-1)^k R_0(VR_0)^k + (R_0V)^n R (VR_0)^n \eneq
 As illustration, let us
consider the first term in this expansion which leads to the free
evolution. It is of the form (for $t>0$)
\begin{align*}
 &\Big|\il_{\R^6} \il_{-\infty}^\infty \sin(t\lambda)
\frac{e^{i\lambda|x-y|}}{|x-y|} \, d\lambda\, f(x)g(y)\,dx dy\Big| =
\f12 t^{-1}\Big|\il_{\R^3} \il_{[|x-y|=t]} f(x)\,\sigma(dx) \,g(y)\, dy\Big| \\
& \less t^{-1} \|\nabla f\|_1 \|g\|_1
\end{align*}
To pass to the second line we used the standard divergence theorem
trick (see eg.~Strauss~\cite{Str})
\begin{align}
&\Big|\il_{[|x-y|=t]} f(x)\,\sigma(dx) \Big|
 = \Big|\il_{[|x-y|=t]} f(x)\frac{x-y}{t}\cdot
\vec{n}\,\sigma(dx)\Big| = t^{-1}\Big|\il_{[|x-y|\le
t]} \nabla(f(x)(y-x))\, dx\Big| \label{eq:strauss}\\
&= t^{-1}\Big|\il_{[|x-y|\le t]} [(y-x)\nabla f(x)- 3f(x)] \,
dx\Big| \less \il_{\R^3} |\nabla f(x)|\, dx + \Big(\il_{\R^3}
|f(x)|^{\frac32}\,
dx\Big)^{\frac23} \nn\\
&\less \il_{\R^3} |\nabla f(x)|\, dx\nn
\end{align}
where the final inequality follows from Sobolev imbedding.

We distinguish between small energies and all other energies. For
the latter, we use~\eqref{eq:res}. Let $\chi_0(\lambda)=0$ for all
$|\lambda|\le \lambda_0$ and $\chi_0(\lambda)=1$ if
$|\lambda|>2\lambda_0$. Here $\lambda_0>0$ is some small parameter.
Fix some $k$ as in~\eqref{eq:res} and consider the contribution of
the corresponding Born term (ignoring a factor of $(4\pi)^{-k-1}$):
\begin{align*} &\il_{\R^{3(k+2)}} \il_{-\infty}^\infty
\chi_0(\lambda) \sin(t\lambda) e^{i\lambda\sum_{j=0}^{k}
|x_j-x_{j+1}|}
 \frac{\prod_{j=1}^k V(x_j)}{\prod_{j=0}^k|x_j-x_{j+1}|}
 \,f(x_0)g(x_{k+1})\,d\lambda\,dx_0\ldots dx_{k+1} \\ &
=\int\il_{\R^{3(k+1)}} \widehat{\chi_0}(\xi)
\il_{[|x_0-x_1|=t-\xi-\sum_{j=1}^{k} |x_j-x_{j+1}|>0]}
\frac{f(x_0)}{|x_0-x_1|}\, \sigma(dx_0)\,
 \frac{\prod_{j=1}^k V(x_j)}{\prod_{j=1}^k|x_j-x_{j+1}|}
 \,g(x_{k+1})\,dx_1\ldots dx_{k+1}\, d\xi
\end{align*}
By definition, $\chi_1=1-\chi_0\in C^\infty_c$ so that
$\widehat{\chi_0}=\delta_0-\widehat{\chi_1}$ with $\widehat{\chi_1}$
a Schwartz function. We start with the argument for $\delta_0$.
Write $\R^{3(k+1)} = A(t)\cup B(t)$ where \beeq\label{eq:AB_def}
A(t) = \Big\{t>\sum_{j=1}^{k} |x_j-x_{j+1}| > t/2 \Big\}, \qquad
 B(t) = \Big\{\sum_{j=1}^{k} |x_j-x_{j+1}| \le t/2 \Big\}
\eneq Then, using the divergence theorem as above and with
$\rho=t-\sum_{j=1}^{k} |x_j-x_{j+1}|$,
\begin{align*}
&\Big|\il_{A(t)} \il_{[|x_0-x_1|=\rho]} \frac{f(x_0)}{|x_0-x_1|}\,
\sigma(dx_0)\,
 \frac{\prod_{j=1}^k V(x_j)}{\prod_{j=1}^k|x_j-x_{j+1}|}
 \,g(x_{k+1})\,dx_1\ldots dx_{k+1} \Big|\\
&\less \il_{A(t)} \Big[\rho^{-1}\il_{[|x_0-x_1|<\rho]} |\nabla
f(x_0)| \, dx_0 + \rho^{-2}\il_{[|x_0-x_1|<\rho]} | f(x_0)| \,
dx_0\Big]
 \frac{\prod_{j=1}^k |V(x_j)|}{\prod_{j=1}^k|x_j-x_{j+1}|}
 \,|g(x_{k+1})|\,dx_1\ldots dx_{k+1}
 \end{align*}
We now use the definition of $A(t)$ to obtain a decay factor of $t^{-1}$ from one
of the denominators $|x_j-x_{j+1}|$. This allows us to further estimate the previous
expression by
\begin{align*}
 &\less k t^{-1} \sum_{\ell=1}^k \il_{\R^{3(k+2)}} |\nabla f(x_0)|\,
 \frac{\prod_{j=1}^k |V(x_j)|}{\prod_{0\le j<\ell}|x_j-x_{j+1}|\,\prod_{\ell<j\le k}|x_j-x_{j+1}|}
 \,|g(x_{k+1})|\,dx_0\ldots dx_{k+1} \\
& + k t^{-1} \sum_{\ell=1}^k \Big\|\il_{\R^3}
\f{|f(x_0)|}{|x_0-x_1|^2}\,dx_0 \Big\|_{L^3_{x_1}} \Big\|
\il_{\R^{3k}}
 \frac{\prod_{j=1}^k |V(x_j)|}{\prod_{1\le j<\ell}|x_j-x_{j+1}|\,\prod_{\ell<j\le k}|x_j-x_{j+1}|}
 \,|g(x_{k+1})|\,dx_2\ldots dx_{k+1} \Big\|_{L^{\f32}(dx_1)}\\
& \less k t^{-1} \|V\|_{\kato}^k\|\nabla f\|_1\|g\|_1  + kt^{-1}
\|V\|_{L^{\f32}} \|V\|_{\kato}^{k-1} \|\nabla f\|_1 \|g\|_1
\end{align*}
where
\[
\|V\|_{\kato} = \sup_y \int \frac{|V(x)|}{|x-y|}\, dx
\]
is the global Kato norm from~\cite{RodSch}, which is finite in our
case. We also used the bound
\[
\Big\|\il_{\R^3} \f{|f(x_0)|}{|x_0-x_1|^2}\,dx_0 \Big\|_{L^3(dx_1)}
\less \|f\|_{L^{\f32}(\R^3)} \less \|\nabla f\|_1
\]
which is obtained by fractional integration and Sobolev imbedding.
The estimate for the integral over $B(t)$ is similar. Here one uses
that $\rho(t)>t/2$ to gain a factor of $t^{-1}$. More precisely,
\begin{align*}
&\Big|\il_{B(t)} \il_{[|x_0-x_1|=\rho]} \frac{f(x_0)}{|x_0-x_1|}\,
\sigma(dx_0)\,
 \frac{\prod_{j=1}^k V(x_j)}{\prod_{j=1}^k|x_j-x_{j+1}|}
 \,g(x_{k+1})\,dx_1\ldots dx_{k+1} \Big|\\
&\less \il_{B(t)} \rho^{-1}\il_{\R^3} |\nabla f(x_0)| \, dx_0 \,
 \frac{\prod_{j=1}^k |V(x_j)|}{\prod_{j=1}^k|x_j-x_{j+1}|}
 \,|g(x_{k+1})|\,dx_1\ldots dx_{k+1} \\
 &\less  t^{-1} \int_{\R^{3(k+2)}} |\nabla f(x_0)|\,
 \frac{\prod_{j=1}^k |V(x_j)|}{\prod_{j=1}^k|x_j-x_{j+1}|}
 \,|g(x_{k+1})|\,dx_0\ldots dx_{k+1} \\
& \less t^{-1}  \|V\|_{\kato}^{k}\|\nabla f\|_1\|g\|_1
\end{align*}
We are done with $A(t)\cup B(t)$ and the $\delta_0$-part of
$\widehat{\chi_0}$.  Parenthetically, we remark that these arguments
(via the infinite expansion in~\eqref{eq:res}) prove the following
small potential result.

\begin{proposition}
\label{prop:small_pot1} Assume that the real-valued potential $V$
satisfies $\|V\|_{\kato}<4\pi$ and $\|V\|_{L^{\f32}}<\infty$. Then
one has the bound
\[ \Big\| \frac{\sin(t\sqrt{H})}{\sqrt{H}} f\Big\|_\infty \less
t^{-1} \|\nabla f\|_{L^{1}(\R^3)}
\]
for all $t>0$.
\end{proposition}

In the general case (i.e., large $V$ as needed in our application),
we need to work with the finite expansion~\eqref{eq:res}. Recall
that we yet have to deal with the contribution by
$\widehat{\chi_1}$, which means obtaining the same estimate as above
for
\[
\int\il_{\R^{3(k+1)}} \widehat{\chi_1}(\xi)
\il_{[|x_0-x_1|=t-\xi-\sum_{j=1}^{k} |x_j-x_{j+1}|>0]}
\frac{f(x_0)}{|x_0-x_1|}\, \sigma(dx_0)\,
 \frac{\prod_{j=1}^k V(x_j)}{\prod_{j=1}^k|x_j-x_{j+1}|}
 \,g(x_{k+1})\,dx_1\ldots dx_{k+1}\, d\xi
\]
Here we need to split the $\xi$-integral into the regions
$\{|\xi|<t/10\}$ and $\{|\xi|>t/10\}$. In the former region the same
argument applies as before, whereas in the latter one uses that
\[ \il_{[|\xi|>t/10]} |\widehat{\chi_0}(\xi)|\, d\xi < C_N\la t\ra^{-N}
\]
for any $N$ by the rapid decay of $\widehat{\chi_0}$.

It remains to bound the contribution by the final term
in~\eqref{eq:res}, the kernel $K(x,y)$ of which can be reduced to the form
\begin{align}
&\int e^{\pm it\lambda}\chi_0(\lambda) \la R(\lambda)
(VR_0(\lambda))^n(\cdot,x),
(VR_0(-\lambda))^n(\cdot,y) \ra \,d\lambda \nn\\
&= \int e^{i\lambda[\pm t+(|x|+|y|)]}\chi_0(\lambda) \la R(\lambda)
(VR_0(\lambda))^{n-1}VG_x(\lambda,\cdot),
(VR_0(-\lambda))^{n-1}VG_y(-\lambda,\cdot) \ra \,d\lambda
\label{eq:Gint}
\end{align}
Here
\[ G_{x}(\lambda,u) := \frac{e^{i\lambda(|x-u|-|x|)}}{4\pi|x-u|} \]
and the scalar product appearing in~\eqref{eq:Gint} is just another way of writing the composition of
the operators.
In \cite{GS} the following bounds were proved. For the sake of
completeness, we reproduce the simple proof.

\begin{lemma}  The derivatives of $G_{x}(\lambda,\cdot)$ satisfy the
estimates
\begin{equation}
\label{eq:Gest}
\begin{aligned}
\sup_{x\in\R^3} \Big\| \frac{d^j}{d\lambda^j} G_{x}(\lambda,
\cdot)\Big\|_{L^{2,-\sigma}} &< C_{j,\sigma} \text{\ \ provided\ \ } \sigma > \frac12 + j \\
\sup_{x\in\R^3} \Big\| \frac{d^j}{d\lambda^j}
G_{x}(\lambda,\cdot)\Big\|_{L^{2,-\sigma}} &<
\frac{C_{j,\sigma}}{\la x\ra} \text{\ \ provided\ \ }
\sigma>\frac32+j
\end{aligned}
\end{equation}
for all $j\ge 0$.
\end{lemma}
\begin{proof}
This follows from the explicit formula
\begin{align*}
\left\|\frac{d^j}{d\lambda^j} \frac{e^{i\lambda(|u-x|-|x|)}}{|x-u|}
\la u\ra^{-\sigma} \, \right\|_{L^2_u}
&=  \left( \int_{\R^3} \frac{(|u-x|-|x|)^{2j}}{|x-u|^2}\, \la u\ra^{-2\sigma}\, du\right)^{\f12} \\
&\le  \left( \int_{\R^3} \frac{\la u\ra^{2(j-\sigma)}}{|x-u|^2}\,
du\right)^{\f12}
\end{align*}
The final estimate on this integral is obtained by dividing $\R^3$
into the regions $|u|<\frac{|x|}2$, $|x-u|<\frac{|x|}2$, and the
complement of these two.  If $\frac12 < (\sigma - j) < \frac32$,
then each of these regions contributes $\la x\ra^{\frac12+j-\sigma}$
to the total. If $\sigma > \frac32 + j$, the first region instead
contributes $\la x \ra^{-2}$, making it the dominant term.
\end{proof}

Let
\[
a_{x,y}(\lambda) = \chi_0(\lambda) \la R(\lambda)
(VR_0(\lambda))^{n-1}VG_x(\lambda,\cdot),
(VR_0(-\lambda))^{n-1}VG_y(-\lambda,\cdot) \ra
\]
Then in view of the preceding one concludes that $a_{x,y}(\lambda)$
has two derivatives in~$\lambda$ and for large $n$,
\begin{equation}
\label{eq:adec}
\Big|\frac{d^j}{d\lambda^j} a_{x,y}(\lambda)\Big| \less
(1+\lambda)^{-2} (\la x\ra\la y\ra)^{-1}
\text{\ \ for\ \ } j = 0,1, \text{\ \ and all\ \ }\lambda>1 
\end{equation}
We need to take $n$ sufficiently large (say, $n>10$) in order to obtain
sufficiently fast decay of $a_{x,y}$ from the limiting absorption principle.
The latter here refers to the bounds for the free and
perturbed resolvents due to Agmon~\cite{Ag}: \begin{align}
\|R_V(\lambda^2\pm i0)\|_{L^{2,\sigma}\to L^{2,-\sigma}} &\less
\lambda^{-1+}, \qquad \sigma>\f12 \label{eq:limap}\\
\|\pr_\lambda^\ell R_V(\lambda^2\pm i0)\|_{L^{2,\sigma}\to
L^{2,-\sigma}} &\less 1, \qquad \sigma>\f12+\ell, \qquad \ell\ge1 \nn
\end{align}
for $\lambda$ separated from zero. Analogous estimates of course
hold for the free resolvent.

Let us assume first that $t>1$. To estimate~\eqref{eq:Gint} we
distinguish between $|t-(|x|+|y|)|<t/10$ and the opposite case.  In
the former case, we conclude that
\[ \max(|x|,|y|)\gtrsim t \]
so that due to \eqref{eq:adec} we obtain
\[ \Big| \int e^{i\lambda[\pm t+(|x|+|y|)]}
a_{x,y}(\lambda)\,d\lambda \Big| \less \chi_{[|x|+|y|>t]}(\la
x\ra\la y\ra)^{-1}\less t^{-1}
\]
In the latter case we integrate by parts once which also gains
$t^{-1}$, and we obtain the  bound
\[ \Big| \int e^{i\lambda[\pm t+(|x|+|y|)]}
a_{x,y}(\lambda)\,d\lambda \Big| \less t^{-1}(\la x\ra\la
y\ra)^{-1}\less t^{-1}
\]
Finally, if $0<t<1$, then we simply put the absolute values inside.
The conclusion is that in all cases \begin{align} & \sup_{x,y}\Big|
\int e^{i\lambda[\pm t+(|x|+|y|)]} \chi_0(\lambda) \la R(\lambda)
(VR_0(\lambda))^{n-1}VG_x(\lambda,\cdot),
(VR_0(-\lambda))^{n-1}VG_y(-\lambda,\cdot) \ra \,d\lambda \Big| \nn\\
& \less (\chi_{[|x|+|y|>t]} + \la t\ra^{-1})(\la x\ra\la y\ra)^{-1}\nn
\end{align}
Careful inspection of these bounds reveals that they only
require $|V(x)|\less \la x\ra^{-\kappa}$ with $\kappa>3$.
For example, consider the term which contains
\[ \pr_\lambda R_0(\lambda) V R_0(\lambda)\]
Then the resolvent on the left requires a weight of $\la x\ra^{\f32+\eps}$,
whereas the one the right requires~$\la x\ra^{\f12+\eps}$, see~\eqref{eq:limap}.
Hence, $V$ needs to absorb the weight $\la x\ra^{2+\eps}$.
On the other hand, by the lemma the term
\[  R_0(\lambda) V \pr_\lambda G_x(\lambda,\cdot)\]
requires weights of $\la x\ra^{\f52+\eps}$ for $\pr_\lambda G_x(\lambda,\cdot)$ (in order to gain $\la x\ra^{-1}$),
whereas $R_0(\lambda)$ needs $\la x\ra^{\f12+\eps}$. In total, this means that
$V$ has to absorb the weight $\la x\ra^{3+\eps}$, whence our assumption $|V(x)|\less \la x\ra^{-\kappa}$ with
$\kappa>3$.
See~\cite{GS} for similar details. In summary, we have obtained

\begin{lemma}
\label{lem:high}
 Assume $|V(x)|\less \la x\ra^{-\kappa}$ with
$\kappa>3$. Then there exists a kernel $K_t(x,y)$ so that
\[ |K_t(x,y)|\less (\chi_{[|x|+|y|>t]} + \la t\ra^{-1})(\la x\ra\la
y\ra)^{-1} \] and such that
\[ \Big\| \big[\frac{\sin(t\sqrt{H})}{\sqrt{H}} \chi_0(H) - K_t\big] f\Big\|_\infty
\less t^{-1} \| \nabla f\|_{L^1(\R^3)}
\]
for all $t>0$. In particular,
\[ \Big\| \frac{\sin(t\sqrt{H})}{\sqrt{H}} \chi_0(H) f\Big\|_\infty
\less t^{-1} \| f\|_{W^{1,1}(\R^3)}
\]
for all $t>0$.
\end{lemma}

We will not make direct use of the more refined bound involving
$K_t$ for the evolution $\frac{\sin(t\sqrt{H})}{\sqrt{H}}$. However,
we will need such a refined bound for the case of~$\cos(t\sqrt{H})$.

In passing, we remark that the previous argument can be easily
adapted to accommodate a  gradient. More precisely, we have

\begin{cor}
\label{cor:high}
 Assume $|V(x)|+|\nabla V(x)|\less \la x\ra^{-\kappa}$ with
$\kappa>3$. Then there exists a kernel $K_t(x,y)$ so that
\[ |K_t(x,y)|\less (\chi_{[|x|+|y|>t]} + \la t\ra^{-1})(\la x\ra\la
y\ra)^{-1} \] and such that
\[ \Big\| \big[\nabla\frac{\sin(t\sqrt{H})}{\sqrt{H}} \chi_0(H) - K_t\big] f\Big\|_\infty
\less t^{-1} (\| \nabla f\|_{L^1(\R^3)} + \| D^2 f\|_{L^1(\R^3)})
\]
for all $t>0$. In particular,
\[ \Big\| \nabla \frac{\sin(t\sqrt{H})}{\sqrt{H}} \chi_0(H) f\Big\|_\infty
\less t^{-1} \| f\|_{W^{2,1}(\R^3)}
\]
for all $t>0$.
\end{cor}
\begin{proof}
The point here is that the proof of the previous lemma is based on a
finite Born series expansion. We can commute the gradient through
the free resolvents in the terms of this series, which leads to
commutators of the gradient and the potential.
 These are harmless,
though, because of our decay assumption on $\nabla V$. Indeed, the
Born terms involving $\nabla V$ are of the same nature as those
arising in Lemma~\ref{lem:high}. As far as the final term in the
Born series (which involves a perturbed resolvent) is concerned, we
note that
\[ \nabla R - R \nabla = - R (\nabla V) R \]
But the right-hand side here does not make any difference to the way
we treated the function $a_{x,y}(\lambda)$ above. Hence the
corollary.
\end{proof}

\noindent It remains to deal with small energies. Recall that we are
assuming that zero energy is a resonance but not an eigenvalue.
Expansions for the perturbed resolvent around zero energy were
obtain by Jensen and Kato~\cite{JK} in that case. Here we will
follow~\cite{ES} in which the method of Jensen and
Nenciu~\cite{JenNen} was implemented for the case of~$\R^3$. Let us
recall the main steps: For $j=0,1,2,...,$ let $G_j$ be the operator
with the kernel
$$
G_j(x,y)=\frac{1}{4\pi j!}|x-y|^{j-1}.
$$
For each $J=0,1,2,...$, \[ R_0(\lambda)=\sum_{j=0}^{J} (i\lambda)^j
G_j+o(\lambda^J), \text{ as }\lambda\rightarrow 0. \] This expansion
is valid in the space, $HS_{L^{2,\sigma}\rightarrow L^{2,-\sigma}}$,
of Hilbert-Schmidt operators between $L^{2,\sigma}$ and
$L^{2,-\sigma}$ for $\sigma>\max((2J+1)/2,3/2)$. Let $U(x)=1$ if
$V(x)\geq 0$ and $U(x)=-1$ if $V(x)<0$, $v=|V|^{1/2}$ and $w=vU$. We
have
$$
V=Uv^2=wv.
$$
We use the symmetric resolvent identity, valid for $\lambda\ne0$,
\beeq\label{eq:res_exp} R(\lambda)=
R_0(\lambda)-R_0(\lambda)vA(\lambda)^{-1}vR_0(\lambda), \eneq with
the  operator
\begin{align*}
A (\lambda ) &= U + v R_0 (\lambda) v=(U + v G_0 v)+\lambda \frac{v[R_0 (\lambda) - G_0 ]v}{\lambda}\\
&=:A_0+\lambda A_1 (\lambda)
\end{align*}
Due to the compactness of $vG_0v$ on $L^2(\R^3)$ we remark that the essential
spectrum of $A_0$ is the same as that of $U$. In particular, if zero is in the
spectrum of $A_0$, then it is an isolated eigenvalue of finite multiplicity.
By inspection,  $A_1(\lambda)$ has the kernel
\begin{align}
A_1(\lambda)(x,y)&=\f{1}{\lambda}v(x)\f{e^{i\lambda|x-y|}-1}{4\pi|x-y|}v(y),\nn\\
|A_1(\lambda)(x,y)|&\leq\f{1}{4\pi}|v(x)|\;|v(y)|\label{eq:A1bd}
\end{align}
Therefore, $A_1(\lambda)\in HS:=HS_{L^2\rightarrow L^2}$ provided
$\langle x \rangle^{\f{3}{2}+} v(x) \in L^\infty$. Also note that
$$
A_1 (0) = iv G_1 v = \f{i\alpha}{4\pi} P_v,\;\;\alpha=\|V\|_1,
$$
where $P_v$ is the orthogonal projection onto span$(v)$.

In our case the operator $A(\lambda)$ is not invertible at
$\lambda=0$ due to the resonance. In fact, let $\psi\in
L^{2,-\frac12-}\setminus\{0\}$ solve $-\Delta\psi+V\psi =0$. Then
$\ker(A_0)=\R\cdot w\psi$ (where $\ker$ is relative to $L^2$). Let
$S_1$ be the orthogonal projection onto $\ker(A(0))$, i.e.,
\beeq
\label{eq:S1def} S_1 = \|w\psi\|_2^{-2} w\psi\otimes w\psi =: \tilde\psi\otimes \tilde\psi
\eneq
and set $A_0:=A(0)$. Then $A_0+S_1$ is invertible on $L^2(\R^3)$,
and $A(\lambda)+S_1$ is invertible, too, for small $\lambda$ (see
below). The following abstract lemma explains how to obtain the singular power of $\lambda$
by inverting $A(\lambda)$.

\begin{lemma}\cite{JenNen}\label{lem:jen-nen}
Let $F\subset \C\setminus\{0\}$ have zero as an accumulation point. Let $A(z)$, $z\in F$,
be a family of bounded operators on some Hilbert space of the form
$$A(z)=A_0+z A_1(z)$$
with $A_1(z)$ uniformly bounded as $z\rightarrow 0$ and $A_0$ self-adjoint. Suppose that $0$ is an isolated point of the spectrum of
$A_0$, and let $S$ be the corresponding Riesz projection. Assume that $\rank(S)<\infty$.
Then for sufficiently small $z\in F$ the operators
$$
B(z):=\frac{1}{z}(S-S(A(z)+S)^{-1} S)
$$
are well-defined and bounded on $\calH$. Moreover, due to $A_0=A_0^*$, they are
 uniformly bounded as $z\rightarrow 0$. The operator $A(z)$ has a bounded inverse in $\calH$ if and only if
$B(z)$ has a bounded inverse in $S\calH$, and in this case
\beeq\label{az-1}
A(z)^{-1}=(A(z)+S)^{-1}+\frac{1}{z}(A(z)+S)^{-1}SB(z)^{-1}S(A(z)+S)^{-1}.
\eneq
\end{lemma}

For the proof see~\cite{ES}. It follows from this lemma that
\beeq\label{S2=0} A(\lambda)^{-1} = (A(\lambda) +
S_1)^{-1}
 + \f{1}{\lambda} (A(\lambda) + S_1)^{-1} S_1
 m(\lambda)^{-1} S_1 (A(\lambda) + S_1)^{-1},
\eneq
provided
\[ m(\lambda)=\lambda^{-1}[S_1-S_1(A(\lambda)+S_1)^{-1}S_1] = m(0)+\lambda m_1(\lambda)\]
is invertible for small $\lambda$ on $S_1L^2$. However, this is indeed the case due to our assumption
that zero is not an eigenvalue of~$H$. In fact,
\[ m(0)= S_1 A_1(0) S_1 = \frac{i\alpha}{4\pi} S_1 P_v S_1 = \frac{i}{4\pi} \la v,\tilde\psi\ra^2 \tilde\psi
\otimes\tilde\psi = \f{i}{4\pi} \Big(\il_{\R^3} V\psi\, dx\Big)^2 \|w\psi\|_2^{-2}
S_1\]
is invertible on $S_1 L^2$ because of~\eqref{eq:intne0}.
Our claims about invertibility of these operators for
small~$\lambda$ will follow once we  justify the $L^2$-convergence
of the Neumann series
\begin{align}
\left(A (\lambda) + S_1\right)^{-1} &=(A_0 + S_1)^{-1} +
\sum_{k=1}^\infty (-1)^k \lambda^k (A_0 + S_1)^{-1}
\left[A_1 (\lambda) (A_0+S_1)^{-1}\right]^k\label{eq:Ainv}\\
&=:(A_0+S_1)^{-1}+\lambda E_1(\lambda), \nn\\
S_1m(\lambda)^{-1}S_1 & = S_1m(0)^{-1}S_1 +\sum_{k=1}^\infty (-1)^k
\lambda^k S_1m(0)^{-1}S_1
\left[m_1(\lambda) m(0)^{-1}S_1\right]^k\, S_1 \label{eq:m0inv}\\
&=:S_1m(0)^{-1}S_1+\lambda E_2(\lambda).\nn
\end{align}
We shall return to this issue later. Thus, using
$S_1(A_0+S_1)=(A_0+S_1)S_1=S_1$, we obtain
\begin{align}
\label{eq:Edef} A (\lambda)^{-1}  & = \frac{1}{\lambda}S_1 m(0)^{-1}
S_1+
 \left(A (\lambda) + S_1\right)^{-1}
\\
 & + E_1(\lambda) S_1
 m(0)^{-1} S_1 \left(A (\lambda) + S_1\right)^{-1}\nonumber \\
 &+  S_1
 m(0)^{-1} S_1 E_1(\lambda)  \nonumber \\
 &+   \left(A (\lambda) + S_1\right)^{-1} S_1
 E_2(\lambda) S_1 \left(A (\lambda) + S_1\right)^{-1} \nonumber\\
& =: -i\f{\beta}{\lambda}S_1+E(\lambda)\nonumber
\end{align}
with $\beta=4\pi \Big(\il_{\R^3} V\psi\, dx\Big)^{-2}\|w\psi\|_2^2$. Plugging
(\ref{eq:Edef}) into (\ref{eq:res_exp}), we have
\begin{align}
\label{eq:3piece}
 R(\lambda) &=i\f{\beta}{\lambda} R_0(\lambda)vS_1v
R_0(\lambda) +R_0 (\lambda) - R_0(\lambda)vE(\lambda)vR_0(\lambda).
\end{align}
Next, we describe the contribution of each of these terms to the
sine-transform~\eqref{eq:sintrans}. We can ignore the second one,
since it leads to the free case. The first term on the right-hand
side of~\eqref{eq:3piece} yields the following expression
in~\eqref{eq:sintrans}:
\begin{align*}
\calS_0(t)(x,y) &:= \f{\beta}{\pi} \int \f{\sin(t\lambda)}{\lambda}
\chi_1(\lambda)
\big[R_0(\lambda) vS_1v R_0(\lambda)\big](x,y) \, d\lambda \\
& = \|w\psi\|_2^{-2}\f{\beta}{2\pi} \il_{\R^6} \int_{-t}^t
e^{i\lambda\tau}\, d\tau\, \chi_1(\lambda)
e^{i\lambda[|x-x'|+|y'-y|]}
\f{V(x')\psi(x')\,V(y')\psi(y')}{4\pi|x-x'|\, 4\pi|y'-y|}\, d\lambda\,dx'dy' \\
&= \|w\psi\|_2^{-2}\f{\beta}{2\pi} \il_{\R^6} \int_{-t}^t
\widehat\chi_1(\tau+|x-x'|+|y'-y|)
\f{V(x')\psi(x')\,V(y')\psi(y')}{4\pi|x-x'|\, 4\pi|y'-y|}\,
d\tau\,dx'dy' \\
&= \|w\psi\|_2^{-2}{\beta} \il_{\R^6}
\f{V(x')\psi(x')\,V(y')\psi(y')}{4\pi|x-x'|\, 4\pi|y'-y|}\, dx'dy'\, \int\widehat{\chi_1}(\xi)\,d\xi \\
& - \|w\psi\|_2^{-2}\f{\beta}{2\pi} \il_{\R^6} \int_{[|\tau|>t]}
\widehat\chi_1(\tau+|x-x'|+|y'-y|)
\f{V(x')\psi(x')\,V(y')\psi(y')}{4\pi|x-x'|\, 4\pi|y'-y|}\,
d\tau\,dx'dy'
\end{align*}
Using that $(-\Delta)^{-1}V\psi=\psi$, we further conclude that
\begin{align}
& \f{\beta}{\pi} \int \f{\sin(t\lambda)}{\lambda} \chi_1(\lambda)
\big[R_0(\lambda) vS_1v R_0(\lambda)\big](x,y) \, d\lambda  =
\|w\psi\|_2^{-2}{\beta}\, \psi(x)\psi(y) \nn \\
&  - \|w\psi\|_2^{-2}\f{\beta}{2\pi} \il_{\R^6} \int_{[|\tau|>t]}
\widehat\chi_1(\tau+|x-x'|+|y'-y|)
\f{V(x')\psi(x')\,V(y')\psi(y')}{4\pi|x-x'|\, 4\pi|y'-y|}\,
d\tau\,dx'dy' \label{eq:VV}
\end{align}
We claim that the integral in~\eqref{eq:VV} is bounded by $t^{-1}$
as an operator from $L^1(\R^3)\to L^\infty(\R^3)$. To see this,
estimate
\begin{align*}
& \Big|\il_{\R^6}\il_{\R^6} \int_{[|\tau|>t]}
\widehat\chi_1(\tau+|x-x'|+|y'-y|)
\f{V(x')\psi(x')\,V(y')\psi(y')}{4\pi|x-x'|\, 4\pi|y'-y|}\,
d\tau\,dx'dy'\, f(x)g(y)\, dxdy\Big|\\
& \less \il_{[|x-x'|+|y-y'|<t/2]} \int_{[|\tau|>t]}
|\widehat\chi_1(\tau+|x-x'|+|y'-y|)|\, d\tau \,
\f{|V(x')\psi(x')|\,|V(y')\psi(y')|}{|x-x'|\,
|y'-y|}\,dx'dy'\, |f(x)g(y)|\, dxdy \\
&\quad +\il_{[|x-x'|+|y-y'|>t/2]}\int
|\widehat\chi_1(\tau+|x-x'|+|y'-y|)|\, d\tau
\,\f{|V(x')\psi(x')|\,|V(y')\psi(y')|}{|x-x'|\, |y'-y|}\, dx'dy'\,
|f(x)g(y)|\, dxdy \\
& \less \la t\ra^{-N} \Big(\sup_{x} \int
\f{|V(x')\psi(x')|}{|x-x'|}\, dx'\Big)^2 \|f\|_1\|g\|_1 + t^{-1}
\sup_{x} \int \f{|V(x')\psi(x')|}{|x-x'|}\, dx'\;\|V\psi\|_1
\|f\|_1\|g\|_1
\end{align*}
Since the expressions involving $V$ are finite constants (note
$\psi\in L^\infty$), we have proved our claim. The conclusion is
that
\begin{align*} \calS_0(t)(x,y) &=\f{\beta}{\pi} \int \f{\sin(t\lambda)}{\lambda} \chi_1(\lambda)
\big[R_0(\lambda) vS_1v R_0(\lambda)\big](x,y) \, d\lambda  = c(V)\,
\psi(x)\psi(y) + K_t(x,y),\\
\|K_t\|_{1\to\infty} &\less t^{-1}
\end{align*}

Finally, we turn to the third term on the right-hand side
of~\eqref{eq:3piece}. We start with the convergence of the Neumann
series~\eqref{eq:Ainv}. In view of~\eqref{eq:A1bd} we see that the
series for $(A(\lambda)+S_1)^{-1}$ converges for small $\lambda$
provided $V$ decays faster than a third power. By definition,
\[ m(\lambda)= -S_1 E_1(\lambda) S_1 \]
so that \begin{align} m_1(\lambda) &=\lambda^{-1}[m(\lambda)-m(0)]
= -\lambda^{-1}S_1(E_1(\lambda)-E_1(0))S_1 \nn\\
&= - \sum_{k=0}^\infty (-1)^{k} \lambda^{k} S_1 \left[A_1 (\lambda)
(A_0+S_1)^{-1}\right]^{k+2}S_1 +
\lambda^{-1}S_1(A_1(\lambda)-A_1(0))S_1 \label{eq:m1series}
\end{align}
The infinite series here again converges in $L^2$ for small
$\lambda$, whereas
\[
\lambda^{-1}S_1(A_1(\lambda)-A_1(0))S_1 = \lambda^{-2}S_1 v
(R_0(\lambda)-G_0-i\lambda G_1) vS_1
\]
admits the point-wise bound on its kernel
\begin{align*}
&\sup_\lambda \Big|\lambda^{-1}S_1(A_1(\lambda)-A_1(0))S_1(x,y)\Big|
\\ &\less |v(x)\psi(x)| \il_{\R^6} |\psi(x')V(x')||x'-y'||\psi(y')V(y')|\,
dx'dy'\; |v(y)\psi(y)| \\ &\less |v(x)\psi(x)| |v(y)\psi(y)|
\end{align*}
which holds provided $|V(x)|\less \la x\ra^{-\kappa}$ with
$\kappa>3$. Finally, the kernel $v(x)\psi(x)v(y)\psi(y)$ has finite
Hilbert-Schmidt norm. In summary, $m_1(\lambda)$ is an $L^2$-bounded
operator uniformly for small~$\lambda$. This proves that the Neumann
series~\eqref{eq:m0inv}  for $m(\lambda)^{-1}$ converges for
small~$\lambda$, as claimed.

We will need to control the contribution that each term in these
Neumann series makes to the sine-transform~\eqref{eq:sintrans}. We
start with the constant term, viz.
\[ E(0)=(A_0+S_1)^{-1} + E_1(0)S_1m(0)^{-1}S_1+S_1E_2(0)S_1 + S_1 m(0)^{-1}S_1 E_1(0) \]
see \eqref{eq:Edef}.  Thus,
\begin{align}
& \il_{-\infty}^\infty \sin(t\lambda) \chi_1(\lambda)
[R_0(\lambda)vE(0)vR_0(\lambda)](x,y)\, d\lambda \nn\\
& = \f{1}{16\pi^2} \il_{\R^6} \il_{-\infty}^\infty \sin(t\lambda)
e^{i\lambda[|x-x'|+|y'-y|]}\;\chi_1(\lambda) \, d\lambda\,
\f{v(x')E(0)(x',y')v(y')}{|x-x'|\,|y-y'|}\, dx'dy' \nn\\
&=\f{1}{32 i\pi^2} \il_{\R^6} \il_{-\infty}^\infty \delta( t + \xi +
[|x-x'|+|y'-y|])\;\widehat\chi_1(\xi) \, d\xi\,
\f{v(x')E(0)(x',y')v(y')}{|x-x'|\,|y-y'|}\, dx'dy' \label{eq:plus}\\
& \quad -\f{1}{32 i\pi^2} \il_{\R^6} \il_{-\infty}^\infty \delta(- t
+ \xi + [|x-x'|+|y'-y|])\;\widehat\chi_1(\xi) \, d\xi\,
\f{v(x')E(0)(x',y')v(y')}{|x-x'|\,|y-y'|}\, dx'dy' \label{eq:minus}
\end{align}
We start with the kernel $K_-(x,y;t)$ given by~\eqref{eq:minus}. Let
$\rho(t,\xi,y-y')=t-\xi-|y-y'|$ and, as usual, $t>0$.  Using
arguments similar to those in the large frequency case, we conclude
that
\begin{align}
& \Big|\il_{\R^6} K_-(x,y;t) f(x)g(y)\, dxdy\Big| \nn\\
& \less \Big| \il_{\R^9} \il_{[|\xi|<t/10]}
\il_{[|x-x'|=\rho(t,\xi,y-y')]} \f{f(x)}{|x-x'|}
\,\sigma(dx)\;\widehat\chi_1(\xi) \, d\xi\,
\f{v(x')E(0)(x',y')v(y')}{|y-y'|}\, dx'dy'\,g(y)\,dy \Big|\label{eq:xiklein}\\
& \quad + \Big| \il_{\R^9} \il_{[|\xi|>t/10]}
\il_{[|x-x'|=\rho(t,\xi,y-y')]} \f{f(x)}{|x-x'|}
\,\sigma(dx)\;\widehat\chi_1(\xi) \, d\xi\,
\f{v(x')E(0)(x',y')v(y')}{|y-y'|}\, dx'dy'\,g(y)\,dy
\Big|\label{eq:xigross}
\end{align}
Using the usual divergence theorem argument yields
\begin{align*}
\eqref{eq:xiklein} & \less \il_{[|y-y'|<t/2]} \il_{[|\xi|<t/10]}
\rho(t,\xi,y-y')^{-1}\,\il_{\R^3}|\nabla f(x)| \,
dx\;|\widehat\chi_1(\xi)| \, d\xi\,
\f{|v(x')E(0)(x',y')v(y')|}{|y-y'|}\, dx'dy'\,g(y)\,dy  \\
& \quad +  \il_{[|y-y'|>t/2]} \il_{[|\xi|<t/10]} \il_{\R^3}
\Big[\f{|\nabla f(x)|}{|x-x'|} +
\f{|f(x)|}{|x-x'|^2}\Big]\,dx\;|\widehat\chi_1(\xi)| \, d\xi\,
\f{|v(x')E(0)(x',y')v(y')|}{|y-y'|}\, dx'dy'\,g(y)\,dy \\
&\less t^{-1} \|\nabla f\|_1\|\widehat{\chi}\|_1 \|v\|_2
\|\,|E(0)(\cdot,\cdot)|\,\|_{2\to2}\;
\sup_{y}\Big\|\frac{v(y')}{|y-y'|}\Big\|_{L^2(dy')} \, \|g\|_1 \\
&\quad + t^{-1} \|\widehat{\chi}\|_1 \; \Big\| \int
\f{|f(x)||v(x')|}{|x-x'|^2}\,
dx\Big\|_{L^2(dx')}\|\,|E(0)(\cdot,\cdot)|\,\|_{2\to2}\; \|v\|_2\,
\|g\|_1 \\
&\less t^{-1} \|\nabla f\|_1 \|g\|_1
\end{align*}
To pass to the last line, we used the fact that $E(0)$ is an
absolutely bounded operator (in the terminology of~\cite{Sch}), see
Lemma~\ref{lem:abs_bd} below. Furthermore, we used that
$\|v\|_2<\infty$, as well as that
\[
\Big\| \int \f{|f(x)||v(x')|}{|x-x'|^2}\, dx\Big\|_{L^2(dx')} \less
\|v\|_{L^6} \Big\| \int \f{|f(x)|}{|x-x'|^2}\,
dx\Big\|_{L^3(dx')}\less  \|v\|_{L^6} \|f\|_{L^{\f32}} \less
\|\nabla f\|_1
\]
via fractional integration and Sobolev imbedding. The argument
for~\eqref{eq:xigross} is similar. Indeed, due to the rapid decay of
$\widehat{\chi_1}$ one obtains the bound
\[ \eqref{eq:xigross} \less \la t\ra^{-N} \|\nabla f\|_1 \|g\|_1 \]
for arbitrary $N\ge1$. In a similar vein, note that
in~\eqref{eq:plus} necessarily $\xi\le -t$. Hence this integral
contributes $\la t\ra^{-N}$. The conclusion is that
\beeq\label{eq:1sin} \Big| \il_{\R^6} \il_{-\infty}^\infty
\sin(t\lambda) \chi_1(\lambda)
[R_0(\lambda)vE(0)vR_0(\lambda)](x,y)\, d\lambda\, f(x) g(y)\, dxdy
\Big| \less t^{-1} \|\nabla f\|_1 \|g\|_1 \eneq For the following
lemma, we call an operator with kernel $K$ {\em absolutely bounded}
on $L^2$, provided the operator with kernel $|K(x,y)|$ is also $L^2$
bounded. The following lemma is quite standard, see~\cite{GS} and~\cite{Sch}
for similar considerations.

\begin{lemma}
\label{lem:abs_bd} The operator $(A_0+S_1)^{-1}$ is absolutely
bounded. In particular, $E(0)$ is also absolutely bounded.
\end{lemma}
\begin{proof}
Since $U^2 = I$ and
\[ A_0 + S_1 = U(I+ UvG_0 v + US_1) \]
we also have
\[ (A_0 + S_1)^{-1} = (I+ UvG_0 v + US_1)^{-1} U \]
The operator in parentheses on the right-hand side is a
Hilbert-Schmidt perturbation of the identity. Hence
\[ (I+ UvG_0 v + US_1)^{-1} - I \]
is again Hilbert-Schmidt, which implies that
$ (A_0 + S_1)^{-1} - U $
is also Hilbert-Schmidt and thus absolutely bounded. Since $U$ is
absolutely bounded, we are done.
\end{proof}

The same  argument which lead to~\eqref{eq:1sin} also yields the
following bound
\[
\Big| \il_{\R^6} \il_{-\infty}^\infty \sin(t\lambda) \chi_1(\lambda)
[R_0(\lambda)vF(\lambda)vR_0(\lambda)](x,y)\, d\lambda\, f(x) g(y)\,
dxdy \Big| \less t^{-1} \|\nabla f\|_1 \|g\|_1
\]
provided we replace $E(0)$ with an operator-valued function
$F(\lambda)$ satisfying\footnote{Below we will refer to this approach,
which is based on checking estimates~\eqref{eq:Fest1}, \eqref{eq:Fest2},
as the $F(\lambda)$ method.}
\begin{align}
\il_{-\infty}^\infty \Big\| \, |\widehat{\chi_1
F}(\xi)(\cdot,\cdot)|\, \Big\|_{2\to2} \, d\xi &<
\infty \label{eq:Fest1}\\
\Big\| \, |\widehat{\chi_1 F}(\xi)(\cdot,\cdot)|\, \Big\|_{2\to2}
&\less \la\xi\ra^{-1} \label{eq:Fest2}
\end{align}
The second estimate~\eqref{eq:Fest2} is needed for the case
$|\xi|>t/10$ to obtain $t$-decay, whereas~\eqref{eq:Fest1} suffices
in case $|\xi|<t/10$. In fact,
\begin{align*}
& \il_{-\infty}^\infty \sin(t\lambda) \chi_1(\lambda)
[R_0(\lambda)vF(\lambda)vR_0(\lambda)](x,y)\, d\lambda \nn\\
&=\f{1}{32 i\pi^2} \il_{[|\eta|<t/10]}\il_{\R^6}
\il_{-\infty}^\infty \delta( t + \xi +\eta +
[|x-x'|+|y'-y|])\;\widehat\chi_1(\xi) \, d\xi\,
\f{v(x')\hat{F}(\eta)(x',y')v(y')}{|x-x'|\,|y-y'|}\, dx'dy' d\eta\\
& \quad -\f{1}{32 i\pi^2} \il_{[|\eta|<t/10]}\il_{\R^6}
\il_{-\infty}^\infty \delta(- t + \xi +\eta +
[|x-x'|+|y'-y|])\;\widehat\chi_1(\xi) \, d\xi\,
\f{v(x')\hat{F}(\eta)(x',y')v(y')}{|x-x'|\,|y-y'|}\, dx'dy' d\eta\\
& \quad + \f{1}{32 i\pi^2}  \il_{[|\eta|>t/10]}\il_{\R^6}
\widehat\chi_1( t +\eta + [|x-x'|+|y'-y|])\;
\f{v(x')\hat{F}(\eta)(x',y')v(y')}{|x-x'|\,|y-y'|}\, dx'dy' d\eta\\
& \quad -\f{1}{32 i\pi^2}  \il_{[|\eta|>t/10]}\il_{\R^6}
\widehat\chi_1(- t +\eta + [|x-x'|+|y'-y|]) \;
\f{v(x')\hat{F}(\eta)(x',y')v(y')}{|x-x'|\,|y-y'|}\, dx'dy' d\eta
\end{align*}
In the first two integrals, we distinguish $|\xi|>t/10$ from
$|\xi|<t/10$ and use~\eqref{eq:Fest1}, whereas in the third and
fourth, we use~\eqref{eq:Fest2} and the integrability of
$\hat{\chi}_1$. The reader should not be confused by the fact that
in the first two integrals the $\xi$-integration has not been carried out
whereas in the final two it has. This is due to the fact that in the first
two integrals the method from before using spherical integration is  needed
in order to gain a factor of~$t^{-1}$, whereas in the final two it is not.
Indeed, for the final two integrals a gain of~$t^{-1}$ is obtained from~\eqref{eq:Fest2},
and the integral in $\eta$ is reduced to integrating out~$\hat{\chi}_1$.

We shall now prove \eqref{eq:Fest1} and \eqref{eq:Fest2} for
$F(\lambda)=E(\lambda)$ where
 \begin{align}
 E(\lambda)  &= \left(A (\lambda) + S_1\right)^{-1}
 + E_1(\lambda) S_1
 m(0)^{-1} S_1 \left(A (\lambda) + S_1\right)^{-1}
 +   S_1
 m(0)^{-1} S_1 E_1(\lambda) \nn\\
 & + \left(A (\lambda) + S_1\right)^{-1} S_1
 E_2(\lambda) S_1 \left(A (\lambda) + S_1\right)^{-1}
 \label{eq:4piece}
 \end{align}
We start with the case of  $F(\lambda)=(A(\lambda)+S_1)^{-1}$. In
view of the Neumann series of $(A(\lambda)+S_1)^{-1}$ it will
further suffice to prove\footnote{Recall that
$[-2\lambda_0,2\lambda_0]$ is  the support of $\chi_1$} that
\begin{align}
\il_{-\infty}^\infty \Big\| \, |[\lambda\chi_1(\lambda)
A_1(\lambda)]^{\wedge}(\xi)(\cdot,\cdot)|\, \Big\|_{2\to2} \, d\xi &
\less \lambda_0 \label{eq:A1est1}\\
\Big\| \, |[\lambda\chi_1(\lambda)
A_1(\lambda)]^{\wedge}(\xi)(\cdot,\cdot)|\, \Big\|_{2\to2} &\less
\lambda_0 \la \xi\ra^{-1} \label{eq:A1est2}
\end{align}
More precisely, to pass from~\eqref{eq:A1est1} to~\eqref{eq:Fest1}
for $F(\lambda)=(A(\lambda)+S_1)^{-1}$ we use the following lemma
from~\cite{ES}, applied to each of the terms in the Neumann series
of $(A(\lambda)+S_1)^{-1}$. For small $\lambda_0$ we then obtain a
summable series.

\begin{lemma}\label{L:convolve}
For each $\lambda\in \R$, let $F_1(\lambda)$ and $F_2(\lambda)$ be
bounded operators from $L^2(\R^3)$ to $L^2(\R^3)$ with kernels
$K_1(\lambda)$ and $K_2(\lambda)$. Suppose that $K_1, K_2$ both have
compact support in $\lambda$ and that $K_j(\cdot)(x,y)\in L^1(\R)$
for a.e.~$x,y\in\R^3$. Let $F(\lambda)=F_1(\lambda)\circ
F_2(\lambda)$ with kernel $K(\lambda)$. Then
$$
\int_{-\infty}^\infty
\left\|\,\left|\widehat{K}(\xi)\right|\,\right\|_{ 2\rightarrow 2}
d\xi \leq \left[\int_{-\infty}^\infty \left\|\, \left|\widehat{K_1
}(\xi)\right|\,\right\|_{ 2\rightarrow 2}d\xi\right]
\left[\int_{-\infty}^\infty \left\|\, \left|\widehat{K_2
}(\xi)\right|\,\right\|_{ 2\rightarrow 2}d\xi\right].
$$
\end{lemma}

This is basically just an operator version of the fact that
convolution does not increase $L^1$-norms.  For more details, see
Lemma~8 in~\cite{ES}. To pass from~\eqref{eq:A1est2}
to~\eqref{eq:Fest2},  we use the same idea, but in addition we need
to take the supports of the convolutions into account. More
precisely, we use that the support of a convolution is the
arithmetic sum of the supports of the individual convolution
factors. Hence, if a $k$-fold convolution is being evaluated at
$|\xi|>L$, then at least one of the factors needs to be evaluated at
$|\xi|>L/k$. On this factor we use the point-wise estimate.  This
leads to a loss of a polynomial factor $k^{2}$ in our
estimate~\eqref{eq:Fest2} of a $k$-fold convolution, which can then
be absorbed into the exponential gain~$\lambda_0^k$.

Let us now prove \eqref{eq:A1est1}, \eqref{eq:A1est2}. First, we
write $\chi_1(\lambda)=\chi(\lambda/\lambda_0)$. Then
\begin{align*}
[\chi(\lambda/\lambda_0)\lambda A_1(\lambda)]^{\wedge}(\xi)(x,y) &=
[v(y)\chi(\lambda/\lambda_0)(R_0(\lambda)-G_0)v(x)]^{\wedge}(\xi) \\
&= \lambda_0 \f{v(x)v(y)}{4\pi |x-y|} \big(
\hat{\chi}(\lambda_0(\xi+|x-y|))-\hat{\chi}(\lambda_0\xi)\big)
\end{align*}
Hence,
\begin{align*}
\Big|[\chi(\lambda/\lambda_0)\lambda
A_1(\lambda)]^{\wedge}(\xi)(x,y)\Big| &\less \lambda_0
\f{|v(x)v(y)|}{|x-y|} \Big|
\int_{\lambda_0\xi}^{\lambda_0(\xi+|x-y|)} \hat{\chi}'(u)\, du\Big|\\
&\less \lambda_0^2 |v(x)||v(y)| \la \lambda_0\xi\ra^{-N} +\lambda_0
\f{|v(x)v(y)|}{|x-y|}\chi_{[|x-y|>|\xi|]}
\end{align*}
and thus, with $HS$ denoting the Hilbert-Schmidt norm,
\begin{align*}
& \int \Big\| [\chi(\lambda/\lambda_0)\lambda
A_1(\lambda)]^{\wedge}(\xi)\Big\|_{HS}\, d\xi \\
\quad &\less  \lambda_0^2 \int \|v(x)v(y)\|_{L^2(\R^6)} \la
\lambda_0\xi\ra^{-N}\, d\xi +\lambda_0 \int \Big\|
\f{|v(x)v(y)|}{|x-y|}\chi_{[|x-y|>|\xi|]} \Big\|_{L^2(\R^6)} \, d\xi
 \less \lambda_0
\end{align*}
where we used that \beeq \Big(\il_{[|x-y|>|\xi|]} \f{|v(x)|^2
|v(y)|^2}{|x-y|^2}\, dxdy\Big)^{\f12} \less \la
\xi\ra^{(1-\kappa)/2} \label{eq:HSxi} \eneq provided $|V(x)|\less
\la x\ra^{-\kappa}$. So as long as $\kappa>3$, which is also needed
for $\|v\|_2<\infty$,  we obtain a bound which is integrable
in~$\xi$. This proves~\eqref{eq:A1est1}. Moreover, the point-wise
bounds also hold:
\[
 \Big\| [\chi(\lambda/\lambda_0)\lambda
A_1(\lambda)]^{\wedge}(\xi)\Big\|_{HS} \less  \lambda_0^2 \la
\lambda_0\xi\ra^{-N} + \lambda_0\la \xi\ra^{(1-\kappa)/2}
 \less \lambda_0 \la \xi\ra^{-1}
\]
if we choose $N=1$.

It remains to deal with the rank-one pieces of $E(\lambda)$, which
are the last three terms in~\eqref{eq:4piece}:
\begin{align*}
F_1(\lambda) &:= E_1(\lambda) S_1
 m(0)^{-1} S_1 \left(A (\lambda) + S_1\right)^{-1} \\
F_2(\lambda) &:= S_1
 m(0)^{-1} S_1 E_1(\lambda) \\
F_3(\lambda) &:= \left(A (\lambda) + S_1\right)^{-1} S_1
 E_2(\lambda) S_1 \left(A (\lambda) + S_1\right)^{-1}
\end{align*}
To analyze $F_1, F_2$, we write
\begin{align*}
E_1(\lambda) &=
-(A_0+S_1)^{-1}A_1(\lambda)(A_0+S_1)^{-1}G_1(\lambda),\qquad
G_1(\lambda):= \sum_{\ell=0}^\infty (-1)^\ell [\lambda A_1(\lambda)
(A_0+S_1)^{-1}]^\ell
\end{align*}
The advantage of this representation lies with the fact that
$G_1(\lambda)$ satisfies~\eqref{eq:Fest1} and~\eqref{eq:Fest2}. This
follows by the same arguments that established the same property for
$(A(\lambda)+S_1)^{-1}$. To show that $F_{1,2}$ satisfy
\eqref{eq:Fest1} and~\eqref{eq:Fest2} it will therefore suffice to
show that $A_1(\lambda)$ does. To this end, compute
\begin{align*}
\widehat{A_1}(\xi)(x,y) &= (4\pi)^{-1}v(y) \Big[\il_0^1
e^{i\lambda|x-y|b}\,db\Big]^{\wedge}(\xi)\,v(x)
= (4\pi)^{-1} v(y)\il_0^1 \delta(\xi-b|x-y|)\, db\, v(x) \\
&= (4\pi)^{-1} v(y) \chi_{[0<\xi<|x-y|]} \frac{v(x)}{|x-y|}
\end{align*}
Hence,
\[
 \Big\|\widehat{A_1}(\xi) \Big\|_{HS} \less \Big\|\chi_{[0<\xi<|x-y|]} \frac{v(x)
 v(y)}{|x-y|}\Big\|_{L^2(\R^6)}
 \less \la \xi\ra^{(1-\kappa)/2},
\]
see \eqref{eq:HSxi}. Since $\kappa>3$, this implies all the desired
properties of $A_1$.

To deal with $F_3$, we use the representation
\begin{align*}
 E_2(\lambda) &= -S_1m(0)^{-1}S_1m_1(\lambda)S_1
m(0)^{-1}S_1 G_2(\lambda), \qquad G_2(\lambda):=
\sum_{\ell=0}^\infty (-1)^\ell \,S_1[\lambda m_1(\lambda)
m(0)^{-1}]^\ell\, S_1
\end{align*}
This shows that it suffices to prove \eqref{eq:Fest1}
and~\eqref{eq:Fest2} for $m_1(\lambda)$. Indeed, the infinite series
for $G_2$ will inherit these properties due to the geometric
factors~$\lambda^\ell$. Returning to $m_1(\lambda)$, we use the
representation~\eqref{eq:m1series} which allows us to write
\[ m_1(\lambda) = -
S_1A_1(\lambda)(A_0+S_1)^{-1}A_1(\lambda)(A_0+S_1)^{-1}G_1(\lambda)S_1
+\lambda^{-2}S_1v(R_0(\lambda)-G_0-i\lambda G_1)vS_1 \] In view of
the preceding, the first term on the right-hand side here has the
desired properties. The kernel of the second term here equals
\[ c\, w(x)\psi(x) \il_{\R^6} \psi(x')V(x')
\f{e^{i\lambda|x'-y'|}-1-i\lambda|x'-y'|}{\lambda^2|x'-y'|}
\psi(y')V(y')\, dx'dy' \; w(y)\psi(y)
\]
To compute the Fourier transform in $\lambda$, we use
\[ \f{e^{i\lambda a}-1-ia\lambda}{\lambda^2 a} = -a\il_0^1
(1-b)e^{i\lambda a b}\, db \] Therefore, for $a>0$,
\[ \Big[\f{e^{i\lambda a}-1-ia\lambda}{\lambda^2 a}\Big]^{\wedge}(\xi) = -a\il_0^1
(1-b)\delta(\xi-ba)\, db = -(1-\xi/a)\chi_{[0<\xi<a]} \] from which
we conclude that
\begin{align*}
& \Big|\il_{\R^6} V(x')\psi(x') \Big[\f{e^{i\lambda
|x'-y'|}-1-i|x'-y'|\lambda}{\lambda^2 |x'-y'|}\Big]^{\wedge}(\xi)
V(y')\psi(y')\, dx'dy'\Big| \\
& = \Big|\chi_{[\xi>0]} \il_{[|x'-y'|>\xi]} V(x')\psi(x')
\Big(1-\f{\xi}{|x'-y'|}\Big) V(y')\psi(y')\, dx'dy' \Big| \\
&\less \la \xi\ra^{2-\kappa+\eps}
\end{align*}
for arbitrary $\eps>0$. It follows that both \eqref{eq:Fest1},
\eqref{eq:Fest2} hold for $m_1$.

To summarize our small energy results, we have proved the following
lemma.

\begin{lemma}
\label{lem:low} Under the assumptions of
Proposition~\ref{prop:infty_decay} there exists a constant $c_0\ne0$
such that
\[ \calS_0(t)(x,y) = c_0\il_{\R^6} \int_{-t}^t
\widehat\chi_1(\tau+|x-x'|+|y'-y|)
\f{V(x')\psi(x')\,V(y')\psi(y')}{4\pi|x-x'|\, 4\pi|y'-y|}\,
d\tau\,dx'dy'\] satisfies
\begin{equation}
\label{eq:sin_1inf} \Big\| \f{\sin(t\sqrt{H})}{\sqrt{H}}\chi_1(H)
P_c f - \calS_0(t) f\Big\|_\infty \less t^{-1} \|f\|_{W^{1,1}(\R^3)}
\end{equation}
for all $t>0$. Also, 
\[ \|\calS_0(t) - c_0 \,(\psi\otimes \psi)  \|_{1\to\infty} \less
t^{-1} \] In particular, $\calS_0(t)\to c_0\,(\psi\otimes\psi)$ as
$t\to\infty$.
\end{lemma}

It may be useful to retain $\calS_0(t)$ in some cases, since its
kernel basically lives on the set \[\{ |x|+|y|\less t\}\] in a weak
sense.
Proposition~\ref{prop:infty_decay}
follows by combining Lemmas~\ref{lem:high} and Lemma~\ref{lem:low}.

Next, we state a corollary for the differentiated evolution. In
conjunction with Corollary~\ref{cor:high} this
proves~\eqref{eq:singrad}.

\begin{cor}
\label{cor:sinder}
 Let $V$ satisfy $|V(x)|+|\nabla V(x)|\less \la
x\ra^{-\kappa}$ with $\kappa>3$. Also, assume that $V<0$ point-wise.
Then
\begin{equation}
\label{eq:1inf2} \Big\| \nabla\f{\sin(t\sqrt{H})}{\sqrt{H}} P_c f -
c_0 ((\nabla\psi)\otimes \psi) f\Big\|_\infty \less t^{-1}
\|f\|_{W^{2,1}(\R^3)}
\end{equation}
for all $t>0$.
\end{cor}
\begin{proof}
We have already settled the high-energy case, see
Corollary~\ref{cor:high}. To deal with the low energies, we need to
pass a gradient through the low energy propagator above. Our
assumption $V<0$ implies that $U=-1$. This allows us to commute
$\nabla$ with $U$, and we can also differentiate $v,w,\psi$.  The
term $\nabla \calS_0(t)$ does not present a problem: It gives
$(\nabla \psi)\otimes\psi$ up to an operator bounded by $t^{-1}$
from $L^1\to L^\infty$. Next, we need to pass a gradient through the
operator
\begin{align*}
 E(\lambda)  &= \left(A (\lambda) + S_1\right)^{-1}
 + E_1(\lambda) S_1
 m(0)^{-1} S_1 \left(A (\lambda) + S_1\right)^{-1}
 +   S_1
 m(0)^{-1} S_1 E_1(\lambda) \nn\\
 & + \left(A (\lambda) + S_1\right)^{-1} S_1
 E_2(\lambda) S_1 \left(A (\lambda) + S_1\right)^{-1}
\end{align*}
First,
\[
[\nabla,(A(\lambda)+S_1)^{-1}] =
(A(\lambda)+S_1)^{-1}([S_1,\nabla]+[vR_0(\lambda)v,\nabla])(A(\lambda)+S_1)^{-1}
\]
and
\[ [S_1,\nabla] = (\nabla \tilde\psi)\otimes\tilde\psi +
\tilde\psi\otimes (\nabla\tilde\psi) \]
as well as
\[
[vR_0(\lambda)v,\nabla] = - (\nabla v)R_0(\lambda) v - vR_0(\lambda)
(\nabla v)
\]
Clearly, similar expressions arise when commuting $\nabla$ with
$A_0+S_1$. Second, using these relations allows us to commute
$\nabla$ through $m(0)^{-1}$, $E_1(\lambda)$, $A_1(\lambda)$,  and
$m_1(\lambda)$. In order to apply the $F(\lambda)$-machinery from
above, we need to check that each of the terms arising as a
commutator satisfies~\eqref{eq:Fest1}, \eqref{eq:Fest2}. For
example, consider the term
\[ (A(\lambda)+S_1)^{-1}[vR_0(\lambda)v,\nabla](A(\lambda)+S_1)^{-1}
= (A(\lambda)+S_1)^{-1}[(\nabla v)R_0(\lambda) v + vR_0(\lambda)
(\nabla v)](A(\lambda)+S_1)^{-1}
\]
We already know that $(A(\lambda)+S_1)^{-1}$ satisfy these bounds,
so by Lemma~\ref{L:convolve} it suffices to check this for $(\nabla
v)R_0(\lambda) v$ and $vR_0(\lambda) (\nabla v)$. Ignoring rapidly
decaying tails, this amounts to showing that
\[
h(\xi) := \Bigl( \int\int \f{|V(x)||V(y)|}{|x-y|^2}\,
\chi_{[||x-y|-|\xi|<1]} \, dxdy \Bigr)^{\half}
\]
satisfies
\begin{equation}
\label{eq:new_piece} \il_{-\infty}^\infty h(\xi)\, d\xi <\infty, \qquad h(\xi)\less \la
\xi\ra^{-1}
\end{equation}
Since
\[ h(\xi)\less \la \xi\ra^{-1} \Bigl( \int\int {|V(x)||V(y)|}\,
\chi_{[||x-y|-|\xi|<1]} \, dxdy \Bigr)^{\half}
\]
the latter bound follows immediately since $V\in L^1$ and the former
follows by applying Cauchy-Schwarz. The other commutators can be
checked similarly.
\end{proof}

\subsection{The cosine evolution}

Next, we state an estimate for $\cos(t\sqrt{H})P_c$.

\begin{prop}\label{prop:cos}
Under the assumptions of Proposition~\ref{prop:infty_decay}, there
is a kernel $K_t$ satisfying~\eqref{eq:Kt} so that the bound
\begin{equation} \label{eq:cos_dec}
\Big\| [\cos(t\sqrt{H})P_c - K_t]f \Big \|_\infty \less t^{-1}
\sum_{1\le|\alpha|\le 2}\|D^\alpha f\|_{L^{1}(\R^3)}
\end{equation}
holds for all $t>0$. Moreover, if $V$ is as in
Corollary~\ref{cor:sinder}, then there is another kernel $\tilde
K_t$ satisfying~\eqref{eq:Kt} so that
\begin{equation} \label{eq:Dcos_dec}
\Big\| [\nabla \cos(t\sqrt{H})P_c - \tilde K_t]f \Big \|_\infty
\less t^{-1} \sum_{1\le|\alpha|\le 3}\|D^\alpha f\|_{L^{1}(\R^3)}
\end{equation}
holds for all $t>0$.
\end{prop}

\begin{remark}
\label{rem:weakL1} The main application of this proposition is the
following estimate, which does require using $K_t$ rather than an
$L^1$ norm of $f$ on the right-hand side:
\begin{equation}
\label{eq:cos2} \|\cos((t-s)\sqrt{H}) P_g^\perp
[\phi_a(\cdot,b(s))-(a(\infty)/b(s))^{\f54}\phi_a(\cdot,a(\infty))] \|_\infty \less
\delta \, \la t-s\ra^{-1}
\end{equation}
 see Section~\ref{sec:stability}. First, if $0<t-s<1$, then via Sobolev's
imbedding, with $\psi=\phi_a(\cdot,a(\infty))$,
\begin{align*}
& \|\cos((t-s)\sqrt{H}) P_g^\perp
[\phi_b(\cdot,b(s))-(a(\infty)/b(s))^{\f54}\psi] \|_\infty \\
&\less \|\cos((t-s)\sqrt{H}) P_g^\perp
[\phi_b(\cdot,b(s))-(a(\infty)/b(s))^{\f54}\psi] \|_2 \\
& \qquad + \|D^2 \cos((t-s)\sqrt{H}) P_g^\perp
[\phi_a(\cdot,a(s))-(a(\infty)/b(s))^{\f54}\psi] \|_2 \\
&\less \| P_g^\perp [\phi_b(\cdot,b(s))-(a(\infty)/b(s))^{\f54}\psi]
\|_2 + \|H \cos((t-s)\sqrt{H}) P_g^\perp
[\phi_b(\cdot,b(s))-(a(\infty)/b(s))^{\f54}\psi] \|_2 \\
&\less \| \phi_b(\cdot,b(s))-(a(\infty)/b(s))^{\f54}\psi] \|_2 + \|H
[\phi_b(\cdot,b(s))-(a(\infty)/b(s))^{\f54}\psi] \|_2 \\
&\less \delta\,\la s\ra^{-1}
\end{align*}
Second, let $t-s>1$.  The problem here is that the function in
brackets decays only like $\la x\ra^{-3}$ and is therefore not in
$L^1$, but only in weak $L^1$. However, its first and second
derivatives are in $L^1$, so the bound on $K_t$ saves us,
see~\eqref{eq:Kt}. Indeed, since
 \[ \int (\chi_{[|x|+|y|>t]}+\la t\ra^{-1}) \la x\ra^{-1}
\la y\ra^{-1} \; \la y\ra^{-3}\, dy \less \la t\ra^{-1}, \] we
obtain~\eqref{eq:cos2} as desired.
\end{remark}

\begin{proof}[Proof of Proposition~\ref{prop:cos}]
As before, we shall only deal with the case when a resonance is present.
The other case is implicit in what we are doing.
We first discuss the estimate \eqref{eq:cos_dec} without the
gradient. The proof proceeds by making appropriate changes to the
preceding proof. The logic here is that $\cos(t\sqrt{H})P_c = \pr_t
\f{\sin(t\sqrt{H})}{\sqrt{H}}P_c$, as in the free case. Since we
have just written
\[
\f{\sin(t\sqrt{H})}{\sqrt{H}}P_c = \calS_0(t) + \calS_1(t)
\]
where $\calS_1(t)$ is dispersive, we  obtain that
\[
\cos(t\sqrt{H})P_c = \dot\calS_0(t) + \dot\calS_1(t)
\]
with
\[ \dot{\calS}_0(t)(x,y)
= c \il_{\R^6} [\widehat\chi_1(t+|x-x'|+|y-y'|) +
\widehat\chi_1(-t+|x-x'|+|y-y'|)]
\f{V(x')\psi(x')\,V(y')\psi(y')}{4\pi|x-x'|\, 4\pi|y'-y|}\,dx'dy'
\]
Via Lemma~\ref{lem:tedious} below one can read off here that \[
|\dot{\calS}_0(t)(x,y)|\less (\la t\ra^{-1}+\chi_{[|x|+|y|>t]}) (\la
x\ra\la y\ra)^{-1}
\]
which implies that $\dot{\calS}_0(t)$ can be included into the
kernel $K_t$. In particular, \[ \|\dot{\calS}_0(t)f\|_\infty \less
\la t\ra^{-1} \|f\|_1
\]
Moreover, it is implicit in the arguments for
$\sin(t\sqrt{H})/\sqrt{H}$ that $\dot\calS_1(t)$ dispersive. To make
this explicit, first recall the estimate for the free propagator:
\[ \|\cos(t\sqrt{-\Delta}) f\|_\infty \less t^{-1} \|D^2 f\|_{L^1(\R^3)} \]
This can be seen by an argument similar to~\eqref{eq:strauss} above:
\begin{align*}
\cos(t\sqrt{H}) f(x) &= \pr_t\; t \il_{S^2} f(x+ty)\, \sigma(dy) \\
& = \il_{S^2} \big[ f(x+ty)+ t(\nabla f)(x+ty)\cdot y\big]\,
\sigma(dy) \\
&= t^{-2} \il_{[|x-y|\le t]} \Big[\nabla f(y)\cdot \f{y-x}{t} + f(y)
\f{3}{t}\Big]\, dy + t^{-1}\il_{[|x-y|\le t]} \Delta f(y)\, dy
\end{align*}
Hence,
\begin{align*}
|\cos(t\sqrt{H}) f(x)| &\less t^{-2} \il_{[|x-y|\le t]} \Big[|\nabla
f(y)| + t^{-1}|f(y)| \Big]\, dy + t^{-1}\il_{[|x-y|\le t]} |\Delta
f(y)|\, dy \\
&\less t^{-1} \Big(\il_{\R^3} |\nabla f(y)|^{\f32}\, dy\Big)^{\f23}
+ t^{-1} \Big(\il_{\R^3} | f(y)|^{3}\, dy\Big)^{\f13} + t^{-1} \|D^2
f\|_1 \\
&\less t^{-1}\|D^2 f\|_1
\end{align*}
To pass to the final inequality we used there the Sobolev imbedding
\[ \|\nabla f\|_{\f32} \less \|D^2 f\|_1, \quad \|f\|_3 \less
\|\nabla f\|_{\f32} \less \|D^2 f\|_1 \] For the perturbed case, we
distinguish small energies from all other energies. We start by
indicating the changes to the argument for energies
$\lambda\in[\lambda_0,\infty)$ that we presented above for
$\f{\sin(t\lambda)}{\lambda}$. First, recall~\eqref{eq:AB_def}. As
before, with $\rho=t-\sum_{j=1}^{k} |x_j-x_{j+1}|$, we now need to
estimate
\[
\Big|\pr_t\; \il_{A(t)} \il_{[|x_0-x_1|=\rho(t)]}
\frac{f(x_0)}{|x_0-x_1|}\, \sigma(dx_0)\,
 \frac{\prod_{j=1}^k V(x_j)}{\prod_{j=1}^k|x_j-x_{j+1}|}
 \,g(x_{k+1})\,dx_1\ldots dx_{k+1} \Big|
\]
The effect of the $\pr_t$ is twofold: We either need to replace $f$
with $\nabla f(\cdot)\cdot \vec{n}$ (with an outward pointing normal
vector $\vec{n}$), or we keep $f$ but replace $|x_0-x_1|$ in the
denominator with $|x_0-x_1|^2=\rho(t)^2$. In the former case, the
only change needed is that the final bound is in terms of $\|D^2
f\|_{L^1(\R^3)}$. In the latter case, the same approach to $A(t)$ as
above leads to  two new terms, viz.
\begin{align}
\int \f{|\nabla f(x_0)|}{|x_0-x_1|^2} \, dx_0 \label{eq:new1} \\
\rho^{-3}\il_{[|x_0-x_1|<\rho]} | f(x_0)| \, dx_0 \label{eq:new2}
\end{align}
To deal with \eqref{eq:new1} we apply fractional integration and
Sobolev imbedding
\[
\Big\|\il_{\R^3} \f{|\nabla f(x_0)|}{|x_0-x_1|^2}\,dx_0
\Big\|_{L^3(dx_1)} \less \|\nabla f\|_{L^{\f32}(\R^3)} \less \|D^2
f\|_1
\]
As for \eqref{eq:new2}, it is dominated by a maximal function so
that
\[
\Big\|\rho^{-3}\il_{[|x_0-x_1|<\rho]} | f(x_0)| \, dx_0\Big\|_{L^3}
\less \|Mf\|_3 \less \|f\|_3 \less \|\nabla f\|_{\f32} \less \|D^2
f\|_1
\]
In passing, we remark that these arguments lead to a small potential
result which is analogous to  Proposition~\ref{prop:small_pot1}, see
Proposition~\ref{prop:small_pot2} below. If the potential is no
longer small, then one can only sum a finite Born series, and the
remainder is controlled by the exact same argument as above (which
does not distinguish between $\f{\sin(t\lambda)}{\lambda}$ and
$\cos(t\lambda)$ or $e^{it\lambda}$). The conclusion is that we have
proved the desired bound for all energies $\lambda\in
[\lambda_0,\infty)$ with $\lambda_0>0$ fixed.

\noindent As for the small energies, we remark that the changes
which need to be made to the analogous argument for
$\f{\sin(t\lambda)}{\lambda}$ are in the spirit of~\eqref{eq:new1}
and~\eqref{eq:new2}. We skip these details.

\noindent Finally, for the gradient estimate \eqref{eq:Dcos_dec} we
argue similarly to the $\sin$-case by passing the gradient through.
As explained previously, this leads to commutator terms that can be
treated by the $F(\lambda)$-method.
\end{proof}

For the sake of completeness, we record the small potential result
for the cosine.

\begin{proposition}
\label{prop:small_pot2} Assume that the real-valued potential $V$
satisfies $\|V\|_{\kato}<4\pi$ and $\|V\|_{L^{\f32}}<\infty$. Then
one has the bound
\[ \Big\| {\cos(t\sqrt{H})}f\Big\|_\infty \less
t^{-1} \|D^2 f\|_{L^{1}(\R^3)}
\]
for all $t>0$.
\end{proposition}

The following technical lemma was needed in the proof of
Proposition~\ref{prop:cos}.

\begin{lemma}
\label{lem:tedious} Let $0\le w(x)\le \la x\ra^{-4-\eps}$ for all
$x\in\R^3$. Then
\[ I(x,y;t):=\il_{[t=|x-x'|+|y-y'|]} \frac{w(x') w(y')}{|x-x'||y-y'|}\, dx'dy'
\less (\la t\ra^{-1} +\chi_{[|x|+|y|>t/4]})(\la x\ra\la y\ra)^{-1}
\]
for all $t\ge0$, $x,y\in\R^3$.
\end{lemma}
\begin{proof}
For the sake of simplicity, we assume that $\min(t,|x|,|y|)>1$. We
leave it to the reader to make the necessary modifications in case
this fails. First, we consider the contribution of $|x'|<\half|x|$
and $|y'|<\half|y|$ to the integral. It is easy to see that this is
bounded by
\begin{align*}
 &\less (\la x\ra\la y\ra)^{-1} \il_{[t=|x-x'|+|y-y'|]}
w(x')w(y')\chi_{[|x'|<\half|x|,\;|y'|<\half|y|]}\, dx'dy' \\
&\less (\la x\ra\la y\ra)^{-1}\;\chi_{[|x|+|y|>t/2]}
\end{align*}
Now let us suppose that $|x'|>\half|x|$ and $|y'|>\half|y|$. The
contribution from this regime is
\begin{align*}
 &\less (\la x\ra\la y\ra)^{-1}\la t\ra^{-1} \il_{[t=|x-x'|+|y-y'|]}
\frac{\la x'\ra w(x')}{|x-x'|} \chi_{[|x-x'|<t/2]}\;\la y'\ra w(y')\, dx'dy' \\
&= (\la x\ra\la y\ra)^{-1}\la t\ra^{-1} \il_{[|x-x'|<t/2]} \frac{\la
x'\ra w(x')}{|x-x'|} \il_{[|y-y'|=t-|x-x'|]} \la y'\ra w(y')\,
\sigma(dy')dx'\\
&\less (\la x\ra\la y\ra)^{-1}\la t\ra^{-1}
\end{align*}
Finally, we turn to the contribution of the regime $|x'|>\half|x|$
and $|y'|<\half|y|$. It is bounded by
\[ \la x\ra^{-1} \la y\ra^{-1}
\il_{[t=|x-x'|+|y-y'|]} \frac{\la x'\ra w(x')}{|x-x'|} \;\la y'\ra
w(y')\, dx'dy'
\]
If $|x-x'|>t/2$ in this integral, then we gain a factor of $\la
t\ra^{-1}$ as desired. So assume that $|x-x'|<t/2$. Then
$|y-y'|>t/2$, which in view of $|y'|<\half|y|$ forces that
$|y|>t/4$. This is again an allowed contribution, and the lemma
follows.
\end{proof}

\subsection{Stability in the potential}

In this final subsection, we comment on the dependence of the bounds
of this section on the potential.
The bounds of this subsection are helpful in obtaining the
contraction step in the nonlinear analysis of Section~\ref{sec:2}. Although that contraction
step can also be carried out by other means, we chose this path since it seems novel  and
of independent interest. For the sake of simplicity,
we will restrict ourselves to the operator
\[ H_a:= -\Delta + V_a, \qquad V_a(x) = -5\phi^4(x,a) \]
rather than trying to formulate this for general potentials.
We also set
\[ \calS_a(t):= \f{\sin(t\sqrt{H_a})}{\sqrt{H_a}} - c_0(\psi_a\otimes \psi_a) \]
where $\psi_a :=\pr_a \phi(\cdot,a)$.
Then we have the following bounds:

\begin{cor}
\label{cor:V_stab}
For any $\half<a,b<2$,
\begin{align}
\big\| [\calS_a(t)-\calS_b(t)]f\|_\infty &\less |a-b|\la t\ra^{-1} \|f\|_{W^{1,1}(\R^3)} \label{eq:diffS}\\
\big\|\nabla [\calS_a(t)-\calS_b(t)]f\|_\infty &\less |a-b|\la t\ra^{-1} \|f\|_{W^{2,1}(\R^3)} \nn\\
\big\|\big\{ [\cos(t\sqrt{H_a})P_{g(\cdot,a)}^\perp -  \cos(t\sqrt{H_b})P_{g(\cdot,b)}^\perp ]-K_t\big\}f\big\|_\infty
&\less |a-b| \la t\ra^{-1} \sum_{1\le |\alpha|\le 2} \|D^\alpha f\|_1 \nn\\
\big\|\big\{ \nabla[\cos(t\sqrt{H_a})P_{g(\cdot,a)}^\perp -  \cos(t\sqrt{H_b})P_{g(\cdot,b)}^\perp ]-\til K_t\big\}f\big\|_\infty
&\less |a-b| \la t\ra^{-1} \sum_{1\le |\alpha|\le 3} \|D^\alpha f\|_1\nn
\end{align}
where
\[ |K_t(x,y)| \less |a-b|(\la t\ra^{-1} +\chi_{[|x|+|y|>t/4]})(\la x\ra\la y\ra)^{-1}
\]
for all $t>0$, $x,y\in\R^3$, and similarly for $\til K_t$.
\end{cor}
\begin{proof}
We first remark that
\[ |V_a(x) - V_b(x)| \less \la x\ra^{-4} |a-b| \]
as well as
\[ |v_a(x) - v_b(x)|\less \la x\ra^{-2} |a-b| \]
where $-v_a^2 = V_a$.
Let us consider the difference of two typical Born-series terms
that arise in the analysis
of the sine evolution:
\begin{align} &\il_{\R^{3(k+2)}} \il_{-\infty}^\infty
\chi_0(\lambda) \sin(t\lambda) e^{i\lambda\sum_{j=0}^{k}
|x_j-x_{j+1}|}
 \frac{\prod_{j=1}^k V_a(x_j)}{\prod_{j=0}^k|x_j-x_{j+1}|}
 \,f(x_0)g(x_{k+1})\,dx_0\ldots dx_{k+1} \nn \\
 &- \il_{\R^{3(k+2)}} \il_{-\infty}^\infty
\chi_0(\lambda) \sin(t\lambda) e^{i\lambda\sum_{j=0}^{k}
|x_j-x_{j+1}|}
 \frac{\prod_{j=1}^k V_b(x_j)}{\prod_{j=0}^k|x_j-x_{j+1}|}
 \,f(x_0)g(x_{k+1})\,dx_0\ldots dx_{k+1} \label{eq:born_diff}
\end{align}
Since
\[ \prod_{j=1}^k V_a(x_j) - \prod_{j=1}^k V_b(x_j) =\sum_{\ell=1}^k \prod_{j=1}^{\ell-1} V_a(x_j)\;
(V_a(x_\ell)-V_b(x_\ell)) \prod_{j=\ell+1}^k V_a(x_j) \]
we can rewrite the difference in \eqref{eq:born_diff} as a sum of terms,
each of which contains the difference of $V_a$ and $V_b$. This
allows us to gain a factor of $|a-b|$ by the same arguments we used to derive
the dispersive estimate above.

The final term of the Born series, which we used to derive the high-energy
estimate of Lemma~\ref{lem:high} for $H_a$, involves the perturbed resolvent
$(-\Delta+V_a-\lambda+i0){-1}$. Hence, we now face the difference of the two
resolvents
\[ (-\Delta+V_a-\lambda+i0){-1} - (-\Delta+V_b-\lambda+i0){-1}
= (-\Delta+V_a-\lambda+i0){-1}[V_b-V_a ] (-\Delta+V_b-\lambda+i0){-1}
\]
Hence, the analysis which involves the bounds in~\eqref{eq:adec} can now be repeated,
and we gain a factor of $|a-b|$ here as well using the limiting absorption bounds~\eqref{eq:limap}.

As far as the low energies are concerned, we need to compare the (regular) parts of the resolvents
as in~\eqref{eq:3piece}. In other words, we need to compute the sine-transform of the difference
\[
 R_0(\lambda) v_a \Ea(\lambda) v_a R_0(\lambda) -  R_0(\lambda) v_b \Eb(\lambda) v_b R_0(\lambda)
\]
where $\Ea, \Eb$ are the operators arising in~\eqref{eq:Edef}.
The only really novel term here is
\[
 R_0(\lambda) v_a (\Ea(\lambda)-\Eb(\lambda)) v_a R_0(\lambda)
\]
To understand the difference $\Ea(\lambda)-\Eb(\lambda)$, we start with the observation that
\[ A_a(\lambda) - A_b(\lambda) = v_aR_0(\lambda) v_a - v_b R_0(\lambda) v_b \]
which implies that
\begin{align*}
 A_a(\lambda)^{-1} - A_b(\lambda)^{-1} &= - A_a(\lambda)^{-1}[ v_aR_0(\lambda) v_a - v_b R_0(\lambda) v_b ]
A_b(\lambda)^{-1} \\
&= - A_a(\lambda)^{-1}[(v_a-v_b)R_0(\lambda) v_a - v_b R_0(\lambda) (v_b-v_a) ]
A_b(\lambda)^{-1}
\end{align*}
We already now that $A_a(\lambda)^{-1}$ and $A_b(\lambda)^{-1}$ satisfy the bounds~\eqref{eq:Fest1}
and~\eqref{eq:Fest2}, i.e., they are amenable to the $F(\lambda)$ method. However, by~\eqref{eq:new_piece}
the middle piece $(v_a-v_b)R_0(\lambda) v_a$ (and its symmetric counterpart) satisfies the same bounds.
Moreover, these bounds come with a factor of $|a-b|$ which is the desired gain.

The other constituents of $\Ea$ (and of $\Eb$) are one-dimensional, and they involve the following basic
building blocks:
\[ \Sa_1, \; \ma(0),\; \Ea_1(\lambda),\; \Ea_2(\lambda) \]
and the latter two can be further reduced to the pieces
\[ \Aa_1(\lambda), \; \ma_1(\lambda) \]
in view of the Neumann series \eqref{eq:Ainv} and~\eqref{eq:m0inv}.
Thus, in order to show that the low-frequency dispersive estimate gains a factor of $|a-b|$,
it will suffice to show that the difference of any two of these pieces with parameter values $a$ and~$b$
satisfies the bounds~\eqref{eq:Fest1} and~\eqref{eq:Fest2} with a gain of $|a-b|$.
By \eqref{eq:S1def},
\[ \Sa_1 = \|v_a\psi_a\|_2^2\; w_a\psi_a \otimes w_a\psi_a \]
so that
\[ \|\Sa_1-\Sb_1 \|_{2\to2} \less |a-b| \]
The operators $\ma(0)-\mb(0)$ are equally easy to deal with.
We leave it to the reader to verify the $F(\lambda)$ property for  the pair-wise (i.e., relative to
$a$ and $b$)
differences of each the
operators $\Aa_1(\lambda), \; \ma_1(\lambda)$, gaining a factor $|a-b|$ in the process.
As for the former, the logic is that
\beeq\label{eq:Aadiff} \Aa_1(\lambda) - \Ab_1(\lambda) = \frac{1}{\lambda}
(v_a-v_b)(x)\frac{e^{i\lambda|x-y|}-1}{4\pi|x-y|}v_a(y) -
\frac{1}{\lambda} v_b(x)\frac{e^{i\lambda|x-y|}-1}{4\pi|x-y|}(v_b(y)-v_a(y))
\eneq
and the two constituents on the right-hand side are again of the same form as $A_1(\lambda)$.
Since we have established the $F(\lambda)$ property for this operator in the proof of Lemma~\ref{lem:low} above,
it follows that the same holds for the difference~\eqref{eq:Aadiff}. Moreover, we gain a factor $|a-b|$ in
the constants. The same type of logic applies to $m_1(\lambda)$, and we leave the details to the reader.
To conclude the proof of~\eqref{eq:diffS}, it only remains to control the difference of  $\calS_0-c_0(\psi\otimes\psi)$
for $a$ and $b$,
see~\eqref{eq:VV}. However, this is an explicit multi-linear expression in the potentials, and can be estimated
just like above.

\noindent The remaining bounds stated in this corollary can be proved by the basically identical arguments,
which we skip.
\end{proof}

\end{document}